\documentclass[11pt,a4paper]{article}
\usepackage{amssymb}
\usepackage{amsthm}
\usepackage{amsmath}
\usepackage{amsbsy}
\usepackage[utf8]{inputenc}
\usepackage[a4paper, left=2cm, right=2cm, top=3cm, bottom=3cm]{geometry}
\usepackage[T1]{fontenc}
\usepackage{lmodern}
\usepackage{extarrows}
\usepackage{graphicx}
\usepackage{caption}
\usepackage{subcaption}
\usepackage{float}
\usepackage[linktocpage,colorlinks=true]{hyperref}
\usepackage{thmtools,thm-restate}

\hypersetup{urlcolor=blue, citecolor=red, linkcolor=blue}

\newtheorem{theorem}{Theorem}[section]
\newtheorem{definition}[theorem]{Definition}
\newtheorem{lemma}[theorem]{Lemma}

\newtheorem{rem}[theorem]{Remark}
\theoremstyle{definition}

\newcommand{\R}{\mathbb{R}}
\newcommand{\N}{\mathbb{N}}
\newcommand{\DS}{\displaystyle}

\begin{document}
\begin{center}
{\Large \bf Computer-assisted proof of heteroclinic connections \\[1ex]
in the one-dimensional Ohta-Kawasaki model}

 \vskip 0.5cm
{\large
Jacek~Cyranka$^{*,\ddag,\dag}$, Thomas~Wanner$^{\S}$
}
\vskip 0.5cm
{\small$^\ddag$ Institute of Computer Science and Computational Mathematics,
Jagiellonian University}\\
{\small  ul. S. {\L}ojasiewicza 6, 30-348 Krak\'ow, Poland}
\vskip 0.5cm
{\small$^\dag$ Institute of Applied Mathematics and Mechanics,
University of Warsaw}\\
{\small  Banacha 2, 02-097 Warszawa, Poland}
\vskip 0.5cm
{\small$^*$ Department of Mathematics, Rutgers, The State University of New Jersey, 110 Frelinghusen Rd, Piscataway, NJ  08854-8019, USA}
\vskip 0.5cm
{\small$^{\S}$ Department of Mathematical Sciences, George Mason University,
Fairfax, VA 22030, USA}
\vskip 0.5cm
{jcyranka@gmail.com, jacek.cyranka@ii.uj.edu.pl, twanner@gmu.edu}
\vskip 0.5cm

\today
\end{center}

\providecommand{\matrixS}{{T}}
\providecommand{\seriesD}[2][l]{d_{#2}^{#1}}
\providecommand{\seriesA}[2][l]{a_{#2}^{#1}}
\providecommand{\seriesB}[2][l]{b_{#2}^{#1}}
\providecommand{\seriesAconst}{2^{2m}}
\providecommand{\operS}{{L_{\lambda}^l}}
\providecommand{\operSS}[2][\lambda]{{L_{#1}^{#2}}}
\providecommand{\Sp}{{Sp}}
\providecommand{\re}{{Re}}
{\bf Abstract.}
We present a computer-assisted proof of heteroclinic connections in
the one-dimensional Ohta-Kawasaki model of diblock copolymers. The
model is a fourth-order parabolic partial differential equation subject
to homogeneous Neumann boundary conditions, which contains as a special
case the celebrated Cahn-Hilliard equation. While the attractor structure
of the latter model is completely understood for one-dimensional domains,
the diblock copolymer extension exhibits considerably richer long-term
dynamical behavior, which includes a high level of multistability.
In this paper, we establish the existence of certain heteroclinic connections between
the homogeneous equilibrium state and local and global energy minimizers.

The proof of the above statement is conceptually simple, and combines several
techniques from some of the authors' and Zgliczy\'nski's works. Central for the
verification is the rigorous propagation of a piece of the unstable manifold of
the homogeneous state with respect to time. This propagation has to lead to
small interval bounds, while at the same time entering the basin of attraction
of the stable fixed point. For interesting parameter values the global attractor
exhibits a complicated equilibrium structure, and the dynamical equation is rather
stiff. This leads to a time-consuming numerical propagation of error bounds, with
many integration steps. This problem is addressed using an efficient algorithm for
the rigorous integration of partial differential equations forward in time. The
method is able to handle large integration times within a reasonable computational
time frame, and this makes it possible to establish heteroclinic at various
nontrivial parameter values.


\section{Introduction}
The goal of this paper is to propose a computer assisted method of constructively proving 
existence of heteroclinic connections between stationary solutions of
\emph{dissipative partial differential equations (PDEs)}. As a case study we 
establish a computer assisted proof of heteroclinic
connections in the one-dimensional \emph{Ohta-Kawasaki diblock copolymer model}.
Precisely, we establish the existence of certain heteroclinic connections between
the homogeneous equilibrium state, which represents a perfect copolymer mixture,
and all local and global energy minimizers. In this way, we show that not every solution
originating near the homogeneous state will converge to the global energy minimizer,
but rather is trapped by a stable state with higher energy. This phenomenon can
not be observed in the one-dimensional \emph{Cahn-Hillard} equation, where generic solutions 
are attracted by a global minimizer.

Our method is general. To achieve the presented goal we develop a framework combining efficient algorithms
and implementation in the C++ programming language performed by the first author~\cite{Cy} and
the theory of \emph{cone conditions} and \emph{self-consistent bounds} developed by Zgliczy\'nski et al.~\cite{BW2, CZ, CzZ, Zattr, ZPer, ZM, Z3}.
We remark that the software is not
hard coded for the particular equation studied in this paper, but is rather a generic tool that can
be easily adapted for study of any dissipative PDE on a one dimensional domain with a polynomial nonlinearity
and periodic/Dirichlet/Neumann boundary conditions.

To the best of our knowledge, this is the first result about
proving constructively the existence of heteroclinic connections for dissipative PDEs.
In addition, it also rigorously shows for the first time that
there are infinitely many orbits which
are trapped by each of the two local and of the two suspected global
energy minimizers. In fact, the proof of
our result below shows that not all solutions
of one dimensional Ohta-Kawasaki diblock copolymer model,
which originate close to the homogeneous state, and satisfty certain space-translational symmetry,
converge to the global energy minimizer.

Our approach differs in spirit from techniques based on Newton-Kantorovich
type arguments in a functional analytic setting, which were used
for example in~\cite{CGL, GLP, BBMJLM}.
We would like to mention that establishing the existence
of connecting orbits in finite-dimensional continuous and discrete
dynamical systems has been the subject of a number of studies, see
for example~\cite{MJM, BBMJLM, Wilczak2006, Wilczak2005, WZ}. In addition,
the global dynamics of the Swift-Hohenberg model has been uncovered
using a Morse decomposition based approach in~\cite{DHMO}, and the
existence of periodic orbits for the ill-posed Boussinesq equation
was established in~\cite{CGL, CzZ}. Finally, the existence of a
heteroclinic orbit between the trivial and a nontrivial stationary
state has been established for the one-dimensional $p$-Laplace equation
in~\cite{CM}. Also, private communication revealed that
Zgliczy\'nski could verify the existence of a heteroclinic connection
in the Kuramoto-Sivashinsky equation using similar techniques in
the recent work in progress~\cite{Zhet}.

In the following two subsections, we first
briefly introduce the Ohta-Kawasaki diblock copolymer model,
structure of its dynamics, relevant bifurcation diagrams, and present our main result (Theorem~\ref{thmmainintro}). After that,
we present the computer assisted framework used to prove the main result.
\subsection{The diblock copolymer model description and the main result}
Phase separation phenomena in materials sciences provide a rich source
for pattern formation mechanisms. In many situations, the resulting
models are dissipative parabolic partial differential equations
which are amenable to rigorous mathematical treatment, and can therefore
explain the underlying pattern formation principles. In the current
paper, we consider one of these models which has been of significant
interest in recent years. This equation is of phase field type, and 
it models microphase separation in diblock copolymers. From a mathematical
perspective, the model has been described in~\cite{nishiura:ohnishi:95a},
in its original form it was proposed by Ohta and Kawasaki~\cite{ohta:kawasaki:86a},
as well as Bahiana and Oono~\cite{bahiana:oono:90a}. See also the discussion
in~\cite{choksi:peletier:09a, choksi:ren:03a}. Consider a material
constrained to the bounded domain~$\Omega \subset \R^d$. Then the model
is based on an energy functional which is comprised of the standard
van der Waals free energy with an additional nonlocal term, given by
\begin{equation} \label{dbcp:energy}
  E_{\epsilon,\sigma}[u] =
  \int_{\Omega} \left( \frac{\epsilon^2}{2} |\nabla u|^2 +
    F(u) \right) dx +
  \frac{\sigma}{2} \int_{\Omega}
    \left|(-\Delta )^{-1/2} (u(x) - \mu)\right|^2 dx \; .
\end{equation}
With this energy one can associate
gradient-like dynamics, which is then used to model the evolution
of the material as a function of time.
The diblock copolymer model in its standard
form considers the evolution equation
\begin{eqnarray}
  u_t & = & -\Delta \left( \epsilon^2 \Delta u + f(u) \right) -
    \sigma(u-\mu)
    \quad\mbox{ in }\; \Omega
    \; , \label{dbcp} \\[1ex]
  & & \mu = \frac{1}{|\Omega|} \int_{\Omega} \! u(x) \, dx \; ,
    \quad\mbox{ and }\quad
    \frac{\partial u}{\partial \nu} =
    \frac{\partial \Delta u}{\partial \nu} = 0
  \quad\mbox{ on }\; \partial \Omega \; ,
  \nonumber
\end{eqnarray}
where the nonlinearity is defined as $f(u) = -F'(u)$, and which
in our situation is given by $f(u) = u - u^3$. This evolution 
equation is based on the gradient in the $H^{-1}$-topology and uses 
homogeneous Neumann boundary conditions for both~$u$ and~$\Delta u$.
The diblock copolymer model~(\ref{dbcp}) is an extension
of the celebrated Cahn-Hilliard equation~\cite{cahn:hilliard:58a},
which corresponds to the special case $\sigma = 0$ and serves as
a fundamental model for the phase separation phenomena spinodal
decomposition~\cite{maier:wanner:98a, maier:wanner:00a,
sander:wanner:99a, sander:wanner:00a, wanner:04a}
and nucleation~\cite{bates:fife:93a, bloemker:etal:10a,desi:etal:11a}.
The quantity $\mu$ average of the solution over the domain is conserved
in time.

\begin{figure}
  \centering
  \setlength{\unitlength}{1 cm}
  \begin{picture}(14.5,12.5)
    \put(0.0,6.5){%
      \includegraphics[width=7.0cm]{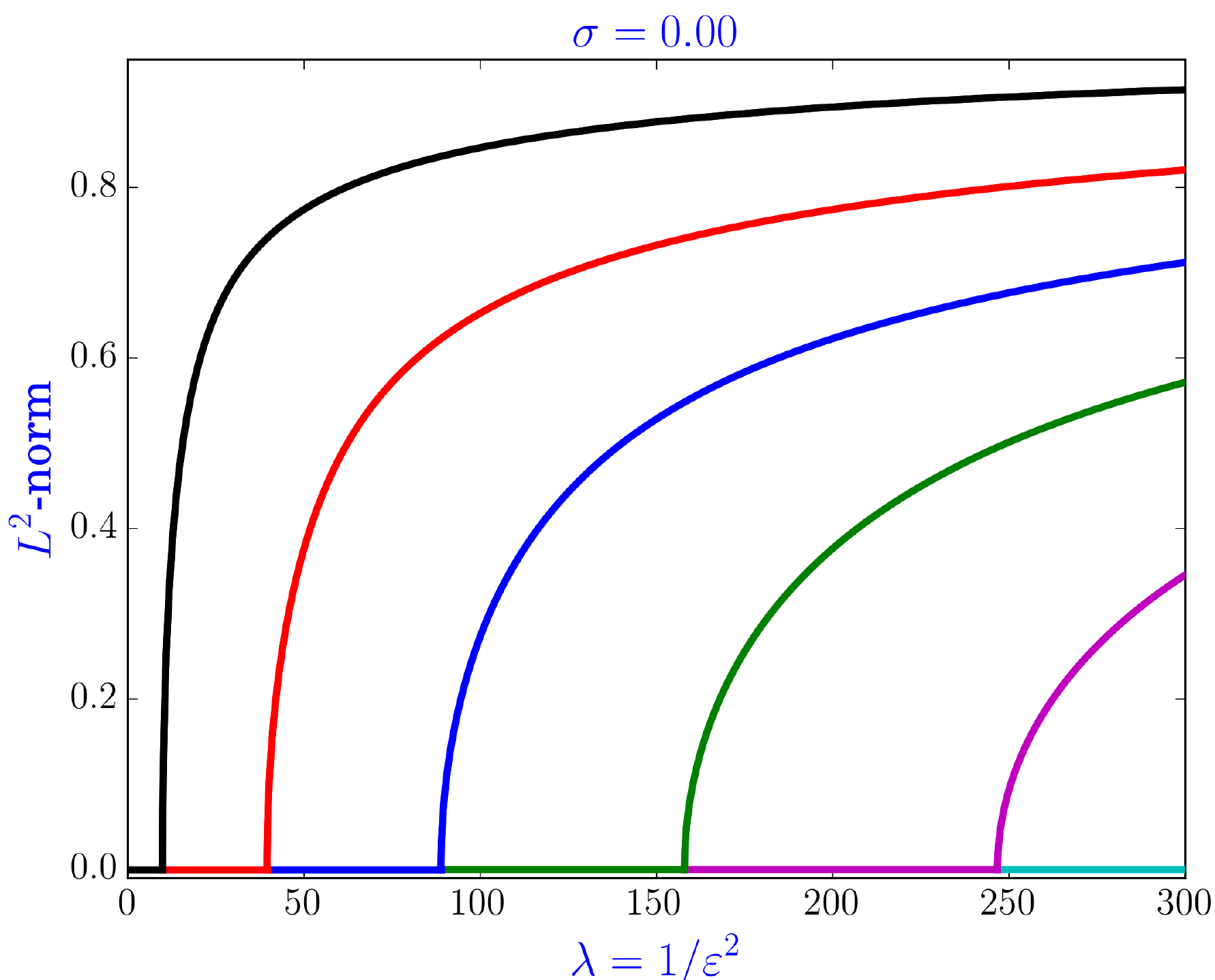}}
    \put(7.5,6.5){%
      \includegraphics[width=7.0cm]{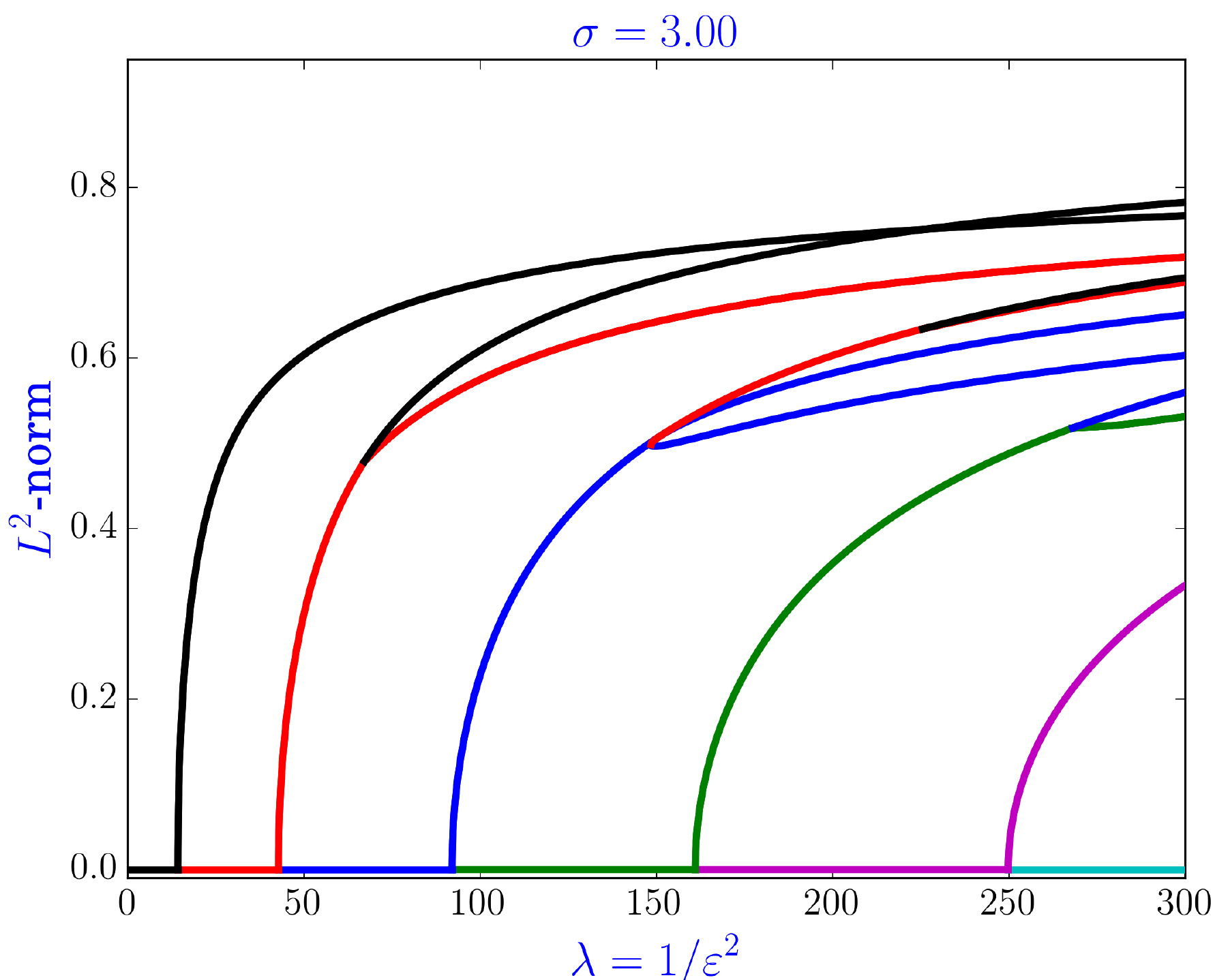}}
    \put(0.0,0.0){%
      \includegraphics[width=7.0cm]{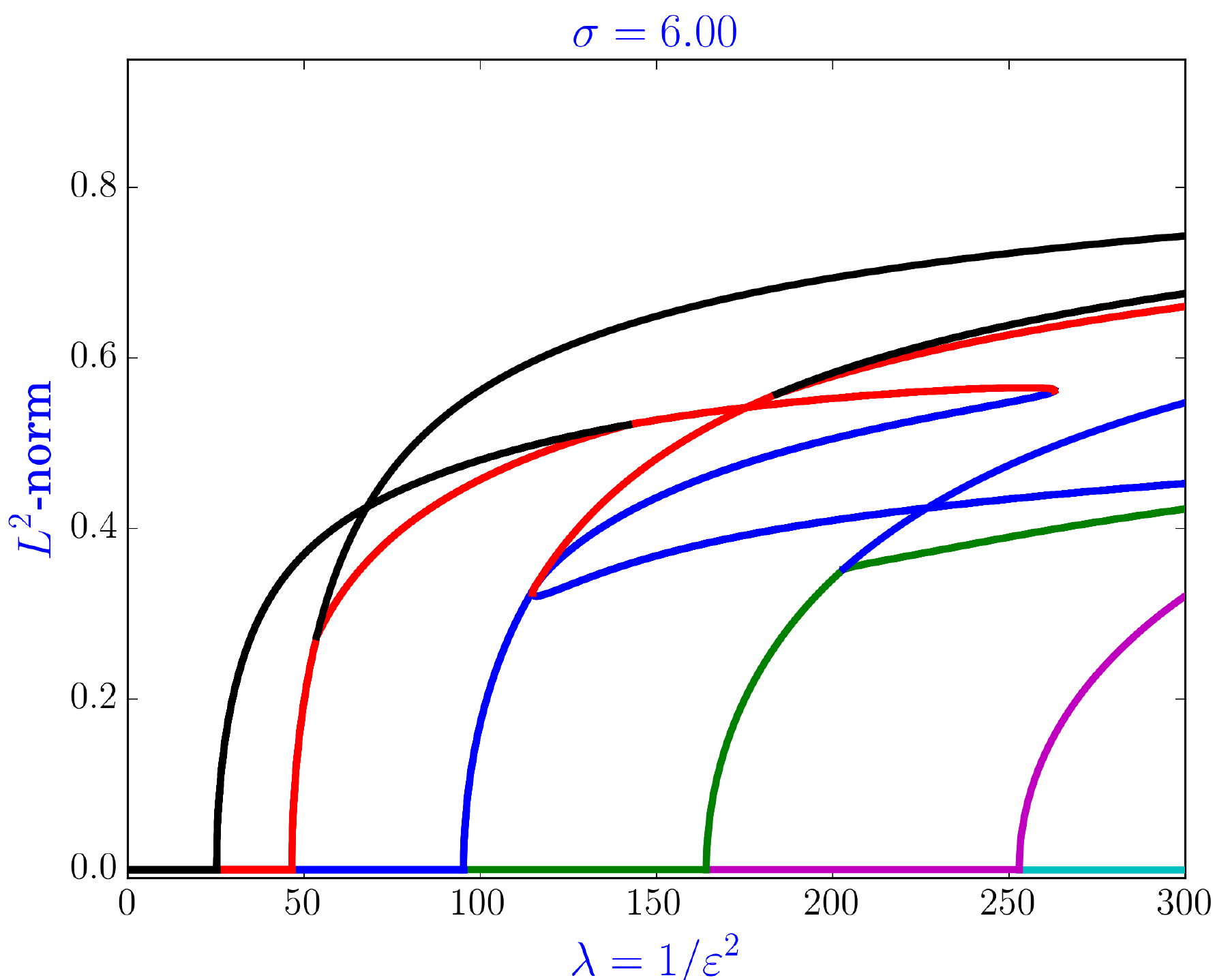}}
    \put(7.5,0.0){%
      \includegraphics[width=7.0cm]{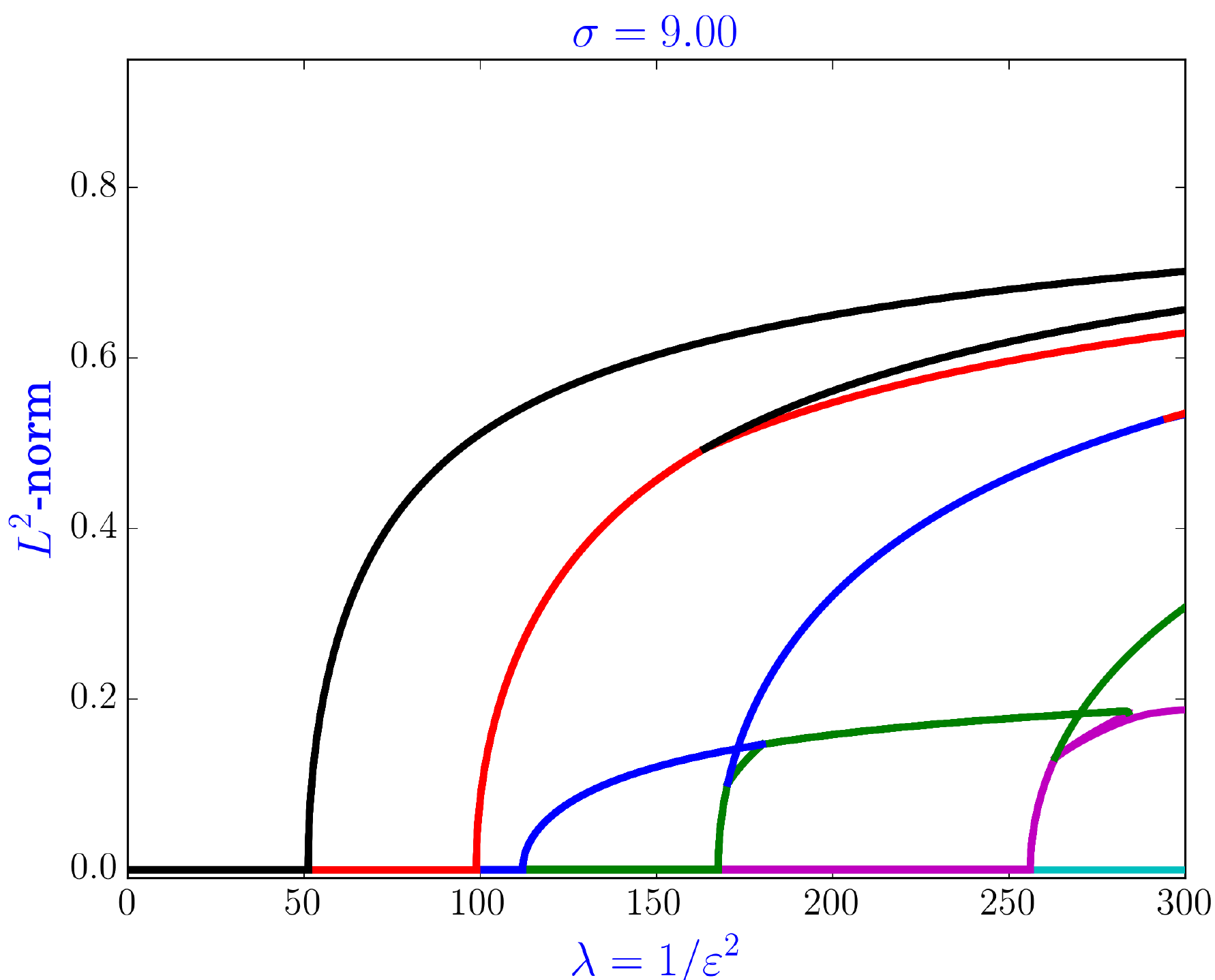}}
  \end{picture}
  \caption{Bifurcation diagrams for the diblock
           copolymer model~(\ref{dbcp}) on the domain~$\Omega
           = (0,1)$ and for total mass $\mu = 0$. From top
           left to bottom right the bifurcation diagrams are
           for $\sigma = 0$, $3$, $6$, and~$9$, respectively.
           In each diagram, the vertical axis measures the
           $L^2(0,1)$-norm of the solutions, and the horizontal
           axis uses the parameter $\lambda = 1 / \epsilon^2$.
           Each colored point corresponds to at least one
           stationary solution of the diblock copolymer model,
           and the color indicates the Morse index of the
           equilibrium in the evolution equation~(\ref{dbcp}),
           as described in the figure caption.
           The solution branches are color-coded by the Morse
           index of the solutions, and black, red, blue, green,
           magenta, and cyan correspond to indices~$0$, $1$,
           $2$, $3$, $4$, and~$5$, respectively.}
  \label{figbifdiags}
\end{figure}

For the present paper, we focus exclusively on the diblock
copolymer model on one-dimensional domains~$\Omega$. Recent
studies in~\cite{johnson:etal:13a} have demonstrated that the
structure of its equilibrium set is considerably richer than
the one for the Cahn-Hilliard equation described
in~\cite{grinfeld:novickcohen:95a}. This is visualized in
Figure~\ref{figbifdiags},
which contains equilibrium bifurcation diagrams for the parameter
values $\sigma = 0$, $3$, $6$, and~$9$, from top left to bottom
right, and for total mass $\mu = 0$ and domain $\Omega = (0,1)$.
Each diagram indicates the structure of the equilibrium set as a
function of the new parameter~$\lambda = 1 / \epsilon^2$.
The interesting regime is the limit
$\epsilon \to 0$, and by switching to~$\lambda$ the diagrams
become easier to visualize. The horizontal straight lines in all images
represent the constant function~$u \equiv \mu = 0$, which is
clearly an equilibrium, and which we call the {\em trivial
solution\/}. As~$\lambda$ increases, nontrivial solutions
bifurcate from the trivial branch in pairs, even though we
only plot one of these branches. This is due to the fact that
the vertical axes of the bifurcation diagrams show the
$L^2(\Omega)$-norm of the stationary state, and the bifurcation
pairs are related by a norm-preserving symmetry.

For the Cahn-Hilliard case $\sigma = 0$, i.e., in the top left
bifurcation diagram, pairs of solutions are created at the bifurcation
points~$\lambda = k^2 \pi^2$ for $k \in \N$, and the solutions
on the $k$-th branch look qualitatively similar to the function
$\cos(k \pi x)$. In this case, each nontrivial solution branch
corresponds to exactly two nontrivial solutions. It was shown
in~\cite{grinfeld:novickcohen:95a} that there are no other 
equilibrium solutions for the equation, i.e., none of the 
bifurcation branches exhibit secondary bifurcations. Moreover,
only the solutions on the first branch are stable equilibria,
and they correspond to solutions with one transition layer.
The situation changes dramatically for $\sigma > 0$. In the
diblock copolymer case, all branches which are created at the
trivial solution exhibit secondary bifurcations. Even more
importantly, these secondary bifurcations can lead to the creation
of multiple stable equilibrium solutions. While as of yet there
are no classical mathematical proofs for these statements, 
computer-assisted proof methods were developed in~\cite{wanner:17a},
which allow for the validation of these equilibrium branches.
See also the work in~\cite{sander:wanner:16a, wanner:17b}.
Some first results for two-dimensional domains can be found
in the recent paper~\cite{BW}.

But how does this multistability manifest itself in the evolution
equation~(\ref{dbcp})? Based on the underlying physical situation,
one is usually interested in studying the long-term behavior of
solutions of the diblock copolymer model which originate close to
the (unstable) homogeneous equilibrium state~$u \equiv \mu$. These
long-term limits are usually periodic in space with a specific period which
depends on the precise parameter combinations. Since the system
exhibits multiple stable states which all are at least local
minimizers of the energy, it is natural to wonder whether typical
solutions converge to the global energy minimizers or are trapped
earlier. A systematic numerical study of this long-term behavior
was performed in~\cite{johnson:etal:13a}, and it was able to show
that for fixed parameters~$\lambda > 0$ and~$\sigma > 0$, most
solutions starting near the homogeneous state usually converge
to the same long-term limit. Moreover, the numerics indicated that
this preferred long-term limit is usually not the stable equilibrium
with the lowest energy. For more details, we refer the reader
to~\cite{johnson:etal:13a}. In addition, the recent
paper~\cite{wanner:16a} contains a heuristic explanation of
these numerical observations. We would like to point out, however,
that to the best of our knowledge there is no rigorous complete
description of the equilibrium set of the diblock copolymer model.
Thus, even though the bifurcation diagrams in Figure~\ref{figbifdiags}
for $\sigma > 0$ show multiple stable equilibria which can be ordered
according to their energies, the lowest energy steady state is
not known to be the global energy minimizer --- there might be 
more equilibria than the ones shown in the diagram. However, since
it is generally assumed that these diagrams are complete in the 
shown parameter ranges, we will refer in the following to the
stable steady state in a diagram with the lowest energy as the
{\em suspected global minimizer}.

Describing the complete attractor in a mathematically rigorous
manner has proven to be elusive, even for the Cahn-Hilliard
special case $\sigma = 0$. While the equilibrium structure
of the Cahn-Hilliard equation is completely known for one-dimensional
base domains~\cite{grinfeld:novickcohen:95a}, the precise structure
of heteroclinic connections is an open question and has been 
answered only partially~\cite{bloemker:etal:10a,
grinfeld:novickcohen:99a}. Even for special
two-dimensional domains such as the unit square, the structure
of the Cahn-Hilliard attractor is extremely complicated, see for
example the discussions in~\cite{fife:kielhofer:etal:97a,
maier:etal:08a, maier:etal:07a, maier:wanner:97a}.
In the case of the full diblock copolymer
model, i.e., for $\sigma > 0$, even less is known. There are mathematical
results on specific types of equilibrium solutions for two- and
three-dimensional domains~$\Omega$, see for example~\cite{ren:wei:06b,
ren:wei:06a, ren:wei:07b, ren:wei:07a} as well as the references
therein. Moreover, numerical studies have been performed of the
long-term behavior of solutions of~(\ref{dbcp}) which originate
close to the homogeneous state~\cite{choksi:etal:11a,
choksi:peletier:09a}. 

\begin{figure}[bt]
  \centering
  \setlength{\unitlength}{1 cm}
  \begin{picture}(15.0,6.0)
    \put(0.0,0.0){%
      \includegraphics[width=7.0cm]{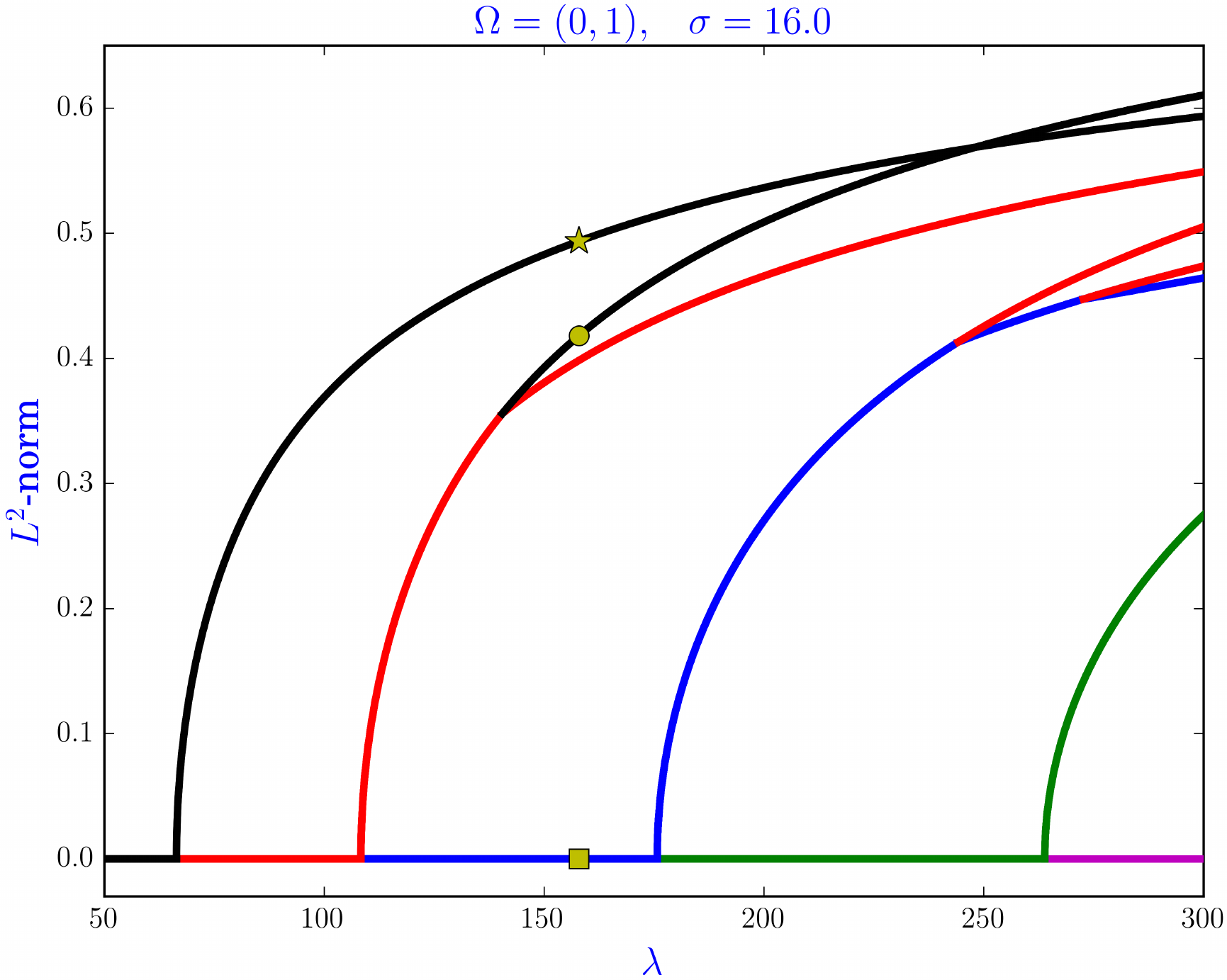}}
    \put(8.0,0.0){%
      \includegraphics[width=7.0cm]{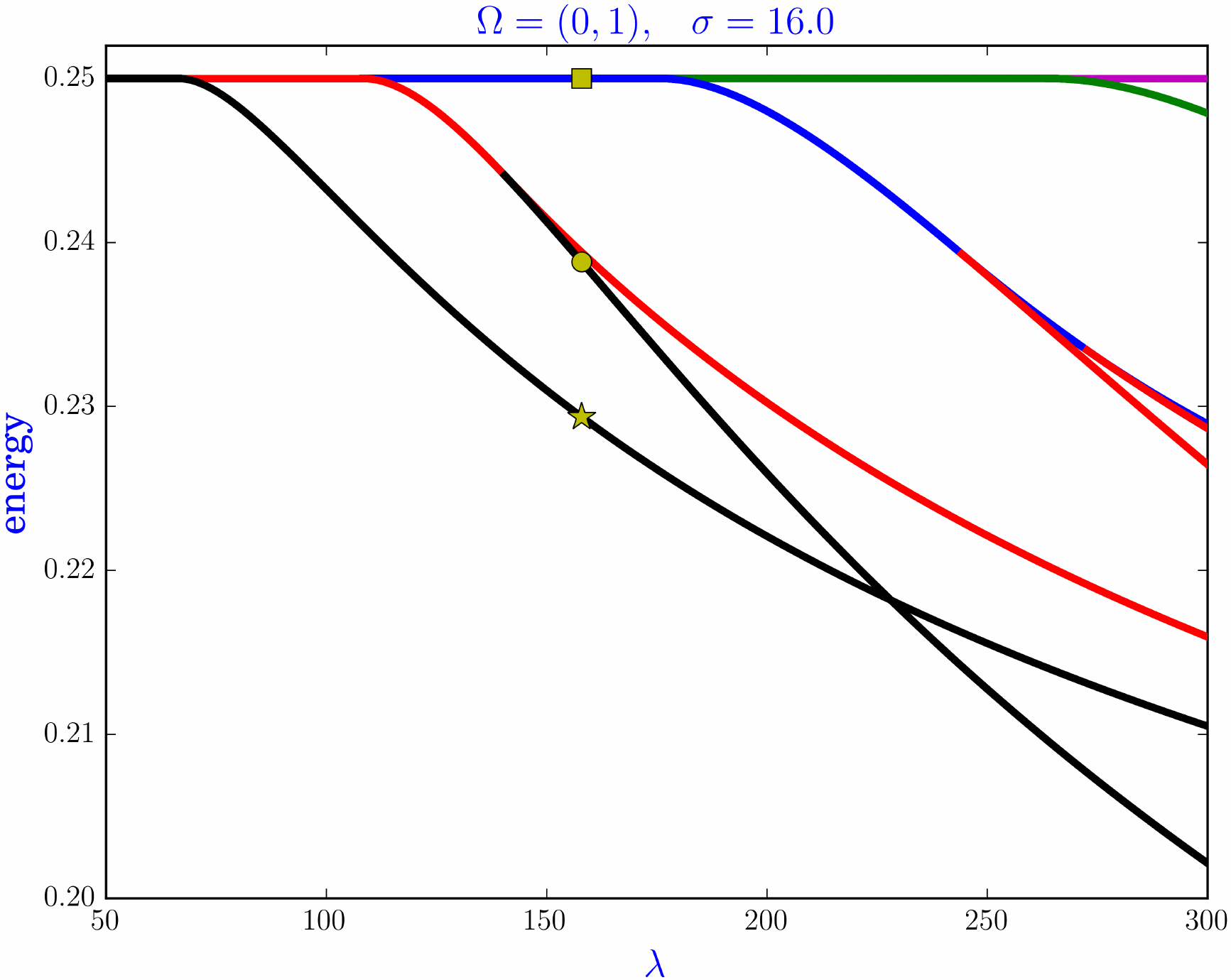}}
  \end{picture}
  \caption{Bifurcation diagrams for the diblock
           copolymer model~(\ref{dbcp}) on the domain~$\Omega
           = (0,1)$, for total mass $\mu = 0$, and for nonlocal
           interaction parameter $\sigma = 16$. In the left
           diagram, the vertical axis measures the $L^2(0,1)$-norm
           of the solutions, while the right diagram shows the
           total energy~(\ref{dbcp:energy}). In both diagrams,
           the horizontal axis uses the parameter $\lambda =
           1 / \epsilon^2$. As in Figure~\ref{figbifdiags}, the
           solution branches are color-coded by the Morse index
           of the solutions, and black, red, blue, green, and
           magenta correspond to indices~$0$, $1$, $2$, $3$,
           and~$4$, respectively.}
  \label{figbifdiagsONE}
\end{figure}

Despite the above numerical evidence, there
has been no mathematical result which shows that large numbers of
solutions of the diblock copolymer model which start near the trivial
equilibrium solution are in fact trapped by a local energy minimizer,
and that at the same time energy minimizers with lower energy can be
reached from the homogeneous state as well.
It is the goal of the present paper to close this gap for the bifurcation
diagram shown in Figure~\ref{figbifdiagsONE}.

This will be accomplished
using a computer-assisted proof technique, and it will establish the
existence of heteroclinic solutions between the homogeneous state
and the local minimizers, as well as between the homogeneous state
and local minimizers with lower energy, which are in fact suspected to
be the global minimizers. More precisely, we verify heteroclinic
connections from the state indicated by a yellow square, to both
the equilibria indicated by yellow circles and the equilibria indicated
by a yellow star. Our main result is the following theorem

\begin{figure}[tb]
  \centering
  \setlength{\unitlength}{1 cm}
  \begin{picture}(14.5,12.5)
    \put(0.0,6.5){%
      \includegraphics[width=7.0cm]{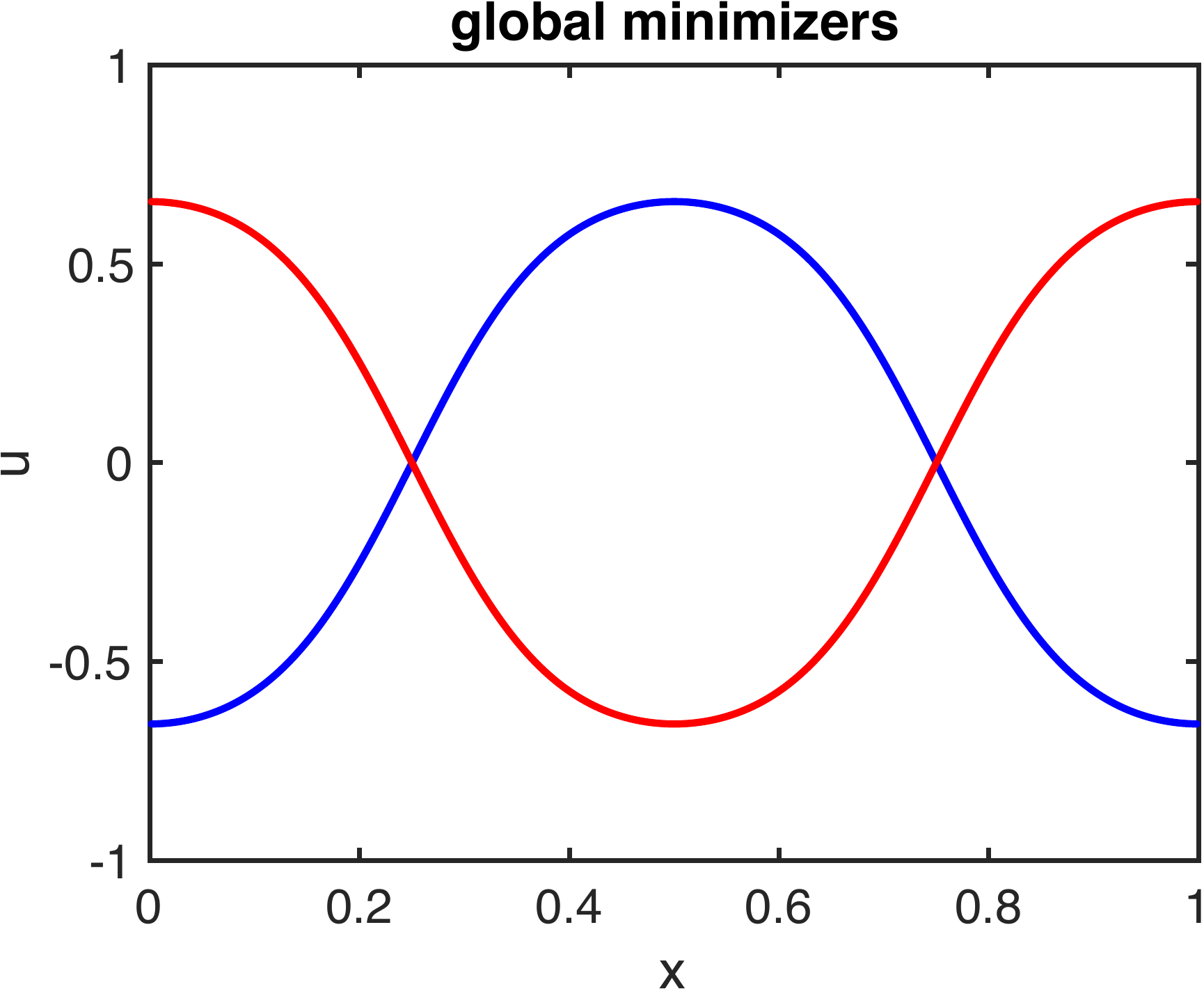}}
    \put(7.5,6.5){%
      \includegraphics[width=7.0cm]{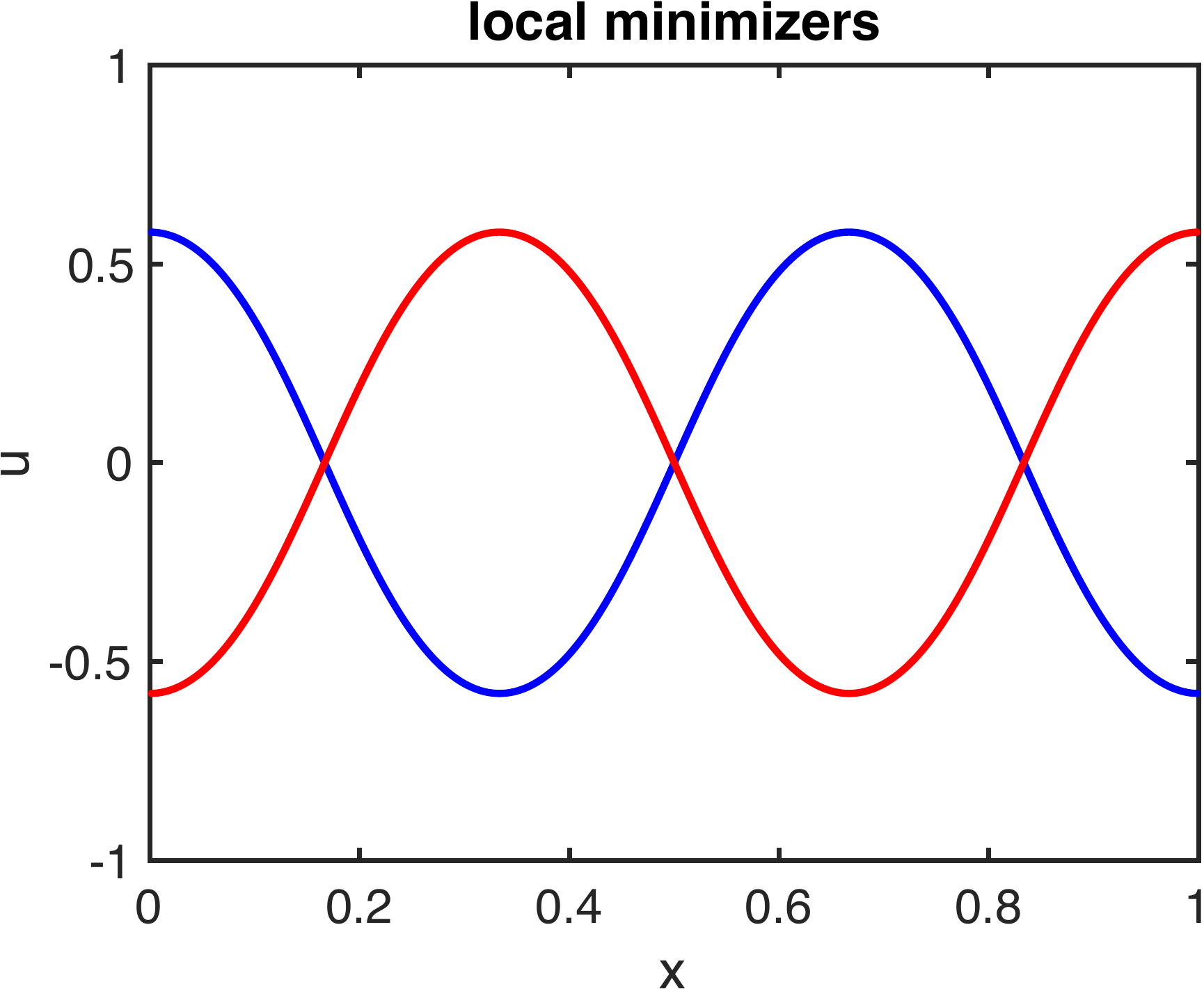}}
    \put(0.0,0.0){%
      \includegraphics[width=7.0cm]{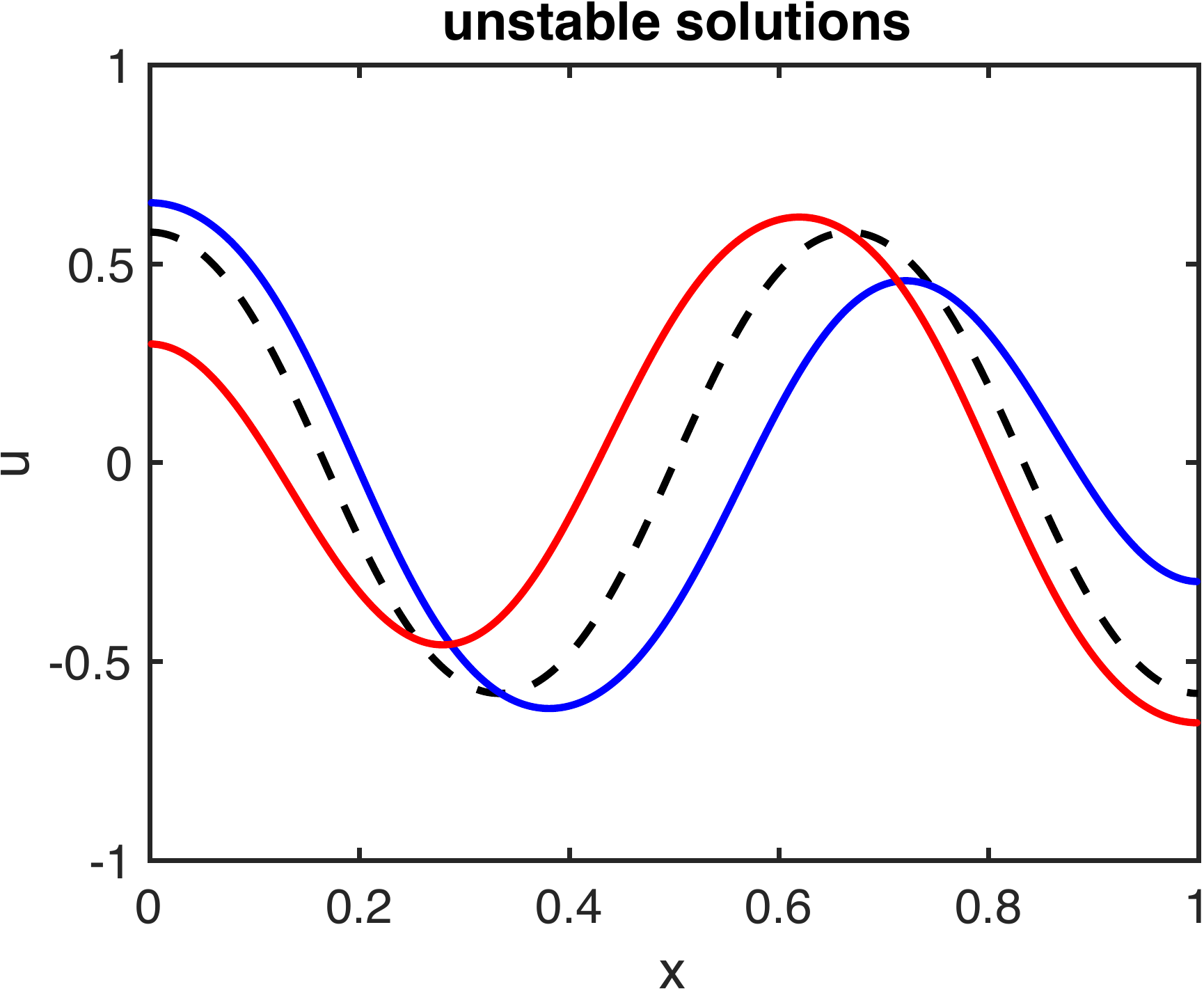}}
    \put(7.5,0.0){%
      \includegraphics[width=7.0cm]{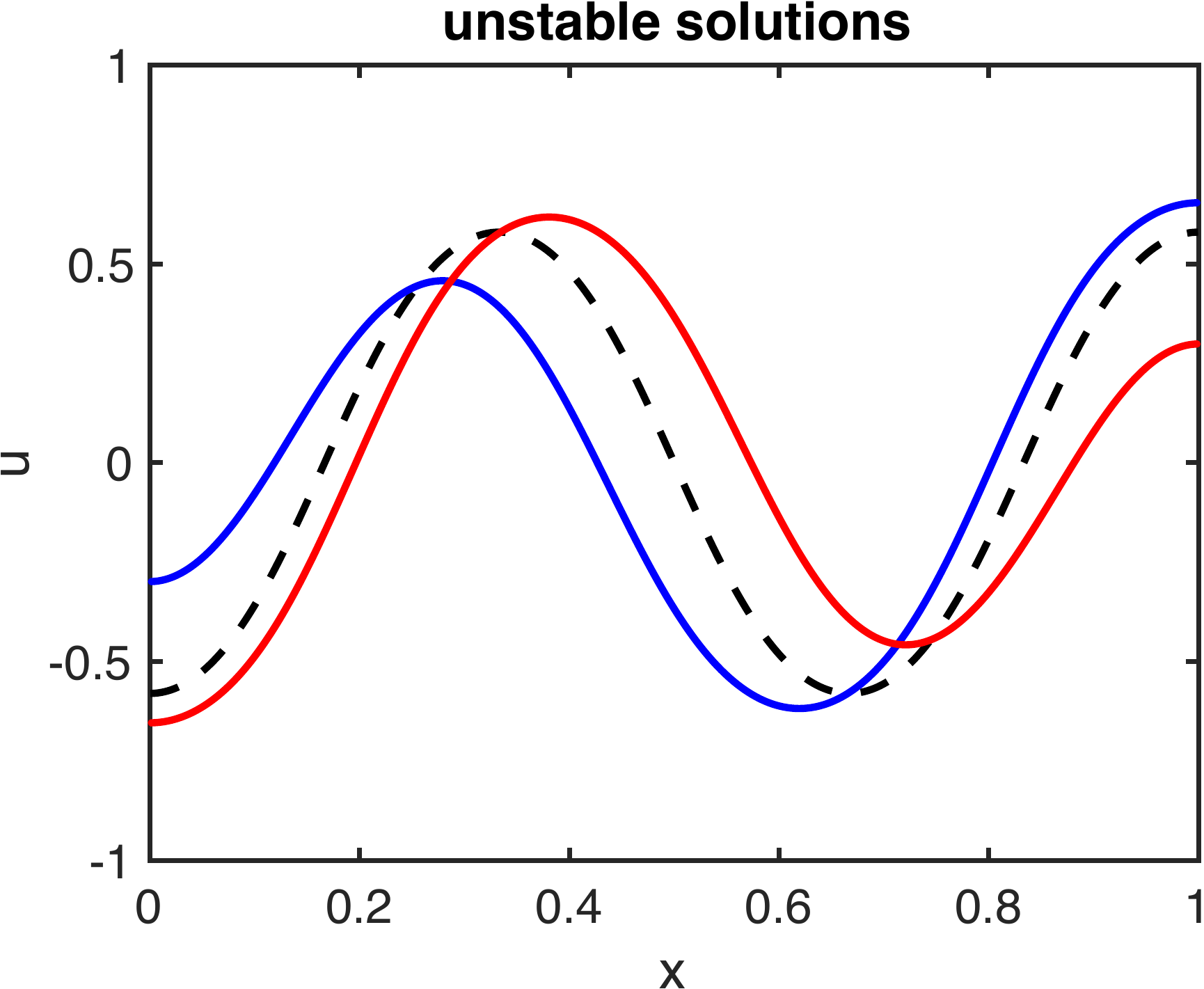}}
  \end{picture}
  \caption{Stationary solutions of the diblock copolymer
           model~(\ref{dbcp}) on the domain~$\Omega = (0,1)$,
           for total mass $\mu = 0$, and for interaction length
           parameters $\lambda = 16 \pi^2 \approx 157.91$ and
           $\sigma = 16$. The top left image contains the two
           suspected global energy minimizers, which are indicated
           by a star in the bifurcation diagrams of
           Figure~\ref{figbifdiagsONE}. The image in the top
           right shows the local energy minimizers, which
           are marked by circles in Figure~\ref{figbifdiagsONE}.
           The two bottom panels show the index one solutions
           which are created through a stabilizing pitchfork
           bifurcation. In each figure, the local minimizer 
           from the respective primary bifurcation branch is
           shown as a dashed black curve.}
  \label{figspecificsols}
\end{figure}

\begin{theorem}[Existence of heteroclinic connections]
\label{thmmainintro}
Consider the diblock copolymer equation~(\ref{dbcp}) on the
one-dimensional domain~$\Omega = (0,1)$, for interaction lengths
$\lambda = 1 / \epsilon^2 = 16\pi^2$ and $\sigma = 16$, and
for total mass $\mu = 0$. Then there exist heteroclinic
connections between the unstable homogeneous stationary state
$u \equiv \mu$ and each of the two local energy minimizers which
are indicated by yellow circles in Figure~\ref{figbifdiagsONE},
and shown in the top right panel of Figure~\ref{figspecificsols}.
In addition, there exist heteroclinic connections between the
unstable homogeneous stationary state and each of the two suspected
global energy minimizers which are indicated by yellow stars in
Figure~\ref{figbifdiagsONE}, and which are shown in the top
left panel of Figure~\ref{figspecificsols}. In other words,
for the above parameter values the diblock copolymer
equation~(\ref{dbcp}) exhibits multistability in the sense
that local or suspected global energy minimizers can be
reached from the homogeneous state.
\end{theorem}
We consider the
specific parameter value~$\lambda = 16\pi^2 \approx 157.91$.
It was shown in~\cite{johnson:etal:13a} that
for this parameter value the diblock copolymer model~(\ref{dbcp}) exhibits
multistability which can be observed in practice. More 
precisely, if one randomly chooses initial conditions near
the homogeneous state~$u \equiv \mu = 0$, then some of the
resulting solutions converge to the suspected global energy
minimizer, while some are trapped in a local energy minimum
At this value, the diblock copolymer model has eight nontrivial
equilibrium solutions, which are shown in Figure~\ref{figspecificsols}.
The top left image shows the two suspected global energy minimizers, 
which are marked by a star in Figure~\ref{figbifdiagsONE},
while the two local minimizers are indicated by a yellow
circle, and are shown in the top right panel of
Figure~\ref{figspecificsols}. Finally, the four colored 
solutions in the two lower panels of Figure~\ref{figspecificsols}
are of index one, and they all lie on the red secondary bifurcation
branch in Figure~\ref{figbifdiagsONE}. We would like to point
out that the equilibrium solutions on the first and second branch 
bifurcating from the horizontal trivial solution line have 
two and three transition layers, respectively. The reason 
for this is described in more detail in~\cite{johnson:etal:13a}.
The above theorem will be established as Theorem~\ref{thmmain},
which for computational reasons considers a rescaled version of the
diblock copolymer model.

We remark that both of the local and global minimizers are symmetric.
  In Figure~\ref{figspecificsols} it is clearly seen that the dominant harmonic of
  the local minimizer is $q=3$,
  whereas the dominant harmonic of the global minimizer
  is $q=2$ (the $a_q\cos(qx)$ term of the solution's expansion in Fourier series dominates).
  Also, the space of sequences $(a_k)$ with $a_k=0$ for $k\neq q\mathbb{N}$
  is invariant under the flow, which is also suggested by the numerical data from
  Appendix~\ref{appdix}. We make use of those symmetries in the construction of
  the stretched unstable isolating blocks, presented in Section~\ref{secproof}.
  Precisely, to obtain the connection to the suspected global minimizer we stretch the block
  in the direction of the second eigenfunction $\cos{2x}$, whereas to obtain the connection to the local minimizer
we stretch the block in the direction of the third eigenfunction $\cos{3x}$.

\subsection{Framework for proving constructively connecting orbits in parabolic PDEs}
In recent years, computer-assisted proof techniques have been used
extensively in the context of nonlinear partial differential equations,
and they are based on a number of different approaches. Here we present
our approach for the computer assisted proof of Theorem~\ref{thmmainintro}.

Our method is general and can be used to achieve similar results for other dissipative PDEs.
The crux of the whole algorithm for the proof of Theorem~\ref{thmmainintro}
is a procedure of rigorous propagation of a piece of the unstable manifold of
the homogeneous state with respect to time. This propagation has to lead to
small interval bounds, while at the same time entering the basin of attraction
of the stable fixed point. For interesting parameter values the global attractor
exhibits a complicated equilibrium structure, and the dynamical equation is rather
stiff. This leads to a numerical propagation of error bounds involving
lots of integration steps. For example the successful proof of the main result required
performing of thousands of numerical integration steps.
To address this problem we have developed an efficient algorithm for
rigorous integration of dissipative PDEs forward in time \cite{Cy}.
Its highlight is an implementation of validated \emph{fast Fourier transforms}, and generic and
efficient implementation of a time-stepping scheme in C++ programming language
combining several techniques, including
\emph{the interval arithmetic, automatic differentiation, Taylor method,
and Lohner's algorithm} \cite{KZ,Lo}. By calling it generic we mean
that we developed a software tool that can be easily adapted to other dissipative PDEs, for example we have
successfully integrated other equations including the Burgers, Swift-Hohenberg and Kuramoto-Sivashinsky PDEs \cite{Cy,C1,CZ}.
Performing all required computations for the proof of the theorem presented in this paper takes only about 10 minutes on a laptop.
All the numerical computations we performed are reproducible, the numerical data from the proofs, and
the software codes with compile/use instructions are published online as a bitbucket repository \cite{codes}.
We are convinced that the method is able to handle large integration times within a
reasonable computational time frame, and in a future work we will explore applicability of our algorithm to parameter regimes
in which the studied model exhibits much more complicated attractor structure and numerical stiffness.
Explicit bounds for the unstable manifold and basins of attraction of
the stationary solutions are computed using now standard techniques based on cone conditions \cite{WZ,Zcc,ZG}
and logarithmic norms \cite{C1,CZ,Zattr}. Computation of explicit bounds for the unstable manifold
and basins of attraction of the stationary solutions is computationally very cheap
compared to the numerical integration part.

\subsection{Structure of the paper}
The remainder of the paper is organized as follows. In
Section~\ref{secsetting} we describe the general setting for our
approach, while Section~\ref{secisolatingblocks} is devoted to its
theoretical foundations. We recall basic definitions and results for
isolating blocks, self-consistent bounds, as well as the cone condition.
Based on this, Section~\ref{secproof} contains the proof of our main
theorem. The next three sections are concerned with the computational
aspects of our work. While the rigorous integration algorithm is described
in Section~\ref{secrigalg}, the following Section~\ref{secccverif} shows 
how the cone condition can be verified in infinite dimensions. General
remarks about the software can be found in Section~\ref{secsoftware}.
Section~\ref{secconclusions} contains conclusions and future
plans. Finally, Appendix~\ref{appendixnumdata} contains some numerical
data from the computer assisted part of the proof. Throughout, from
now on, we will drop the label ``suspected'' when we talk about
the two suspected global minimizers in the above situation.
\section{Problem formulation and basic setting}
\label{secsetting}
We begin by presenting the basic setup for our computer-assisted
proof. As was mentioned in the introduction, our goal is to consider
the diblock copolymer model~(\ref{dbcp}) on the domain $\Omega = (0,1)$.
Unfortunately, for the parameter range of interest this choice of domain
leads to extreme stiffness in the discretized equations. We will address
this issue through rescaling techniques, and as a consequence our basic
setup will be for a general one-dimensional domain of the form~$\Omega = (0,L)$.
The length~$L$ will be chosen later in such a way to reduce the stiffness
of the problem. Thus, we study problem \eqref{dbcp} in one dimension in
the form
\begin{equation} \label{dbcp1d}
  u_t = - \left( u_{xx} + \lambda f(u) \right)_{xx} -
        \lambda\sigma(u-\mu) \; ,
\end{equation}
which is equivalent to~(\ref{dbcp}) up to a rescaling of time.
In this formulation, we use the parameter $\lambda = 1 / \epsilon^2$,
and consider the partial differential equation subject to the mass and
boundary constraints given by
\begin{displaymath}
  \mu = \frac{1}{L} \int_0^L{u(t,x)\,dx}
  \quad\mbox{ and }\quad
  u_x(t,x) = u_{xxx}(t,x) = 0
  \quad\mbox{ for }\quad
  t > 0 \; , \;\; x \in \{ 0, L \} .
\end{displaymath}
As before, the nonlinearity is given by $f(u)=u-u^3$, and we consider
the case $\mu=0$ of zero total mass, which means that we have equal 
amounts of the two polymers~A and~B. 

Our basic discretization of this infinite-dimensional problem is
based on the Fourier series expansion. More precisely, due to the 
imposed homogeneous Neumann boundary conditions we consider the
expansion
\begin{equation} \label{uexpansion}
  u(t,x) = a_0(t) + 2\sum_{k=1}^\infty {a_k(t)\cos{\frac{\pi kx}{L}}}
  \; ,
\end{equation}
which is based on the \emph{cosine Fourier basis}. We would like to
point out that due to our choice $\mu = 0$ for the total mass we
automatically have
\begin{displaymath}
  a_0(t) = 0
  \quad\mbox{ for all }\quad
  t \ge 0 \; .
\end{displaymath}
Notice also that the domain length~$L$ enters the expansion. We will
see later that by choosing~$L$ large enough the stiffness of the
problem~\eqref{dbcp1d} can be reduced, since the eigenvalues of the
linear part of a suitable finite-dimensional approximation are more
averaged. To simplify notation, we will use the abbreviation
\begin{displaymath}
  e_k(x) := \cos{\frac{\pi kx}{L}} \; ,
  \quad\mbox{ where }\quad k \geq 0 \; ,
\end{displaymath}
for the basis functions in the above expansion. By default we assume
that the coefficients~$a_k$ are dependent on time~$t$, and we therefore
frequently drop the explicit mentioning of the temporal argument.
Upon substituting
the series representation of~$u$ into the cubic part of the nonlinearity~$f$
one obtains
\begin{eqnarray*}
  u^3 & = &
    \left(2 \sum_{k \ge 1} {a_k e_k(x)}\right)^3
    \quad = \quad
    8 \sum_{j_1 \geq 1} \sum_{j_2 \geq 1} \sum_{j_3 \geq 1}
      {a_{j_1} a_{j_2} a_{j_3} \cdot \cos{\frac{\pi j_1 x}{L}} \cdot
      \cos{\frac{\pi j_2 x}{L}} \cdot\cos{\frac{\pi j_3 x}{L}}} \\[2ex]
  & = &
    2 \sum_{j_1,j_2,j_3\geq 1} a_{j_1} a_{j_2} a_{j_3} \left(
      \cos{\frac{\pi}{L}(j_1-j_2-j_3) x} + \cos{\frac{\pi}{L}(j_1+j_2-j_3) x}
      \right. \\
  & & \qquad\qquad\qquad\qquad \left. +
      \cos{\frac{\pi}{L}(j_1-j_2+j_3) x} + \cos{\frac{\pi}{L}(j_1+j_2+j_3) x}
      \right) \, .
\end{eqnarray*}
The last expression can be further simplified using the fact
that
$\cos{\frac{\pi kx}{L}} = \frac12\left( \exp\frac{i\pi kx}{L}
    + \exp\frac{-i\pi kx}{L} \right) \; $,
%
%
as long as we recall $a_0 = 0$ and define~$a_{-k} = a_k$ for
$k \in \mathbb{N}$. By reordering one further
obtains
%
\begin{eqnarray*}
  u^3 = \sum_{\substack{j_1,j_2,j_3\in\mathbb{Z}\\ j_1+j_2+j_3=0}}
    a_{j_1} a_{j_2} a_{j_3}
    + 2 \sum_{k \in \mathbb{N}} \left(
    \sum_{\substack{j_1,j_2,j_3\in\mathbb{Z}\\ j_1+j_2+j_3=k}}
    a_{j_1} a_{j_2} a_{j_3} \right)
    \cos{\frac{\pi k x}{L}}
\end{eqnarray*}
If we finally substitute this expression for~$u^3$, together with
the original series expansion for~$u$, into the differential
equation~(\ref{dbcp1d}), then one obtains after division by~$2$
and extraction of the coefficients of the basis functions~$e_k(x)$
the identities
\begin{equation} \label{eq:infDim}
  \frac{d a_k}{d t} = 
  \left(-\frac{k^4 \pi^4}{L^4} + \frac{\lambda k^2 \pi^2}{L^2} -
    \lambda\sigma \right) a_k - \frac{\lambda k^2 \pi^2}{L^2}
    \sum_{\substack{j_1,j_2,j_3 \in \mathbb{Z}\\j_1+j_2+j_3=k}}
    {a_{j_1} a_{j_2} a_{j_3}}
  \quad\text{ for all } \quad k > 0 \; .
\end{equation}
In this way, we have reformulated the original parabolic partial
differential equation as an infinite system of ordinary differential
equations. This system will be used throughout the remainder of
the paper. Recall one more time that the coefficient sequence~$a_k$
satisfies
\begin{displaymath}
  a_0 = 0 \; ,
  \quad\mbox{ as well as }\quad
  a_k = a_{-k} \in \mathbb{R}
  \quad\mbox{ for all }\quad
  k \in \mathbb{Z} \; .
\end{displaymath}
To close this section, we present two central definitions.
The first one is concerned with the function space which
forms the basis of our computer-assisted proof.
\providecommand{\subspace}{H^\prime}
\begin{definition}[Sequence space with algebraic coefficient decay]
\label{defhprime}
Let~$H$ denote the space $\ell^2(\mathbb{Z},\mathbb{R})$, i.e.,
elements~$a \in H$ are sequences $a \colon \mathbb{Z} \to
\mathbb{R}$ such that $\sum_{k\in\mathbb{Z}}{|a_k|^2}\leq\infty$.
In addition, let $\widetilde{H} \subset H$ denote the subspace
of~$H$ which is defined by
\begin{displaymath}
  \widetilde{H} := \left\{ \left\{ a_k \right\}_{k\in\mathbb{Z}}
  \in H \colon
  \text{ there exists a constant } C \ge 0
  \text{ such that } \left| a_k \right| \leq
    \frac{C}{\nmid k\nmid^6}
    \text{ for } k\in\mathbb{Z} \right \} \; ,
\end{displaymath}
where~$\nmid x \nmid = |x|$ for $x \neq 0$ and
$\nmid 0 \nmid = 1$. Finally, let the space~$H^\prime$
be given by
\begin{displaymath}
  \subspace := \widetilde{H} \cap 
    \left\{ \left\{ a_k \right\}_{k\in\mathbb{Z}} \colon
    a_0 = 0 \text{ and }
    a_{k} = a_{-k} \text{ for all } k\in\mathbb{Z}
    \right\} \; .
\end{displaymath}
Notice that all three spaces are Hilbert spaces.
\end{definition}
From now on, all computations we perform in infinite dimensions
will be constrained to the space~$\subspace$. Needless to say,
any actual numerical computations performed on a computer have
to take place in finite dimensions, and for this we make
use of \emph{Galerkin approximations} of \eqref{eq:infDim},
which are the subject of the following final definition of
this section.
\begin{definition}[Galerkin approximation]
Consider the infinite system of ordinary differential equations
given in~(\ref{eq:infDim}), and let $n \in \mathbb{N}$ be 
arbitrary. Then the $n$-th \emph{Galerkin approximation} of
\eqref{eq:infDim} is defined as
\begin{equation} \label{eq:galerkin}
  \frac{d a_k}{d t} =
  \left( -k^2\left(\frac{\pi}{L}\right)^2
    \left(k^2\left(\frac{\pi}{L}\right)^2 - \lambda\right) -
    \lambda\sigma\right) a_k - 
  \lambda k^2 \left(\frac{\pi}{L}\right)^2
    \sum_{\substack{|j_1|,|j_2|,|j_3| \leq n\\j_1+j_2+j_3=k}}
    {a_{j_1} a_{j_2} a_{j_3}} \; ,
\end{equation}
for any coefficient index $k = 1,\dots,n$. In other words,
the $n$-th Galerkin approximation is an $n$-dimensional system
of ordinary differential equations, which depends only on the
coefficients~$a_1, \ldots, a_n$.
\end{definition}
\section{Isolating blocks, self-consistent bounds, and the cone condition}
\label{secisolatingblocks}
In this section we briefly recall the theoretical framework that is
used to establish the main result of this paper. More precisely, we
review the concepts of \emph{cone conditions} and \emph{self-consistent
bounds} as introduced by Zgliczy\'nski. In Section~\ref{secfindim},
we present a finite-dimensional version of \emph{$h$-sets}, \emph{isolating
blocks}, and \emph{cone conditions}. This is followed in Section~\ref{secscb}
by a discussion of \emph{self-consistent bounds}, which allows one to extend
the notion of $h$-sets, isolating blocks, and cone conditions to the Hilbert
space setting. Finally, in Section~\ref{secccinfdim} we present the
concept of infinite-dimensional cone conditions, i.e. we show how to verify
cone conditions on infinite-dimensional self-consistent bounds. Since the results presented
in this section are known, our presentation will be short, albeit as self-contained
as possible. For more detailed discussions of $h$-sets, isolating blocks, and
cone conditions we refer the reader to~\cite{WZ, Zcc, ZG}, and references
therein.

The concepts introduced in the present section will be used in two different
ways in the proof of our main theorem, which is the subject of the next
Section~\ref{secproof}. On the one hand, we use the fact that the construction
of an \emph{unstable} infinite-dimensional isolating block which satisfies
a cone condition provides us with explicit bounds for the unstable manifold
of the unstable steady state contained in the block. On the other hand, the
construction of a \emph{stable} infinite-dimensional isolating block which
satisfies a cone condition guarantees the local uniqueness of the stable
stationary solution contained in the block, as well as explicit bounds
for the size of its basin of attraction.
\subsection{Isolating blocks and cone conditions in finite dimensions}
\label{secfindim}
We begin by reviewing the concepts of isolating blocks and cone
conditions in a finite-dimensional setting. For this, let $n \in
\mathbb{N}$ be arbitrary. Then for any vector $z \in\mathbb{R}^n$
we let~$\|z\|$ denote a norm on~$\mathbb{R}^n$, which does not have
to be the standard Euclidean norm. The basic types of sets necessary
for our discussion are interval sets. For this, consider real 
numbers~$x_k^\pm$ for $k = 1,\ldots,n$. Then we define the
\emph{interval set}~$[x]$ by
\begin{displaymath}
  [x] = \prod^n_{k=1}{\left[ x^-_k, x^+_k \right]}
  \subset\R^n \; ,
  \quad\mbox{ as long as }\quad
  x_k^- \leq x_k^+
  \quad\mbox{ for }\quad
  k = 1,\ldots,n \; .
\end{displaymath}
For any point $x_0 \in \R^n$ and arbitrary $r > 0$, we let
$B_n(x_0, r) = \{ z \in \R^n\ |\ \|z-x_0\| < r \}$ denote the
open ball of radius~$r$ centered at~$x_0$, and we use the 
abbreviation $B_n = B_n(0,1)$ for the unit ball centered at
the origin. For a linear map $A \colon \R^n \to \R^n$, we
let~$\Sp(A)$ denote the spectrum of~$A$. When considering
product spaces such as~$\R^u\times\R^s$, we write elements
$z \in \R^u \times \R^s$ in the form $z = (x,y)$, where~$x
\in \R^u$ and~$y \in \R^s$ denote the first and second
component vectors of~$z$, respectively, and the superscripts
$s$, and $u$ refer to the stable and unstable spaces
respectively. Moreover, if
$f \colon \mathbb{R}^n \to \mathbb{R}^n$ denotes a continuously
differentiable function, and if $Z \subset \R^n$ is a compact
set, then we define
\begin{displaymath}
  [d f(Z)] = \left\{ M = (M_{ij})_{i,j=1}^n \in \R^{n\times n}
    \; \left| \;
    M_{ij} \in \left[
      \inf_{z \in Z}{\frac{\partial f_i}{\partial x_j}(z)} \, , \;
      \sup_{z \in Z}{\frac{\partial f_i}{\partial x_j}(z)} \right]
      \right. \right\} \; .
\end{displaymath}
Consider now an ordinary differential equation of the form
\begin{equation} \label{ode}
  z^\prime = f(z) \; ,
  \qquad z \in \R^n \; ,
  \qquad f \in C^2(\R^n,\R^n) \; .
\end{equation}
Then we denote by~$\varphi(t,z_0)$ the solution of \eqref{ode}
which satisfies the initial condition $z(0) = z_0$. Notice that
since we assume the twice continuous differentiability of~$f$,
the solution of this initial value problem is in fact uniquely
determined.

After these preparations, we now turn our attention to the 
concept of $h$-sets. Informally speaking, an $h$-set is a subset
of Euclidean space, which becomes a cube after a suitable choice
of coordinate system. In addition, the letter~$h$ alludes to the
fact that this choice of coordinate system will lead to a phase
portrait for the underlying ordinary differential equation which
exhibits typical hyperbolic-like form. The concept of $h$-sets
originated in the work on covering relations for multi-dimensional
dynamical systems, see~\cite{ZG}. From the perspective of our 
specific application, the considered $h$-sets have a trivial
structure in the sense that the cube is always expressed with
respect to canonical coordinates. Consequently, also the cone
conditions which will be introduced below are expressed in the
canonical coordinate system. More precisely, in our application
to the unstable homogeneous state, two of the dimensions will be
distinguished as unstable, while the homeomorphism~$c_N$ appearing
in the following definition of an $h$-set is just a scaling. Despite
this fact, we now recall the appropriate definitions in their
original form.
\begin{definition}[$h$-set~{\cite[Definition~1]{ZG}}]
\label{defhset}
An \emph{$h$-set} is a quadruple of the form $(N, u(N), s(N),
c_N)$ which satisfies the following three conditions:
\begin{itemize}
\item The set~$N$ is a compact subset of~$\mathbb{R}^n$,
\item the natural numbers~$u(N), s(N) \in \mathbb{N}_0$
satisfy the identity $u(N) + s(N) = n$, and
\item the mapping $c_N \colon \mathbb{R}^n \to \mathbb{R}^n$,
where we write $\mathbb{R}^n = \mathbb{R}^{u(N)} \times
\mathbb{R}^{s(N)}$, is a homeomorphism which satisfies
\begin{displaymath}
  c_N(N) = \overline{B_{u(N)}} \times \overline{B_{s(N)}}
  \; .
\end{displaymath}
\end{itemize}
In the following, we will slightly abuse notation and
refer to an $h$-set as just~$N$, rather than always listing
the full quadruple. In addition, with each $h$-set we associate
the sets
\begin{displaymath}
  \begin{array}{lclclclclcl}
    \DS N_c^- & := & \DS \partial {B_{u(N)}}\times\overline{B_{s(N)}}
      \; ,
    & \quad &
    \DS N_c^+ & := & \DS \overline{B_{u(N)}}\times\partial B_{s(N)}
      \; ,
    & \quad &
    \DS N_c   & := & \DS \overline{B_{u(N)}}\times\overline{B_{s(N)}}
      \; , \\[2ex]
    \DS N^- & := & \DS c_N^{-1}\left(N_c^-\right) \; ,
    & \quad &
    \DS N^+ & := & \DS c_N^{-1}\left(N_c^+\right) \; ,
    & \quad & N & := & c_N^{-1}\left(N_c\right).
  \end{array}
\end{displaymath}
and we define~$\dim(N) := n$.
\end{definition}
The above definition shows that an $h$-set~$N$ is the product of
two closed unit balls in some appropriate coordinate system. The
two integers~$u(N)$ and~$s(N)$ are called the \emph{nominally
unstable} and \emph{nominally stable dimensions}, respectively.
The subscript~$c$ refers to the new coordinates given by the
homeomorphism~$c_N$. Note that if we have $u(N) = 0$, then one
immediately obtains $N^- = \emptyset$, and the identity $s(N) = 0$
implies $N^+ = \emptyset$. The next definition introduces the
concept of a horizontal disk.
\begin{definition}[Horizontal disk~{\cite[Definition~10]{WZ}}]
\label{defhorizontaldisk}
Let~$N$ denote an arbitrary $h$-set. In addition, consider a
continuous mapping $b\colon \overline{B_{u(N)}}\to N$, and define
$b_c = c_N \circ b$. Then we say that~$b$ is a \emph{horizontal
disk} in~$N$, if there exists a homotopy $h \colon [0,1] \times
\overline{B_{u(N)}} \to N_c$ such that
\begin{eqnarray*}
  h(0,x) = b_c(x)
    & \quad\mbox{ for all }\quad &
    x \in \overline{B_{u(N)}} \; , \\
  h(1,x) = (x,0)
    & \quad\mbox{ for all }\quad &
    x\in \overline{B_{u(N)}} \; , \\
  h(t,x) \in N_c^-
    & \quad\mbox{ for all }\quad &
    t \in [0,1]
    \quad\mbox{ and }\quad
    x \in \partial \overline{B_{u(N)}} \; .
\end{eqnarray*}
\end{definition}
In order to handle the hyperbolic structure of $h$-sets,
we use the notion of \emph{cones}, which are defined through
a quadratic form defined on the $h$-set in the following
way.
\begin{definition}[$h$-set with cones~{\cite[Definition~8]{ZG}}]
Let $N \subset \mathbb{R}^n$ denote an arbitrary $h$-set,
and consider the associated splitting $\mathbb{R}^n =
\mathbb{R}^{u(N)}\times\mathbb{R}^{s(N)}$. Furthermore,
let $Q \colon \mathbb{R}^{u(N)}\times\mathbb{R}^{s(N)}
\to \mathbb{R}$ be a quadratic form which is given as the
difference
\begin{displaymath}
  Q(x,y) = \alpha(x) - \beta(y) \; ,
  \quad\mbox{ for arbitrary }\quad
  (x,y) \in \mathbb{R}^{u(N)} \times \mathbb{R}^{s(N)} \; ,
\end{displaymath}
where $\alpha \colon \mathbb{R}^{u(N)} \to \mathbb{R}$ and
$\beta \colon \mathbb{R}^{s(N)} \to \mathbb{R}$ are positive
definite quadratic forms. Then the pair~$(N,Q)$ is called
an \emph{$h$-set with cones}.
\end{definition}
\begin{definition}[Cone condition~{\cite[Definition~10]{ZG}}]
\label{defconecondition}
Consider an~$h$-set with cones given by~$(N,Q)$,
as introduced in the previous definition. In addition, let
$b \colon \overline{B_{u(N)}} \to N$ denote a horizontal disk
as in Definition~\ref{defhorizontaldisk}. Then we say
that~$b$ \emph{satisfies the cone condition (with respect
to~$Q$)}, if and only if for arbitrary points $x, y
\in \overline{B_{u(N)}}$ we have
\begin{displaymath}
  Q(b_c(x) - b_c(y)) > 0
  \quad\mbox{ as long as }\quad
  x \neq y \; .
\end{displaymath}
Recall that we have~$b_c\colon\R^n\to\R^{u(N)}\times\R^{s(N)}$, and $b_c = c_N \circ b$.
\end{definition}
\providecommand{\isoBlck}{\mathcal{W}}
The remaining definitions of this section are concerned
with standard dynamical systems notions. We begin by
recalling the concept of \emph{isolating block} from
\emph{Conley index} theory~\cite{MM}.
\begin{definition}[Isolating block]
Let~$N$ denote an arbitrary $h$-set. We say that~$N$ is
an \emph{isolating block} for the vector field~$f : \R^n \to \R^n$,
if and only if~$c_N$ is a diffeomorphism, if the sets~$N^+$
and~$N^-$ are local sections for~$f$, i.e., the vector field~$f$
is transversal to $N^\pm$, and if the following two conditions
are satisfied:
\begin{itemize}
\item For all points $x \in N^-$ there exists a
constant $\delta>0$ such that $\varphi(t,x) \not\in N$
for all $t \in (0,\delta]$.
\item For all points $x \in N^+$ there exists a
constant $\delta>0$ such that $\varphi(t,x) \not\in N$
for all $t \in [-\delta,0)$.
\end{itemize}
In other words, the closed set~$N^-$ consists of all exit
points for the flow associated with~(\ref{ode}), and the
closed set~$N^+$ is the set of entry points.
\end{definition}
\begin{lemma}[Graph representation~{\cite[Lemma~5]{Zcc}}]
Let~$(N,Q)$ denote an arbitrary $h$-set with cones, and let
$b \colon \overline{B_{u(N)}} \to N$ be a horizontal disk which
satisfies the cone condition as in Definition~\ref{defconecondition}.
Then there exists a Lipschitz continuous function
$y \colon \overline{B_{u(N)}} \to \overline{B_{s(N)}}$
such that $b_c(x) = (x,y(x))$.
\end{lemma}
\begin{definition}[Hyperbolic fixed point]
Let $z_0\in\R^n$ and consider again the ordinary differential
equation~(\ref{ode}). Then~$z_0$ is called a \emph{hyperbolic
fixed point} for~\eqref{ode} if $f(z_0) = 0$ and if
$\re{\lambda}\neq 0$ for all eigenvalues $\lambda \in
\Sp(df(z_0))$, where~$df(z_0)$ is the Jacobian matrix of the
vector field~$f$ at the point~$z_0$, and $\re{\lambda}$ denotes
the real part of the eigenvalue~$\lambda$.

Now consider a hyperbolic fixed point~$z_0$ which is contained
in a given set~$Z \subset \R^n$. Then we define the stable and
unstable manifolds of~$z_0$ as
\begin{eqnarray*}
  W_Z^s(z_0,\varphi) & = &
    \left\{ z \in Z \; \left| \; \varphi(t,z) \in Z
    \quad\mbox{ for all }\quad t \ge 0 \quad\mbox{ and }\quad
    \lim_{t \to +\infty} \varphi(t,z) = z_0 \right.
    \right\} \; \mbox{ and} \\[1.5ex]
  W_Z^u(z_0,\varphi) & = &
    \left\{ z \in Z \; \left| \; \varphi(t,z) \in Z
    \quad\mbox{ for all }\quad t\leq 0 \quad\mbox{ and }\quad
    \lim_{t\to-\infty}{\varphi(t,z) = z_0} \right.
    \right\} \; ,
\end{eqnarray*}
respectively.
\end{definition}
Important for our later applications will be the following
theorem. It shows that for hyperbolic fixed points, one
can always find an $h$-set~$N$ with cones in such a way that
the unstable manifold~$W^u_N(z_0)$ is a horizontal disk in~$N$,
and that the stable manifold~$W^s_N(z_0)$ is a vertical disk.
\begin{theorem}[Invariant manifolds as disks~{\cite[Theorem~26]{Zcc}}]
\label{thmzcc}
Suppose that the point~$z_0$ is a hyperbolic fixed point for
the ordinary differential equation~\eqref{ode}. Moreover, let
$Z \subset \R^n$ denote an open set with $z_0 \in Z$. Then there
exists an $h$-set~$N$ with cones such that the following hold:
\begin{itemize}
\item We have both $z_0 \in N$ and $N\subset Z$.
\item The $h$-set~$N$ is an isolating block for the vector
field~$f$.
\item The unstable manifold~$W^u_N(z_0)$ is a horizontal disk
in~$N$ which satisfies the cone condition.
\item The stable manifold~$W^s_N(z_0)$ is a vertical disk in~$N$
which also satisfies the cone condition.
\end{itemize}
Finally, the unstable manifold~$W^u_N(z_0)$ can be represented as
the graph of a Lipschitz continuous function over the unstable space 
of the linearization of~$f$ at~$z_0$, and this graph is tangent to
the unstable space at~$z_0$. An analogous statement holds for the
stable manifold~$W_N^s(z_0)$.
\end{theorem}
In our application below we are interested in a slightly different
situation. In our case, the existence of the fixed point is not
known ahead of time. Its existence is, however, a consequence of
the existence of an isolating block. Such a theorem in the context
of maps has been proven in~\cite[Theorem~10]{Zcc}. In the context
of differential equations it has been first established
in~\cite{mccord:88b}, see also~\cite{maier:etal:08a, ZM}.
\subsection{Self-consistent bounds}
\label{secscb}
We now turn our attention to the concept of \emph{self-consistent bounds},
which have been used extensively in recent years~\cite{CZ, CzZ, maier:etal:08a,
Zattr, ZPer, ZM, Z3}. In the following, we recall the basic definitions and
fundamental ideas behind this concept, and describe an application which will
be used for our main result.

In its most general form, the method of self-consistent bounds applies to
arbitrary \emph{dissipative evolution equations} in a suitable subspace
of the sequence space~$\ell^2$, subject to certain admissibility conditions
which depend on the specific notion of dissipativity. In other words, one
can suppose that the underlying evolution equation takes the abstract form
\begin{equation} \label{absevo}
  \frac{d u}{d t} = F(u) \; ,
\end{equation}
where~$F$ is defined on a suitable subspace of~$\ell^2$ and takes values
in the sequence space. Rather than presenting the concept of self-consistent
bounds in this general setting, we right away specialize to the situation
of the present paper. Our base space is the Hilbert space $\subspace \subset
\ell^2$ which was introduced in Definition~\ref{defhprime}, and we consider the
infinite system of ordinary differential equations given by~\eqref{eq:infDim}.
In other words, in our situation the abstract equation~(\ref{absevo}) can
be written in coordinate form as
\begin{equation} \label{eq:infDim2}
  \frac{d a_k}{d t} = F_k(a) =
  \left(-\frac{k^4 \pi^4}{L^4} + \frac{\lambda k^2 \pi^2}{L^2} -
    \lambda\sigma \right) a_k - \frac{\lambda k^2 \pi^2}{L^2}
    \sum_{\substack{j_1,j_2,j_3 \in \mathbb{Z}\\j_1+j_2+j_3=k}}
    {a_{j_1} a_{j_2} a_{j_3}}
  \quad\text{ for all } \quad k > 0 \; ,
\end{equation}
where $a_0(t) \equiv 0$. Without going into details, we note that the
right-hand side is well-defined for all sequences~$a \in \subspace$. In
fact, the algebraic decay imposed through the subspace~$\widetilde{H}$
in Definition~\ref{defhprime} guarantees the four-time continuous 
differentiability of the function~$u$ given by~(\ref{uexpansion}),
and therefore standard regularity results for parabolic partial
differential equations ensure the well-posedness of~(\ref{eq:infDim2})
in the space~$\subspace$.

The basic idea behind the method of self-consistent bounds is the
construction of a suitable subset of the phase space of~(\ref{absevo}),
in our case the Hilbert space~$\subspace \subset \ell^2$, for which
one can easily determine the flow across its boundary. Usually,
such a set is constructed as an infinite product of intervals for
each of the coefficients~$a_k$ of the sequences~$a \in \subspace$,
where, in order to be able to handle this set computationally,
 a bound for infinite number of them is provided by an algebraic decay.
More precisely, let~$a_k^-$ and~$a_k^+$ denote the lower and upper 
bounds, respectively, for the $k$-th coefficient~$a_k$, i.e., we will
assume that elements of the constructed set satisfy $a_k \in [a_k^-,a_k^+]$,
as long as~$k$ is sufficiently large. In detail, self-consistent bounds
are defined as follows.
\begin{definition}[Self-consistent bounds]
Consider a set of the form
\begin{displaymath}
  W \oplus T = \left\{ a \in \subspace \; \left| \;
    \left(a_1,\ldots,a_m\right) \in W
    \;\mbox{ and }\;
    \left(a_{m+1},a_{m+2},\ldots\right) \in T
    \right. \right\} \subset \subspace \; ,
\end{displaymath}
which is the product of a finite-dimensional set $W \subset \R^m$ and
an infinite-dimensional \emph{tail} $T \subset\R^\infty$. Then the
set~$W \oplus T$ forms \emph{self-consistent bounds} for the evolution
equation~\eqref{eq:infDim2} if and only if the following conditions
are satisfied.
\begin{itemize}
\item The infinite-dimensional tail~$T$ is given by an infinite
product of intervals in the form
\begin{displaymath}
  \R^\infty \supset T = \prod_{k>m}{\left[ a_k^-, a_k^+ \right]}
  \; .
\end{displaymath}
Moreover, there exists a usually large integer~$M>m$ such that
$a_k^- \le 0 \le a_k^+$ for all $k \ge M$, i.e., all but finitely
many of the intervals in~$T$ contain zero.
\item There exists a constant~$C>0$ and an exponent
$s \geq 6$ such that~$|a_k^\pm| \le C / |k|^s$ for all $k > m$.
\item The vector field~$F$ points inwards on the parts of the
boundary of~$W \oplus T$ which are associated with the components
in the tail~$T$ in the sense that for all sequences $x \in 
W \oplus T$ we have
\begin{displaymath}
  F_k(x) > 0 \;\;\text{ if }\;\; x_k = a_k^- \; ,
  \quad\mbox{ and }\quad
  F_k(x) < 0 \;\;\text{ if }\;\; x_k = a_k^+ \; ,
  \quad\text{ for all }\quad k > m \; .
\end{displaymath}
\end{itemize}
\end{definition}
The evolution equation~\eqref{eq:infDim2} represents a dissipative
partial differential equation in the sense discussed in~\cite{Zattr,
ZPer, Z3}. Moreover, one can show that the right-hand side operator
in~\eqref{eq:infDim2} is in fact continuous on the set~$W \oplus T$.
The \emph{block decomposition} of the underlying Hilbert
space~$\subspace$ is given by the infinite sum
\begin{equation} \label{blockdcmp}
 \subspace = \subspace_1 \oplus \subspace_2 \oplus \ldots \; , 
\end{equation}
where the subspace~$\subspace_k$ is the orthogonal projection onto
the subspace of~$\subspace$ which corresponds to the $k$-th cosine
basis vector via~(\ref{uexpansion}). For all of our considerations
we will not use the standard induced norm on~$\subspace$, but 
rather the so-called \emph{block-infinity} norm given by
\begin{displaymath}
  |x|_{\infty} := \max_{k \in \mathbb{N}}{\left| x_k \right|}
  \quad\mbox{ for all }\quad
  x = \{ x_k \}_{k \in \mathbb{Z}} \in \subspace \; ,
\end{displaymath}
This norm was introduced in~\cite{Zattr}. For later applications,
we also define the projection~$P_n \colon \subspace \to \mathbb{R}^n$
which maps~$a$ to~$(a_1,\ldots,a_n)$ for all $n \in \mathbb{N}$.
Finally, within the set of self-consistent bounds $W \oplus T \subset
\subspace$ for~\eqref{eq:infDim2} we have
\begin{itemize}
\item uniform convergence and existence of a solution for the
  infinite-dimensional system~\eqref{eq:infDim2} in the sense
  that the orbits that stay in the set during a finite time interval
  converge in the block-infitnity norm,
\item the solutions of~\eqref{eq:infDim2} are uniquely determined,
and
\item one can derive a bound for the~\emph{Lipschitz constant}
of the associated flow.
\end{itemize}
For more details on the justifications for each of these
statements we refer the reader to~\cite{Zattr}.
\subsection{Cone conditions in infinite dimensions}
\label{secccinfdim}
In this section we briefly describe how the finite-dimensional
cone conditions from Section~\ref{secfindim} can be interpreted
in the infinite-dimensional self-consistent bounds setting. The
involved results have appeared in two recent works on computer-assisted
proofs for partial differential equations. On the one hand, they were
used to establish a heteroclinic connection for the Kuramoto-Sivashinsky
partial differential equation in~\cite{Zhet}. In addition, in the recent
thesis~\cite{Ma} they are fundamental for the computation of rigorous
enclosures of unstable manifolds in the Cahn-Hillard equation.

As a first step, we have to extend the notion of $h$-set to the
infinite-dimensional self-consistent bounds setting. This can be
accomplished as follows.
\begin{definition}[Self-consistent bounds and $h$-sets]
Suppose that $W \oplus T \subset \subspace$ establishes
self-consistent bounds for the evolution equation~\eqref{eq:infDim2}.
Then the set $W \oplus T$ is an $h$-set, if and only if the finite-dimensional
part~$W$ is an $h$-set in the sense of Definition~\ref{defhset}.
\end{definition}
Our primary interest is the verification of cone conditions for
self-consistent bounds for the evolution equation~\eqref{eq:infDim2}
for the purposes of establishing the existence of an equilibrium
solution in the self-consistent bound set, and to recognize the
unstable manifold at this stationary state as a suitable horizontal
disk. This is the subject of the next theorem, which states that if
the cone condition has been verified for self-consistent bounds
$W \oplus T$ and a specific matrix~$Q$, then there exists a
stationary point $z_0 \in W \oplus T$ for~\eqref{eq:infDim2}, and
the unstable manifold~$W^u_N(z_0)$ is a horizontal disk in
$W \oplus T$ which satisfies the cone condition.
\begin{theorem}[Consequences of the cone condition]
\label{ccthm}
Consider self-consistent bounds $W \oplus T \subset \subspace$
for~\eqref{eq:infDim2}, and assume that the orthogonal projection
$P_n(W \oplus T)$ is an isolating block for the $n$-th Galerkin
approximation~\eqref{eq:galerkin} of~\eqref{eq:infDim2} for
all $n > M$. In particular, there is a finite number of
unstable directions, all of them are confined to the finite
dimensional part of $W \oplus T$, and all the directions
corresponding to the infinite dimensional tail are stable.
In addition, let~$Q$ denote an infinite diagonal
matrix with~$Q_{jj} = 1$ for all indices~$j$ corresponding to
the unstable directions and $Q_{jj} = -1$ for the stable ones.
Using this matrix, we assume that
\begin{equation} \label{ccinfd}
  w^t \left( [d\,P_mF(W)]^T Q + Q[d\,P_mF(W)] \right)w > 0
\end{equation}
for all $w\neq 0$ of the form $w = x-y$, where $x,y \in W$.

Then there exists an equilibrium $z_0 \in W \oplus T$
for~\eqref{eq:infDim2} and the unstable manifold~$W^u_N(z_0)$
is a horizontal disk in~$W \oplus T$ which satisfies the cone
condition, where we set $N = W \oplus T$.
We use here an extended notion of $h$-set as a product of
a compact subset of $\R^M$ and an infinite tail.
Moreover, the
manifold~$W^u_N(z_0)$ can be represented as the graph of a
Lipschitz continuous function over the unstable space.
\end{theorem}
For a detailed proof we refer the reader to~\cite{Ma, Zhet}. Basically,
this theorem follows from the finite-dimensional version given in
Theorem~\ref{thmzcc}, combined with the fact that within self-consistent
bounds $W \oplus T \subset \subspace$ both the Fourier series of all terms
appearing in the right-hand side~$F$ of~\eqref{eq:infDim2} and all
partial derivatives in~$dF$ converge uniformly. This allows one to
establish the existence of unstable manifolds for all Galerkin
approximations of~\eqref{eq:infDim2} with uniform bounds, which in
addition are contained in a compact set. One can now employ a standard
convergent subsequence argument to establish the existence of an
unstable manifold for the infinite-dimensional system~\eqref{eq:infDim2}.
Moreover, the verification of positive definiteness of the infinite-dimensional
matrix $[d\,P_mF(W)]^T Q + Q[d\,P_mF(W)]$ reduces to the verification of positivity
of a finite number of inequalities, due to the polynomial growth of the
diagonal which stems from the quartic increase of the eigenvalues of the
linear part of~\eqref{eq:infDim2}. Finally, since the equation is dissipative,
for sufficiently large~$k$ the eigenvalues of the linear part of~\eqref{eq:infDim2}
are all negative, and therefore we have~$Q_{kk} = -1$ for all these indices.
In order to verify~\eqref{ccinfd}, we make use of the following lemma.
\begin{lemma}[Verifying the cone condition]
\label{lemcomput}
If for some $\varepsilon > 0$ and for all $k \in \mathbb{N}$
one can show that
\begin{equation} \label{cccomput}
  2 \inf_{x \in W} {\left|\frac{\partial F_k}{\partial x_k}(x)\right|} -
  \sum_{\ell \neq k}{\sup_{x \in W}{\left|Q_{\ell\ell}
  \frac{\partial F_\ell}{\partial x_k}(x) +
  Q_{kk}\frac{\partial F_k}{\partial x_\ell}(x)\right|}}
  \geq \varepsilon \; ,
\end{equation}
is satisfied, then the inequality~\eqref{ccinfd} holds.
\end{lemma}
It will be the subject of Section~\ref{secccverif} below how
this infinite-dimensional condition can be verified computationally,
and this will provide the unstable equilibrium which forms the
initial point for our heteroclinic solutions.

Finally, we would like to mention that the existence of a stable
steady state of~\eqref{eq:infDim2} can be accomplished by checking
suitable estimates for \emph{logarithmic norms}, which were
introduced in~\cite{D, L}. In addition, such logarithmic norm
estimates will provide explicit bounds for the size of the basin
of attraction of the stable steady state. For this, we basically
apply the theory, algorithms, and code from~\cite{C1, CZ, Zattr},
and the interested reader is referred to these papers for more details.
We do not further describe the details this verification procedure,
since the logarithmic norm condition is equivalent to the cone
condition~\eqref{cccomput} with only stable directions, i.e., with
the matrix~$Q$ being given by~$-\text{Id}$, where $\text{Id}$ denotes
the infinite-dimensional identity matrix.
\section{Proof of the main result}
\label{secproof}
\providecommand{\fss}{z_0^{s,loc}}
\providecommand{\fssg}{z_0^{s,glob}}
\providecommand{\fsu}{z_0^{u}}
\providecommand{\ws}{\mathcal{W}_{\text{s}}^{loc}}
\providecommand{\wsg}{\mathcal{W}_{\text{s}}^{glob}}
\providecommand{\wu}{\mathcal{W}_{\text{u}}}
\providecommand{\us}{\overline{u}_{\text{s}}}
\providecommand{\uu}{\overline{u}_{\text{u}}}
With this section we begin proving our main result, which
was described in Theorem~\ref{thmmainintro} in the introduction.
While this result guarantees the existence of heteroclinic
connections between an unstable index two stationary solution,
which in fact is the homogeneous steady state, and two stable
local energy minimizers, its formulation concentrated on the
diblock copolymer model on the domain~$\Omega = (0,1)$. This
choice was motivated by the studies in~\cite{johnson:etal:13a},
yet for numerical reasons it is not optimal. In fact, by naively
choosing the parameters as in Theorem~\ref{thmmainintro}
one ends up with an extremely stiff Galerkin
approximation~(\ref{eq:galerkin}).

In order to overcome this stiffness problem, we make use of
standard rescaling arguments. For a fixed approximation dimension
taking larger domain size results in 'averaging' the set of
linear part eigenvalues appearing in the approximation of a fixed dimension.
For this, assume that~$u$ denotes
a solution of the diblock copolymer model~(\ref{dbcp}) on the
one-dimensional domain~$\Omega = (0,1)$, and for fixed parameters
$\lambda = 1 / \epsilon^2 > 0$ and $\sigma > 0$. Furthermore,
let~$L > 0$ be a fixed interval length, and define a
function~$v : \mathbb{R}_0^+ \times (0,L) \to \mathbb{R}$
via
\begin{displaymath}
  v(\tau,y) = u\left( \frac{\lambda \tau}{L^4} \, , \;
    \frac{y}{L} \right)
  \quad\mbox{ for all }\quad
  \tau \ge 0
  \quad\mbox{ and }\quad
  y \in (0,L) \; .
\end{displaymath}
Then one can easily verify that~$v$ solves the diblock copolymer
equation
\begin{displaymath}
  v_\tau = -\left( v_{yy} + \frac{\lambda}{L^2} \cdot f(v)
    \right)_{yy} - \frac{\sigma}{L^2} \cdot \frac{\lambda}{L^2}
    \cdot (v-\mu)
  \quad\mbox{ on the domain }\quad
  \mathbb{R}_0^+ \times (0,L) \; ,
\end{displaymath}
which is exactly of the form given in~(\ref{dbcp1d}), but with
rescaled parameters~$\lambda$ and~$\sigma$.
\begin{figure}[tb]
  \centering
  \setlength{\unitlength}{1 cm}
  \begin{picture}(15.0,6.0)
    \put(0.0,0.0){%
      \includegraphics[width=7.0cm]{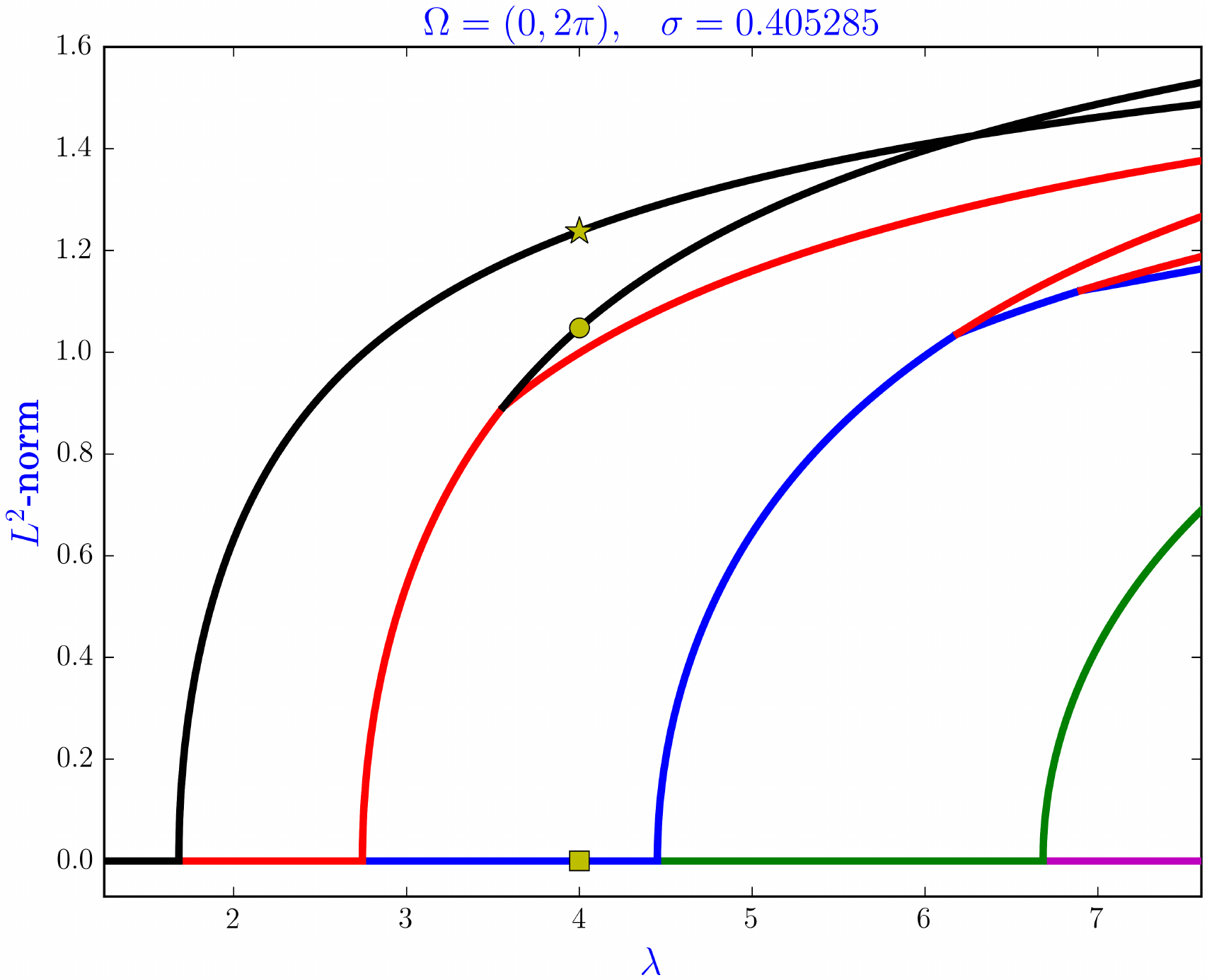}}
    \put(8.0,0.0){%
      \includegraphics[width=7.0cm]{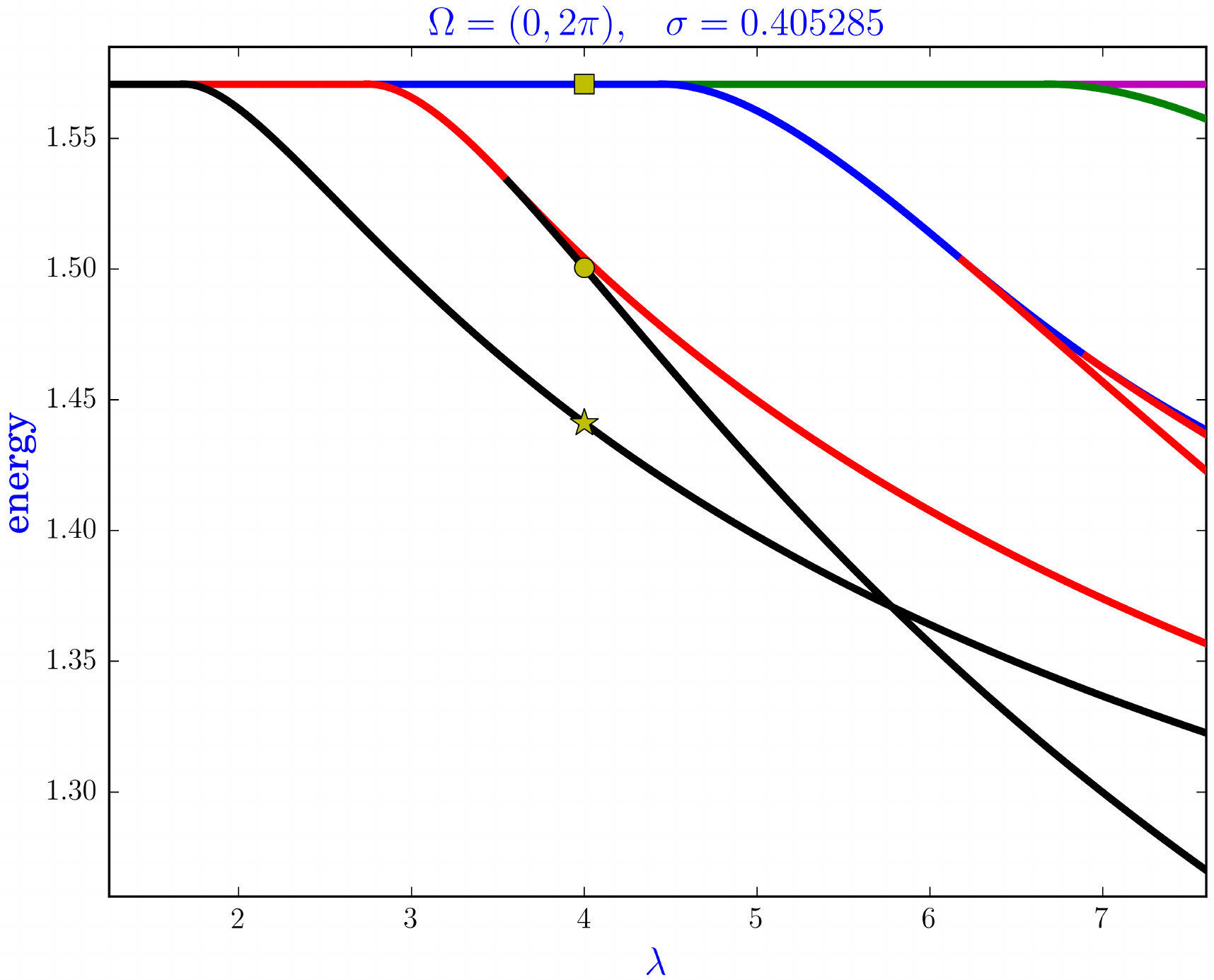}}
  \end{picture}
  \caption{Bifurcation diagrams for the diblock
           copolymer model~(\ref{dbcp}) on the (rescaled) domain~$\Omega
           = (0,2\pi)$, for total mass $\mu = 0$, and for nonlocal
           interaction parameter $\sigma = 16 / (2\pi)^2 \approx 0.405285$.
           In the left diagram, the vertical axis measures the
           $L^2(0,2\pi)$-norm of the solutions, while the right diagram
           shows the total energy~(\ref{dbcp:energy}). In both diagrams,
           the horizontal axis uses the parameter $\lambda =
           1 / \epsilon^2$. As in Figure~\ref{figbifdiags}, the
           solution branches are color-coded by the Morse index
           of the solutions, and black, red, blue, green, and
           magenta correspond to indices~$0$, $1$, $2$, $3$,
           and~$4$, respectively.}
  \label{figbifdiags2PI}
\end{figure}

It turns out that for the parameter combinations given in
Theorem~\ref{thmmainintro} it is most convenient to choose
the new domain length~$L = 2\pi$, since this reduces the size
of the eigenvalues of the linearized Galerkin
approximation~(\ref{eq:galerkin}). Thus, in the following,
we will establish the following rescaled result.
\begin{theorem}[Rescaled version of the main theorem]
\label{thmmain}
Consider the diblock copolymer equation~(\ref{dbcp1d}) on the
domain~$\Omega = (0,L)$ with $L = 2\pi$, for
interaction lengths $\lambda = 4$ and $\sigma = 4 / \pi^2
\approx 0.405284735$, and for total mass $\mu = 0$. Then there
exist heteroclinic connections between the unstable homogeneous
stationary state $u \equiv 0$, which is indicated by a yellow
square in Figure~\ref{figbifdiags2PI}, and 
\begin{enumerate}
\item[(a)] Each of the two local energy minimizers which are
indicated by yellow circles, as well as
\item[(b)] Each of the two global energy minimizers which are
indicated by yellow stars.
\end{enumerate}
In other words, for the above parameter combinations, the
diblock copolymer equation exhibits multistability, and all
energy minimizers can be reached from any arbitrarily small
neighborhood of the homogeneous state.
\end{theorem}
To illustrate the formulation of the above theorem, we have
included the bifurcation diagrams for the new, rescaled
situation in Figure~\ref{figbifdiags2PI}. Compare also
their original counterparts in Figure~\ref{figbifdiagsONE}.
\begin{figure}[tb]
  \begin{center}
    \includegraphics[width=0.5\textwidth]{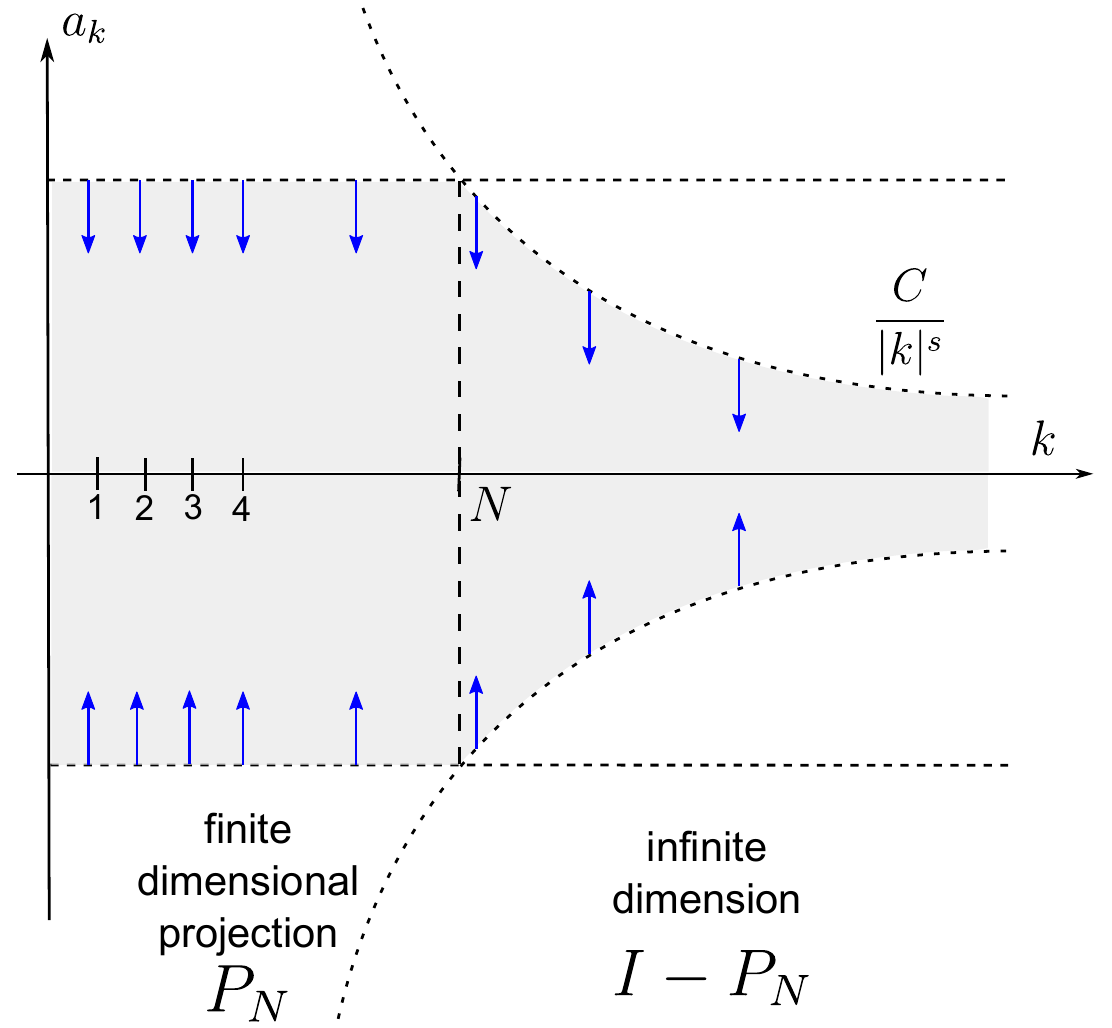}
  \end{center}
  \caption{The stable isolating block~$\ws$ and its
           infinite-dimensional self-consistent bounds structure.
           The isolating block has the topological property that the
           vector field in~(\ref{eq:infDim2}) points inwards in all
           coordinate directions.}
  \label{figws}
\end{figure}
\begin{figure}[tb]
  \begin{center}
    \includegraphics[width=0.5\textwidth]{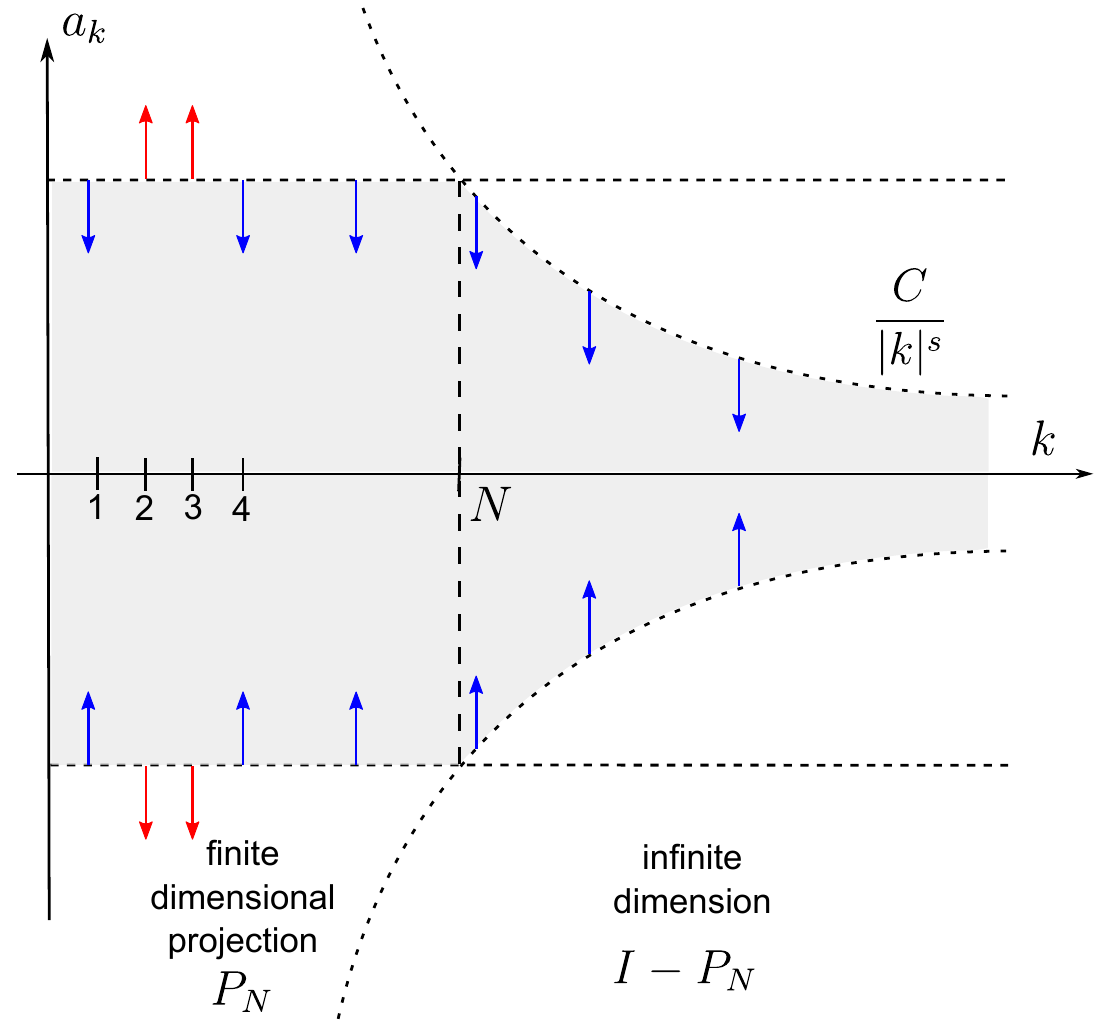}
  \end{center}
  \caption{The unstable isolating block~$\wu$ which contains the unstable
           homogeneous equilibrium of~(\ref{dbcp1d}) at the origin, and
           its infinite-dimensional self-consistent bounds structure.
           This isolating block has the topological property that the
           vector field in~(\ref{eq:infDim2}) points outwards for the
           second and third coordinate directions, while pointing inwards
           for all remaining directions. Once the cone condition given in
           Lemma~\ref{lemcomput} has been verified for this block, the
           existence of a two-dimensional unstable manifold follows, as
           well as the fact that it is a horizontal disk within the
           depicted set.}
\label{figwu}
\end{figure}

Before turning our attention to the proof of Theorem~\ref{thmmain},
we need to introduce some notation. In the following, we denote 
the unstable homogeneous steady state of~(\ref{dbcp1d}) at origin
by $\fsu \in \subspace$. This equilibrium is indicated by a yellow
square in Figure~\ref{figbifdiags2PI}. In addition, let
$\fss \in \subspace$ denote one of the two stable local energy
minimizers of~(\ref{dbcp:energy}), which are indicated by
a yellow circle in Figure~\ref{figbifdiags2PI}. 
Let $\fssg \in \subspace$ denote one of the two global energy 
minimizers of~(\ref{dbcp:energy}), which are indicated by yellow
star in Figure~\ref{figbifdiags2PI}. 
One can easily
see that the theorem only has to be verified for one of these 
local/global minima, 
the result for the second one follows from symmetry arguments.

For each of the involved equilibrium solutions, we
need to construct self-consistent bounds containing them. Thus,
we let
\begin{displaymath}
  \ws \subset \subspace
\end{displaymath}
denote an infinite-dimensional stable isolating block which
contains the stable steady state~$\fss$, and which forms
self-consistent bounds. We also let
\begin{displaymath}
  \wsg \subset \subspace
\end{displaymath}
denote the infinite-dimensional stable isolating block which
contains the stable steady state~$\fssg$

As long as it also satisfies the cone
condition, the existence of stable steady states in~$\ws$ and $\wsg$
follows then from Theorem~\ref{ccthm} with all directions
being stable directions, as illustrated in Figure~\ref{figws}.
Now consider the unstable equilibrium involved in the heteroclinic.
Let
\begin{displaymath}
  \wu \subset \subspace
\end{displaymath}
denote the infinite-dimensional unstable isolating block which
contains the known unstable equilibrium~$\fsu = 0$, and which forms
self-consistent bounds. In this case, if the block satisfies
the cone condition, then the existence of a two-dimensional
unstable manifold is guaranteed at the origin by Theorem~\ref{ccthm},
as is the fact that the manifold is a horizontal disk in~$\wu$.
See also the illustration in Figure~\ref{figwu}. After these
preparations, we can now proceed with the proof of
Theorem~\ref{thmmain}.
\begin{proof}
The verification of Theorem~\ref{thmmain} is divided into three parts,
which are described in more detail below. We would like to point out
already now that all sets used in the proof are constructed using
rigorous numerical programs which are based on \emph{interval
arithmetic}, see for example~\cite{Mo}. In this way, all computations
done by a computer are validated in a strict mathematical sense.
In particular, if the program verifies the cone condition for the
whole set under consideration, this automatically implies that the
condition is satisfied for all elements of this set. For more
illustrations of this approach, also in other contexts, we refer
the interested reader to~\cite{T}.
\begin{figure}[tb]
  \begin{center}
    \includegraphics[width=0.7\textwidth]{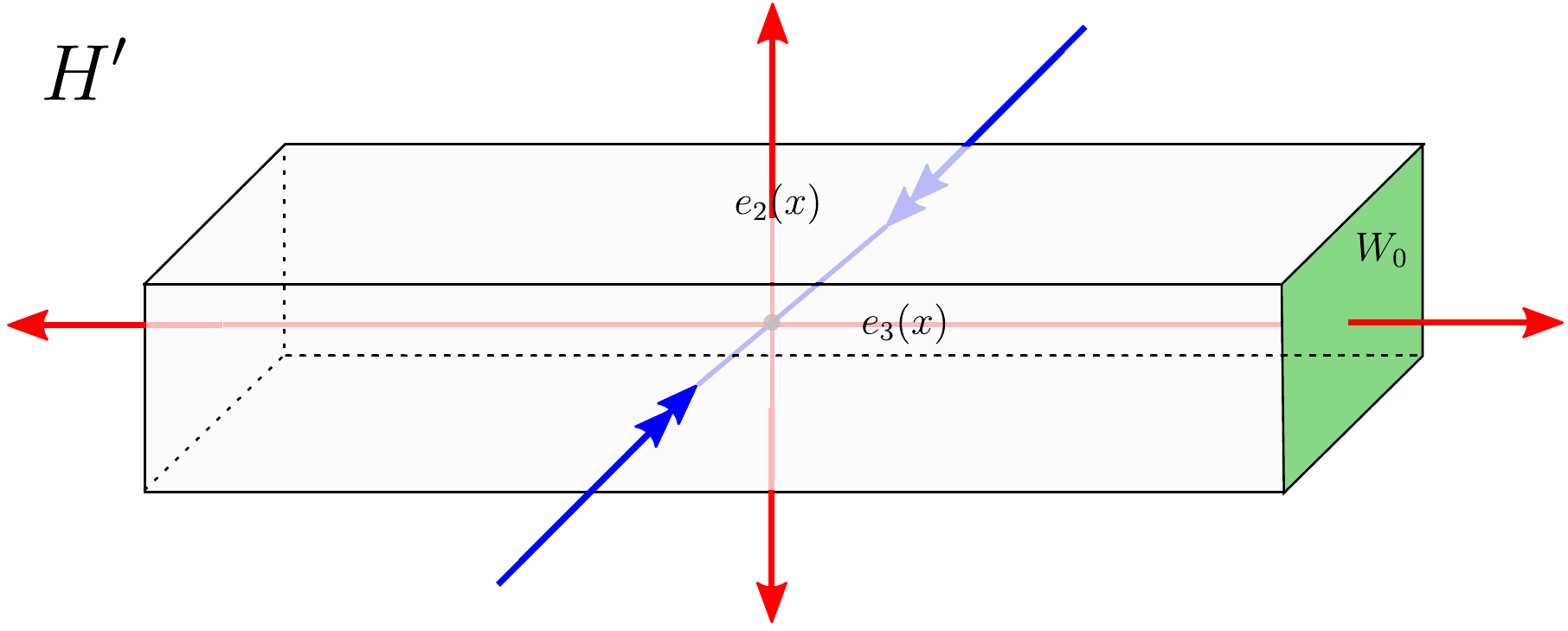}
  \end{center}
  \caption{Sketch of the isolating block~$\wu$ centered at the unstable
           equilibrium~$\fsu$, which is the homogeneous zero state. The
           eigenfunctions~$e_2(x)$ and~$e_3(x)$ correspond to the two unstable
           directions, as can be easily seen by determining the linearization
           of~\eqref{eq:infDim2} at zero. For establishing the existence 
           of the heteroclinic, the face~$W_0$ should be preferably small,
           and as far away from the unstable stationary point as possible,
           as it will be propagated forward
           in time through rigorous integration. We illustrate the unstable
           isolating block that we construct in the case of the connection
           to the local minimizer. Observe that the set is elongated in
           the direction of the third eigenvector, and is very much thinner
           in all the other directions.
           For the other case of the connection to the
           global minimizer the isolating block looks qualitatively the same,
           but it is elongated along $e_2$ eigendirection instead of $e_3$.}
  \label{figublock}
\end{figure}

{\bf First step:} Using a computer program based on interval arithmetic,
we first construct self-consistent bounds $\wu \subset \subspace$ in
such a way that~$\wu$ is an unstable isolating block for~\eqref{eq:infDim2}.
For this, the set~$\wu$ is centered at the origin, which represents our
known unstable equilibrium~$\fsu$ for~\eqref{eq:infDim2}. More specific
information regarding the specific choice of~$\wu$ is presented in
Appendix~\ref{appendixnumdata}, which contains verified numerical
data from the program. These results show that the topological
structure of~$\wu$ is exactly of the form shown in Figure~\ref{figwu}.
In addition, the set~$\wu$ is verified to satisfy the cone
condition~\eqref{cccomput}, using the method presented in Section~\ref{secccverif}.
Now Theorem~\ref{ccthm} implies that
the unstable manifold~$W_N^u(z_0^u)$ is a horizontal disk in~$\wu$,
and therefore the unstable manifold at the origin crosses one of
the faces of~$\wu$. Let~$W_0$ denote such a face, see also the diagram
in Figure~\ref{figublock}. The whole construction of~$\wu$ relies on the
fact that the unstable equilibrium~$\fsu$ is given by the origin. We can
therefore determine the eigenfunctions of the linearization
of~\eqref{eq:infDim2} at~$\fsu$ explicitly, and this easily shows
that all eigenfunctions are just the basis functions~$e_k(x)$
for $k \in \mathbb{N}$. Furthermore, recall from our brief
discussion in the introduction that for our choice of~$\sigma$ 
the two unstable eigendirections correspond to the
eigenfunctions~$e_2(x)$ and~$e_3(x)$. In order to make the
subsequent parts of the proof computationally faster,
the set~$\wu$ is chosen in such a way that it is much more elongated in the second eigenfunction direction compared
to all the other directions, which is given
by~$e_3(x)$ (we denote the face of $\wu$ in this direction by $W_0(e_3)$) 
for the case of the connection to the local minimizer.
Additionally, the set~$\wu$ is chosen in such a way that
it is much more elongated in the third eigenfunction $e_2(x)$ direction for the case of the connection to the global 
minimizer (we denote the face of $\wu$ in this direction by $W_0(e_2)$). See again Figure~\ref{figublock}.

{\bf Second step:} Similar to the first step, we use interval
arithmetic to construct self-consistent bounds $\ws,\wsg \subset
\subspace$ in such a way that~$\ws,\wsg$ is a stable isolating block
for~\eqref{eq:infDim2}. The block~$\ws$ is centered at a numerically
computed approximation to the actual stable equilibrium~$\fss$,
and the topological structure of~$\ws$ is as indicated in
Figure~\ref{figws}. Analogously, the block $\wsg$ is centered at 
a numerically computed approximation to the stable equilibrium~$\fssg$
-- the global minimizer.
The finite-dimensional part of this set
is determined in an appropriate coordinate system, which is 
given by the numerically determined eigenbasis of the Jacobian
matrix at the numerical approximation of~$\fss$, $\fssg$ respectively. 
As a consequence,
the block decomposition~\eqref{blockdcmp} of~$\subspace$ takes
on the following form. The first~$m$ coordinates correspond
to the eigenbasis of the $m$-th Galerkin
approximation~\eqref{eq:galerkin}, while the remaining 
coordinate directions are given by the standard basis
functions~$e_k(x)$ for $k > m$. In other words, we have to
take into account that the Jacobian at~$\fss$, $\fssg$ is considered
in two different coordinate systems, and this is handled as
outlined in more detail in~\cite{C1, Zattr}. Next, the computer
program is used to compute a verified upper bound on the
logarithmic norm of the Jacobian in block coordinates~\eqref{blockdcmp}.
As can be seen from the numerical data in Appendix~\ref{appendixnumdata},
this bound turns out to be negative. Notice that determining an upper
bound for a logarithmic norm is equivalent to computing the cone
condition~\eqref{cccomput} in the situation of only stable directions,
and our program makes use of the techniques described in more detail
in~\cite{C1, Zattr}. Thus, the numerical data presented in the
appendix ensures the existence of an attracting equilibriums
for~\eqref{eq:infDim2}, namely $\fss$ in~$\ws$, and $\fssg$ in~$\wsg$. 
In addition, we obtain that the
complete sets~$\ws$, $\wsg$ are contained in the basin of attraction of~$\fss$, 
and $\fssg$ respectively.
\begin{figure}[tb]
  \begin{center}
    \includegraphics[width=0.9\textwidth]{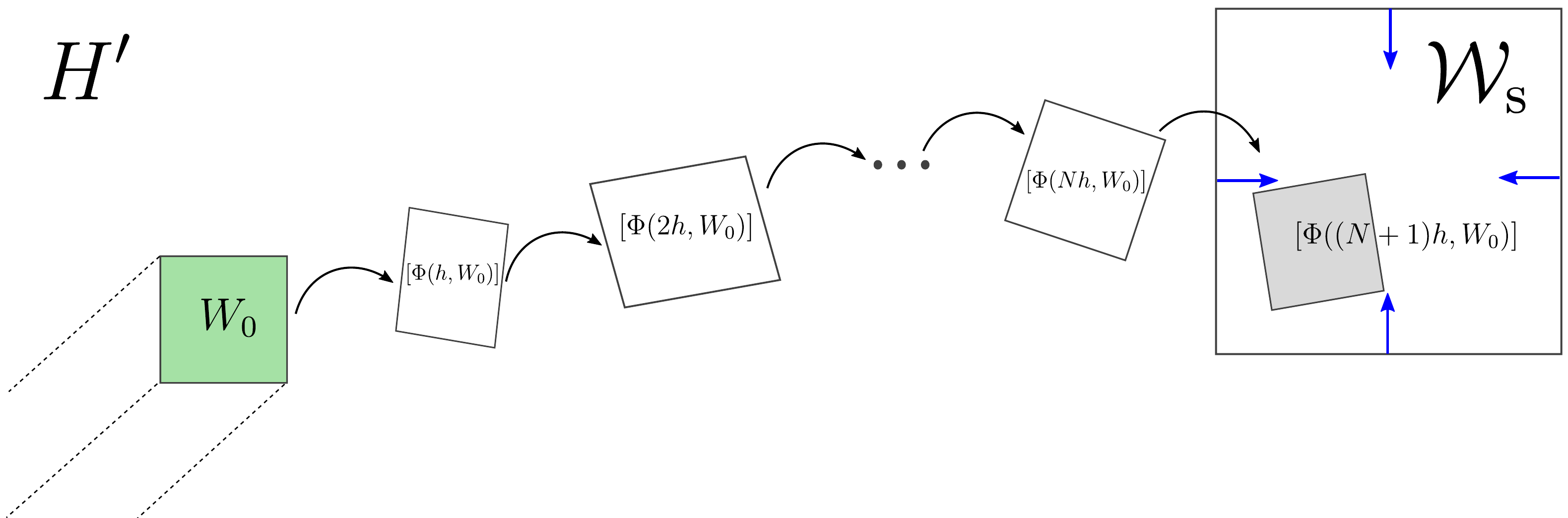}
  \end{center}
  \caption{The last part of the proof described in Step~3 establishes
           the existence of a heteroclinic connection. This is accomplished
           through rigorous integration forward in time of the face~$W_0$
           of~$\wu$, which was indicated in Figure~\ref{figublock}. Notice
           that this numerical integration has to be performed for the
           infinite-dimensional system~\eqref{eq:infDim2} with phase
           space~$\subspace$.}
\label{figintegration}
\end{figure}

{\bf Third step:} Since the first two steps of the proof
have established sets~$W_0(e_2), W_0(e_3)$ which contains part of the unstable 
manifold of~$\fsu$, as well as a subsets~$\ws$, $\wsg$ which are contained 
in the basin of attraction of~$\fss$, $\fssg$ respectively, it only remains to show that
all solutions of~\eqref{eq:infDim2} which originate in~$W_0(e_3)$ end
up in~$\ws$, and all solutions of~\eqref{eq:infDim2} which originate in~$W_0(e_2)$ 
end up in $\wsg$. This is accomplished in the last step of the proof
through rigorous forward integration of the flow associated
with~\eqref{eq:infDim2}. For this, we use the algorithm for
rigorous integration described in the series of works~\cite{C1,
KZ, ZPer, Z3}. In fact, we make use of an optimized version of
this algorithm which was developed in~\cite{Cy}, and this version
will be described in more detail in the next Section~\ref{secrigalg}.
In this way, we are able to rigorously propagate the sets~$W_0$
forward in time in such a way that all computational errors are
accounted for, and we obtain that the flow of~\eqref{eq:infDim2}
in~$\subspace$ transports $W_0(e_3)$ into~$\ws$, and $W_0(e_2)$ into~$\wsg$. 
A cartoon of this
procedure is shown in Figure~\ref{figintegration}. In the figure,
the sets~$[\Phi(kh,W_0)]$ denote rigorous enclosures for the 
images of~$W_0$ under the flow at time~$kh$, for integers $k > 0$.
After~$N+1$ steps of this forward integration by time~$h$, the
final enclosure is shown to be contained in~$\ws$, and therefore
the basin of attraction of $\fss$.
In the actual computation performed by a computer program 
any initial condition in the face $W_0(e_3)$ of $\wu$
 at time $T=3.02$ (after performing $1510$ numerical integration steps)
is verified to be within $\ws$, 
and any initial condition in the face $W_0(e_2)$ of $\wu$
at time $T=4.62$ (after performing $2310$ numerical integration steps)
is verified to be within $\wsg$.
As a consequence, one obtains the existence of a manifold
$W^u(\fsu) \cap W^s(\fss)$ which connects the unstable equilibrium~$\fsu$
with the stable equilibrium~$\fss$. As well as the existence of 
a manifold $W^u(\fsu) \cap W^s(\fssg)$ which connects the unstable equilibrium~$\fsu$
with the stable equilibrium~$\fssg$
For more details on the numerical
data from the actual computations we refer the reader to
Appendix~\ref{appendixnumdata} and the numerical data available in
an online repository \cite{codes}. Needless to say, the third and
last step of the proof is the most time consuming part.
\end{proof}

We would like to close this section with a brief comment on
the applicability of the above approach. 
The proof given above defines an algorithm for proving 
Theorem~\ref{thmmain}. The algorithm is general.
While in the proof we have considered a very specific
instance of the diblock copolymer model, for a particular set
of parameters, in principle 
the algorithm applies to the verification of any saddle to sink connection for the
diblock copolymer model~\eqref{eq:infDim2}.
The algorithm describes numerical verification  needed
to be performed
with a clear separation of the computation into three 
consecutive steps.
We have implemented a computer program, realizing this algorithm,
which is part of the publication \cite{codes} and available online.
 The input to the program
are the parameters in the equation, and if the program terminates
successfully we have proved the existence of a heteroclinic connection.
The main prerequisite for termination is of course that the 
set propagated by the rigorous numerical integrator is mapped
completely inside the stable isolating block~$\ws$. However,
the process of obtaining such a proof for arbitrary parameters
is not fully automated. In general, one has to try to make the
set~$\ws$ as large as possible, while at the same time keeping~$W_0$
small. In this way, the accumulated effect of roundoff errors
and inherent expansion in the differential equation can often
be controlled.
\section{Rigorous integration algorithm}
\label{secrigalg}
To rigorously integrate~\eqref{eq:infDim2} forward in time we apply
an algorithm based on the rigorous \emph{Taylor method}, as well as
the \emph{Lohner algorithm} for an efficient representation of sets
and to avoid the \emph{wrapping effect}. The Lohner algorithm is
presented in detail in~\cite{KZ, Lo}. We specifically use the algorithm
introduced in~\cite{Cy}, which allows for the efficient implementation of
the Taylor method for differential equations in combination with the
variational equations arising in partial differential equations. This
algorithm combines \emph{automatic differentiation} with the \emph{fast
Fourier transform algorithm}. In the remainder of this section, we briefly
explain the ideas behind this algorithm, and discuss why it is much more
efficient than the direct method of computation of normalized derivatives.
For a more detailed discussion, we refer the reader to~\cite{Cy}, where a
computational complexity analysis can be found in the context of partial 
differential equations with a quadratic nonlinearity. We will demonstrate
below that this analysis can be extended to the diblock copolymer
model~\eqref{dbcp1d} with a cubic nonlinearity. In fact, in the situation
of the present paper, the method of~\cite{Cy} leads to even greater
speedups. One of the main advantages of the method of~\cite{Cy} is that
it allows one to use a high-order Taylor method as opposed to the direct
integration method. In fact, in the actual proof of the results of this
paper we use a Taylor method of order~$16$.

To describe our approach in more detail, consider again the right-hand
side~$F$ of~\eqref{eq:infDim2}. In addition, we denote the right-hand
side of the $m$-th Galerkin approximation~\eqref{eq:galerkin} by~$F_m$.
We begin by giving a brief description of how the actual solution of
an initial value problem for~\eqref{eq:infDim2} can be enclosed in an
efficient way. This is achieved by bounding the solution set of the
\emph{differential inclusion} given by
\begin{equation} \label{difInc}
  \frac{d z}{d t}(t) \in F_m( z(t) ) + [\delta] \; ,
\end{equation}
where~$z(t)$ takes values in~$\mathbb{R}^m$. More precisely, the
solution set of~\eqref{difInc} is bounded by first enclosing the
solution of the \emph{Cauchy problem} associated with \eqref{difInc},
and then adding appropriate \emph{a-posteriori} bounds, see for
example~\cite{C1, KZ} for more details.

In order to rigorously enclose the solution of the Cauchy problem
associated with~\eqref{difInc}, we make use of the classical Taylor
method. Thus, the enclosure for the solution is propagated by one
time step through the knowledge of the \emph{normalized derivatives}
of orders up to~$p$, where the integer~$p$ denotes the order of the
Taylor method. These normalized derivatives are defined as follows.
\begin{definition}[Normalized derivative] \label{defnormder}
Let $n>0$ be an arbitrary integer, and let $w \colon
[0,t_{max}) \to \mathbb{R}$ denote a sufficiently smooth function.
Then the \emph{$n$-th normalized derivative} of~$w$ is defined
as
\begin{displaymath}
  w^{[n]}(t) = \frac{1}{n!} w^{(n)}(t) \; ,
\end{displaymath}
where~$w^{(n)}(t)$ denotes the classical $n$-th derivative
of~$w$.
\end{definition}
According to this definition, the normalized derivatives are
precisely the terms occurring in a standard Taylor expansion.
They can be determined efficiently via a well-known recursive
formula from \emph{automatic differentiation}, which is
recalled in the next remark.
\begin{rem}[Recursion for normalized derivatives]
For an ordinary differential equation of the form
\begin{displaymath}
  \frac{d z}{d t} = G(z) \; ,
  \quad\mbox{ where }\quad
  z \in \mathbb{R}^m \; ,
  \quad\mbox{ and }\quad
  G \colon \mathbb{R}^m \to \mathbb{R}^m \; ,
\end{displaymath}
the normalized derivatives of~$z$ can be determined as 
follows. If we set $z^{[0]} = z^{(0)}$, then for arbitrary
differentiation orders $j=1,\dots, q-1$ the normalized
derivative~$z^{[j]}$ is given by the recursive formula
\begin{displaymath}
  z^{[j+1]} = \frac{1}{j+1}(G\circ z)^{[j]} \; ,
\end{displaymath}
as long as~$G$ is $q$-times differentiable.
\end{rem}
After these preparations we can now analyze the cost of computing
the $(p+1)$-st order Taylor method for the finite system of ordinary
differential equations given by the Galerkin approximation~\eqref{eq:galerkin}.
If we assume that the normalized derivatives of~$u$ of orders $1,2,\dots,p$
have already been computed, then the cost of determining the $(p+1)$-st
normalized derivative is equivalent to the cost of evaluating the
convolutions
\begin{displaymath}
  (u \ast u \ast u)^{[p]}_k
  \qquad\mbox{ for all indices }\qquad
  |k| \leq m \; ,
\end{displaymath}
where the function~$u$ is given by its time-dependent Fourier coefficients,
and~$[p]$ denotes $p$-th normalized derivative with respect to time.

\providecommand{\uvar}{\mathbf{u}}
Since we are interested in using the Lohner algorithm~\cite{Lo}, we not
only have to keep track of the Fourier coefficients $\{a_k\}_{k=1}^m$ of the solution~$u$,
but also of the associated partial derivatives vector
\[
  \frac{\partial u_j}{\partial a_k},\quad j,k=1,\dots,m.
\]
for each $j=1,\dots,m$ the first order multi-variate polynomial with respect to
an auxiliary variable $\xi = \left(\xi_1,\dots,\xi_m\right)$
\[
\mathbf{u}_j(\xi) = \mathbf{u}_j = a_j + \sum_{k=1}^m{\frac{\partial u_j}{\partial a_k}\xi_k}
\]
is usually referred
to as the \emph{jet}. Thus, in order to emphasize that we work with
jets, not scalar values we will from
now on write~$\mathbf{u}$ instead of~$u$. An algebra of jets can be formed defining
addition and multiplication as the usual operations on polynomials. But
the terms of the second order in $\xi$ are truncated after the multiplication.
Our goal is the computation of the convolutions
\begin{displaymath}
  (\uvar \ast \uvar \ast \uvar)^{[p]}_k
  \qquad\mbox{ for all indices }\qquad
  |k| \leq m \; ,
\end{displaymath}
which can be expanded as
\begin{equation} \label{convolution}
  (\uvar \ast \uvar \ast \uvar)^{[p]}_k =
  \sum_{\substack{0 \leq p_1,p_2,p_3 \leq p\\p_1+p_2+p_3 = p}} \;\;
    \sum_{\substack{|k_1|,|k_2|,|k_3|\leq m\\k_1 + k_2 + k_3 = k}}
    {\uvar_{k_1}^{[p_1]}\uvar_{k_2}^{[p_2]}\uvar_{k_3}^{[p_3]}}
  \quad\mbox{ for all }\quad
  |k| \leq m \; .
\end{equation}
If we now let~$N$ denote the number of coefficients that are calculated
in~\eqref{convolution}, along with their partial derivatives for the
variational equations, then one can easily see that we basically have
\begin{displaymath}
  N = m^2+m \; ,
\end{displaymath}
up to a constant, which has to be large enough to avoid the 
\emph{aliasing effect} (in our algorithm the constant is equal to  $2$). 
With this, one can finally analyze the computational
cost of evaluating~\eqref{convolution}. While we do not present the 
details, the result is collected in the following remark.
\begin{rem}[Computational cost of evaluating the convolution]
\label{dircost}
The total cost of the direct evaluation of the convolution given
in~\eqref{convolution} is proportional to~$p^2N^3$, where~$N = m^2 + m$
up to a constant. An example analysis can be found in \cite{Cy}.
\end{rem}
Besides the above direct computation, it is well-known that one can
evaluate convolutions as in~\eqref{convolution} indirectly through
the use of the \emph{fast Fourier transform algorithm}. If we assume
again that all normalized derivatives of~$u$ of order $1,2,\dots,p$
have already been computed, then the discrete functional values of
the normalized derivatives of~$u$, which are also called the
\emph{$L$-coefficients} in~\cite{Cy}, have already been determined,
and they can therefore be cached in memory. We denote these values
by
\providecommand{\fft}[2][j]{\mbox{FFT}_{#1}\left(#2\right)}
\providecommand{\ffti}[2][k]{\mbox{FFT}^{-1}_{#1}\left(#2\right)}
\providecommand{\luvar}{\widehat{\uvar}}
\begin{displaymath}
  \luvar_j^{[r]} = \fft{\uvar^{[r]}}
  \quad\mbox{ for all }\quad
  r=0,\dots,p \; ,
  \quad j=0,\dots,c\cdot m \; ,
\end{displaymath}
where~$\fft{\cdot}$ denotes the $j$-th component of the discrete
Fourier transform computed by the fast Fourier transform method,
and the constant~$c$ is specified below. Using these values, one
can then determine~\eqref{convolution} by computing the convolution
for the $L$-coefficients in the form
\begin{equation} \label{fftconv}
  \left(\widehat{\uvar^3}\right)^{[p]}_j =
  \sum_{\substack{0\leq p_1,p_2,p_3\leq p\\p_1+p_2+p_3 = p}}
    {\luvar_j^{[p_1]}\ast \luvar_j^{[p_2]}\ast \luvar_j^{[p_3]}}
  \quad\mbox{ for all }\quad
  j=0,\dots,c\cdot m \; ,
\end{equation}
and then applying the inverse fast Fourier transform, which we denote
by~$\mbox{FFT}^{-1}$, in order to obtain
\begin{equation} \label{ffti}
  \left(\uvar^3\right)^{[p]}_k =
  \ffti{\left(\widehat{\uvar^3}\right)^{[p]}}
  \quad\mbox{ for all }\quad
  |k| \leq N \; .
\end{equation}
After this, we can also calculate the $L$-coefficients of~$\luvar^{[p+1]}$
and cache them in memory, provided we intend to increase the order further,
i.e., we evaluate
\begin{displaymath}
  \luvar^{[p+1]} = \mbox{FFT}\left(\uvar^{[p+1]}\right) \; .
\end{displaymath}
While at first glance the approach based on the fast Fourier transform
seems to be more involved, its computational savings are significant. To see
this, recall that the cost of performing a fast Fourier transform is of
the order~$N\log{N}$, and therefore the same is true for the computation
of~\eqref{ffti}. Moreover, the cost of computing the convolution for the
$L$-coefficients~\eqref{fftconv} is just the cost of multiplying~$p^2cN$
values, and this leads to the following observation.
\begin{rem}[Computational cost of the convolution via the fast Fourier transform]
\label{fftcostrem}
If we assume that the $L$-coefficients of~$\luvar$ have already been computed
and cached in memory, then the total cost of the fast Fourier transform approach
to the computation of~\eqref{convolution} is proportional to
\begin{equation} \label{fftcost}
  N\log{N} + p^2cN \; .
\end{equation}
An example analysis can be found in \cite{Cy}.
\end{rem}
Typically, the integer~$N$ is significantly larger than the order~$p$,
since it counts the number of equations and variationals.
Clearly the cost provided in Remark~\ref{dircost}, i.e. $p^2N^3$ dominates
both of the terms appearing in \eqref{fftcost} -- the cost provided in
Remark~\ref{fftcostrem} for any $p$ and $N$.
Also, increasing the order $p$ in case of the cost given in
Remark~\ref{fftcostrem} does not result in as severe increase of the
total computational cost, as in the case of Remark~\ref{dircost}.
The typical choice of the constant~$c$ in the above discussion is $c=3$.
Using some optimizations that were described in~\cite{Cy}, one can
in fact reduce the constant to~$c=2$.

Needless to say, in practice the fast Fourier transform technique
greatly reduces the total running time of our computer-assisted proof.
\begin{rem}
The actual running time of the integration step of the proof
with parameters as reported in  Appendix~\ref{appendixnumdata}
is \emph{triple} reduced. The real time of the 
integration procedure on a laptop is 3m55s FFT vs 11m55 direct (without FFT) 
in case of the connection to the local minimizer, and
6m34s FFT vs 17m29s in case of the connection to the global minimizer. 

We remark that the actual proof was performed
with a moderate Galerkin approximation dimension ($m=15$), 
which allows to achieve the proof in a reasonable time on a
laptop. In 
particular performing of an analogous proof, but for a parameter
regime further from zero will require larger approximations, 
due to much more complicated dynamics on the attractor.
The analysis 
of the algorithm presented in this section shows that for larger 
approximation dimension even a higher speedup is anticipated. 
\end{rem}
\section{Verifying the cone condition in infinite dimensions}
\label{secccverif}
This section is devoted to illustrating how the cone
condition~\eqref{cccomput} can be verified in infinite dimensions.
For the sake of brevity, this will only be done in the context of
the set~$\wu$. As mentioned earlier, the verification of the cone
condition for the remaining set~$\ws$ basically amounts to obtaining
logarithmic norm estimates, and this can be achieved directly with the 
algorithms from~\cite{C1, Zattr}, and is therefore omitted in the
following. To begin with, we would like to remind the reader that
according to~\eqref{cccomput}, we have to verify the estimate
\begin{equation} \label{ccrec}
  2\inf_{x \in \wu}{\left|\frac{\partial F_k}{\partial x_k}(x)\right|} -
    \sum_{\substack{\ell \in \mathbb{Z}\\ \ell \neq k}}
    {\sup_{x \in \wu}{\left|Q_{\ell\ell}\frac{\partial F_\ell}{\partial x_k}(x) +
    Q_{kk}\frac{\partial F_k}{\partial x_\ell}(x) \right|}}
  \geq \varepsilon \; ,  
\end{equation}
for the unstable isolating block~$\wu$, for some $\epsilon > 0$
and for all $k \in \mathbb{N}$.

In order to accomplish this, the block $\wu \subset \subspace$
is divided into two parts, and we refer to this representation again
as self-consistent bounds. More precisely, we have that
\begin{itemize}
\item the \emph{finite part} for $k=1,\dots,M$ is represented
as a finite-dimensional interval product, i.e., the
projection~$\left(\wu\right)_k$ is an interval, and
\item the \emph{tail}, which corresponds to $k = M+1, \dots$,
is represented by an algebraically  decaying bound, i.e., for
all $x \in \left(\wu\right)_k$ the inequality $|x| \leq C/|k|^6$
is satisfied. In order to explicitly relate the constant~$C$ to
the block~$\wu$, we will use the notation~$C(\wu)$ in the
following.
\end{itemize}
Notice that $\wu$ forms self-consistent bounds with the fixed
polynomial decay, therefore the
tail in the algorithm is represented by only the value of the
constant~$C > 0$.
\providecommand{\conv}{\mathcal{C}}

In order to establish the cone condition, we first determine
self-consistent bounds~$\mathcal{C} \subset \subspace$, which
bound all possible convolutions of elements from the block~$\wu$,
i.e., such that for all $k = 0,\dots,M$ we have
\begin{displaymath}
  \conv_k
  \quad\mbox{ is such that }\quad
  \conv_k \supset (a\ast b)_k
  \quad\mbox{ for all }\quad
  a, b \in \wu \; .
\end{displaymath}
Since the nonlinearity in~\eqref{eq:infDim2} is cubic, one therefore
immediately obtains an enclosure for the partial derivatives matrix~$DF$
in the form
\begin{displaymath}
  \left(D F\right)_{k\ell} =
  \frac{\partial F_k}{\partial x_\ell}(x) \subset
  \delta_{k\ell} \left( -\frac{k^4\pi^4}{L^4} + \frac{\lambda k^2\pi^2}{L^2} -
    \lambda\sigma \right) -
    \frac{3 \lambda k^2 \pi^2}{L^2} \left(\conv_{k-\ell} +
      \conv_{k+\ell} \right)
  \quad\mbox{ for all }\quad
  x \in \wu \; .
\end{displaymath}
Before proceeding, we need to introduce some notation. if~$I \subset
\mathbb{R}$ denotes an arbitrary interval, then we denote the maximal
absolute value of all points from~$I$ by~$|I|_{\max}$. We can now estimate
the second expression in~(\ref{ccrec}) by
\begin{displaymath}
  \sum_{\substack{\ell\in\mathbb{Z}\\ \ell \neq k}}
    {\sup_{x \in \wu}{\left|Q_{\ell\ell}\frac{\partial F_\ell}{\partial x_k}(x) +
    Q_{kk}\frac{\partial F_k}{\partial x_\ell}(x) \right|}} \leq
  \frac{3\lambda \pi^2}{L^2} \sum_{\substack{\ell\in\mathbb{Z}\\ \ell \neq k}}
    {\left| \ell^2 Q_{\ell\ell} \left( \conv_{\ell-k} + \conv_{k+\ell}\right) +
    k^2 Q_{kk} \left(\conv_{k-\ell} + \conv_{k+\ell} \right) \right|_{\max}}
  \; .
\end{displaymath}
The right-hand side in this estimate can be bounded further.
For this, we split it into two parts to obtain
\begin{multline*}
  \frac{3\lambda \pi^2}{L^2} \sum_{\substack{\ell\in\mathbb{Z}\\ \ell\neq k}}
    {\left| \ell^2 Q_{\ell\ell} \left(\conv_{\ell-k} + \conv_{k+\ell}\right) +
    k^2 Q_{kk} \left(\conv_{k-\ell} + \conv_{k+\ell}\right) \right|_{\max}} = \\
  \frac{3\lambda \pi^2}{L^2} \sum_{\substack{|\ell| \leq M+k\\ \ell\neq k}}
    {\left| \ell^2 Q_{\ell\ell} \left(\conv_{\ell-k} + \conv_{k+\ell}\right) +
    k^2 Q_{kk} \left(\conv_{k-\ell} + \conv_{k+\ell}\right)\right|_{\max}} + \\
  \frac{3\lambda\pi^2}{L^2} \sum_{\substack{|\ell| > M+k\\ \ell\neq k}}
    {\left| \ell^2 Q_{\ell\ell} \left(\conv_{\ell-k} + \conv_{k+\ell}\right) +
    k^2 Q_{kk} \left(\conv_{k-\ell} + \conv_{k+\ell}\right)\right|_{\max}}
    \; .
\end{multline*}
The first part of this sum consists of expressions which contain at least
one term from the finite part. In addition, the first summand is a finite 
sum, and we can therefore compute it term by term using interval arithmetic.
In contrast, the second summand is infinite and consists exclusively of
terms from the tail part, which can therefore be estimated. Using standard
results from rigorous numerics for partial differential equations,
due to a high regularity of the Fourier coefficients sequences all
of the appearing sums have a finite value, which can be estimated efficiently.

Finally this yields the following estimate
\begin{multline*}
  \frac{3\lambda \pi^2}{L^2} \sum_{\substack{|\ell| > M+k\\ \ell\neq k}}
    {\left| \ell^2 Q_{\ell\ell} \left(\conv_{\ell-k} + \conv_{k+\ell}\right) +
    k^2 Q_{kk} \left(\conv_{k-\ell} + \conv_{k+\ell}\right)\right|_{\max}} \leq \\
  \frac{3\lambda \pi^2}{L^2} \sum_{\substack{|\ell| > M+k\\ \ell\neq k}}
    {\left|\ell^2 \left(\conv_{\ell-k} + \conv_{k+\ell}\right)\right|_{\max} +
    \frac{3\lambda \pi^2}{L^2} \sum_{\substack{|\ell| > M+k\\ \ell\neq k}}
    \left| k^2 \left(\conv_{k-\ell} + \conv_{k+\ell}\right)\right|_{\max}} \leq \\
  \frac{3\lambda \pi^2}{L^2} \left( \frac{2C(\conv) k^2}{5M^{5}} +
    \frac{2C(\conv)}{3M^{3}}(1 + k / M) \right) \; ,
\end{multline*}
whose right-hand side is given explicitly --- and from this one can easily
establish the validity of~(\ref{ccrec}) for $k = 1,\dots, M$.

For $k>M$ using now standard arguments, see e.g. \cite{Z3},
it can be shown that (\ref{ccrec}) is decreasing with $k$
for $k>M$, due to the fast decay of eigenvalues. Hence,
the negativity of (\ref{ccrec}) for $k>M$ follows from the
negativity of (\ref{ccrec}) for $k=M$.

Finally the cone condition can be verified.
\section{Software}
\label{secsoftware}
The computer program performing all computational steps required to prove
Theorem~\ref{thmmain}, as described in more detail in the respective proof
in Section~\ref{secproof}, has been written in the C++ programming language. 
The source codes, some compiled binary files, and the numerical output from
the actual proofs are available in an online bitbucket repository~\cite{codes}.
The published source codes for this work utilizes some parts of the rigorous
numerics software package developed by the first author, which was dedicated
to the study of dissipative partial differential equations based on Fourier
series expansions. This package has been used in a number of different
contexts as well, see for example~\cite{Cy, C1, CZ}.

The software package is written in such a way that it should be easily
adaptable to other situations as well. More precisely, by performing only
small modifications in the code, it can be applied to a range of dissipative
partial differential equations which are written in terms of infinite
Fourier series. Examples include the viscous Burgers equation, the 
Kuramoto-Sivashinsky equation, and the Swift-Hohenberg equation.
The main part of the software package is a dedicated algorithm for
performing time integration using Taylor's method, and the efficient
computation of the (higher order) time derivatives of the flow. The
latter are required by Taylor's method, and they are determined using
the FFT algorithm described earlier.
\section{Conclusion and future work}
\label{secconclusions}
In this paper, we have presented a first instance in which heterolinic
solutions can be rigorously determined for the diblock copolymer model
on one-dimensional domains. Moreover, our results have shown that the
model exhibits strong multistability in the sense that all (local or
global) energy minimizers can be reached from arbitrarily small 
neigborhoods of the homogeneous steady state.
In the near future we plan to extend the
presented research to prove heteroclinic connections for a wide range
of other parameter regimes, with particular emphasis on regions further
from the homogeneous steady state. Moreover, we plan to perform similar
studies on simple two-dimensional domains. Another intriguing possibility
for future work is to apply the parametrization method~\cite{CFL} in the
first step of the proof of Theorem~\ref{thmmain}. Recently, the
parametrization method has been generalized to allow for the validation
of unstable manifolds for parabolic partial differential equations,
see for example~\cite{RMJ, vandenBerg2016}. Using the parametrization
method, bounds for the unstable manifold are obtained by computing its
expansion in the form of a finite Taylor polynomial. In contrast, in
the present work, we obtain bounds for the unstable manifold by verifying
the cone condition on a linear segment. For nonlinear partial differential
equations and nontrivial parameter values, it is expected that the
unstable manifold behaves linearly only very close to the unstable
steady state. We are convinced that, using the parametrization method,
especially for high Taylor orders, one can compute bounds for the unstable 
manifold of the diblock copolymer model that reach much further than
those obtained using a linear block and cone conditions. This in turn
is expected to greatly reduce the work needed for performing the final
step of the proof, i.e., using the rigorous integrator to demonstrate that
the unstable manifold enters the basin of attraction of the stable steady
state. This should make it possible to then establish a result analogous
to Theorem~\ref{thmmain}, but for large parameter values.
\section{Acknowledgments}
The presented work was performed while the first author held
post-doctoral positions at the Warsaw Center of Mathematics and
Computer Science, Poland, and subsequently at Rutgers, The State
University of New Jersey, USA. The idea for the study arose during
the workshop ``Rigorous computation for infinite-dimensional nonlinear
dynamics'' at the American Institute of Mathematics in Palo Alto,
California, in August~2013. J. Cyranka was supported in part by 
NSF grant DMS 1125174 and DARPA contract HR0011-16-2-0033.
The work of the second author was partially
supported by NSF grants DMS-1114923 and DMS-1407087.

\appendix
\section{Numerical data from the computer proof}
\label{appdix}
\subsection{Heteroclinic to the local minimizer}
\label{appendixnumdata}
In this appendix we present some numerical data from the
computer-assisted part of the proof presented in this paper.
For the sake of brevity we list only a small number of the
relevant modes which are involved in the constructed
infinite-dimensional bounds. In addition, only three
decimals of all numerical values are shown.
Complete and detailed numerical data presented with larger
precision can be found online \href{http://bitbucket.org/dzako/dbcp\_proof}{http://bitbucket.org/dzako/dbcp\_proof}.

\paragraph{Unstable isolating block~$\mathcal{W}_u\subset \subspace$
around zero.} This is the first step of the proof of the main theorem
from Section~\ref{secproof}. We mark with text in bold the 
numerical values for the eigendirection $e_3$ along which the unstable manifold of the  homogeneous state 
connects to the stable \emph{local} minimizer ($\mathcal{W}_u$ is stretched along this direction)..

\[
\footnotesize{
\begin{array}{|c|c|}\hline\mathbf{k} & \mathbf{a_k}\text{ interval in the form center + radius} \\ \hline\hline
1 & 0+[-2.365,2.365]10^{-16}\\ 
2 & 0+[-1,1]10^{-12}\\ 
\mathbf{3} & \mathbf{[-0.075,0.075]}\\ 
4 & 0+[-1.538,1.538]10^{-15}\\ 
5 & 0+[-3.038,3.038]10^{-16}\\ 
6 & 0+[-1.478,1.478]10^{-16}\\ 
7 & 0+[-9.664,9.664]10^{-17}\\ 
8 & 0+[-5.623,5.623]10^{-17}\\ 
9 & 0+[-5.193,5.193]10^{-17}\\ 
10 & 0+[-3.771,3.771]10^{-17}\\ 
11 & 0+[-3.332,3.332]10^{-17}\\ 
12 & 0+[-2.248,2.248]10^{-17}\\ 
13 & 0+[-2.33,2.33]10^{-17}\\ 
14 & 0+[-1.926,1.926]10^{-17}\\ 
15 & 0+[-1.69,1.69]10^{-17}\\ \hline
16-75&\text{small intervals of width }\leq 10^{-20}\\ \hline
\geq 76 & <1.393\cdot 10^{-46}/k^{6}\\\hline\end{array}}
\]
The infinite-dimensional cone condition \eqref{ccrec} is satisfied with
\[
\varepsilon = 0.08478.
\]
The file \emph{unstable\_box\_log.txt} which can be found online \cite{codes}
contains more detailed numerical data related to this computation.
The actual running time of the program computing this set is negligible
compared to the last step, i.e., the numerical integration. Therefore,
we do not report it.

\paragraph{Stable isolating block $\ws\subset \subspace$ }
This is the second step of the proof of the main theorem from
Section~\ref{secproof}. The block~$\ws$ is centered at a numerically
computed approximation to the actual stable equilibrium~$\fss$ which
was loaded from the provided file \emph{fixedPoint.in}.

\[
\footnotesize{
\begin{matrix}
\begin{array}{|c|c|}\hline\mathbf{k} & \mathbf{a_k} \\ \hline\hline
1 & 3.321\cdot 10^{-18}+[-6.895,6.895]10^{-6}\\ 
2 & 8.18\cdot 10^{-17}+[-1.691,1.691]10^{-4}\\ 
\mathbf{3} & \mathbf{0.2956+[-1.031,1.031]10^{-5}}\\ 
4 & -6.551\cdot 10^{-17}+[-9.807,9.807]10^{-5}\\ 
5 & -1.241\cdot 10^{-17}+[-9.732,9.732]10^{-6}\\ 
6 & -1.411\cdot 10^{-16}+[-2.994,2.994]10^{-6}\\ 
7 & -3.635\cdot 10^{-18}+[-2.875,2.875]10^{-6}\\ 
8 & -9.215\cdot 10^{-18}+[-1.417,1.417]10^{-5}\\ 
9 & -5.613\cdot 10^{-3}+[-1.518,1.518]10^{-6}\\ 
10 & -2.005\cdot 10^{-18}+[-6.141,6.141]10^{-6}\\ 
11 & -4.048\cdot 10^{-18}+[-1.659,1.659]10^{-6}\\ 
12 & 4.056\cdot 10^{-16}+[-1.001,1.001]10^{-6}\\ 
13 & -2.184\cdot 10^{-18}+[-1.367,1.367]10^{-6}\\ 
14 & -1.914\cdot 10^{-18}+[-1.559,1.559]10^{-6}\\ 
15 & 1.062\cdot 10^{-4}+[-6.264,6.264]10^{-7}\\ 
16 & -1.771\cdot 10^{-19}+[-7.908,7.908]10^{-7}\\ 
17 & -1.616\cdot 10^{-19}+[-5.853,5.853]10^{-7}\\ 
18 & 1.351\cdot 10^{-16}+[-2.593,2.593]10^{-7}\\ 
19 & -2.224\cdot 10^{-19}+[-3.73,3.73]10^{-7}\\ 
20 & -1.555\cdot 10^{-19}+[-2.83,2.83]10^{-7}\\ 
21 & -2.021\cdot 10^{-6}+[-2.355,2.355]10^{-7}\\ 
22 & 6.774\cdot 10^{-20}+[-3.8,3.8]10^{-7}\\ 
23 & 6.011\cdot 10^{-20}+[-3.309,3.309]10^{-7}\\ \hline
\end{array}&
\begin{array}{|c|c|} \hline
24 & -1.441\cdot 10^{-17}+[-1.763,1.763]10^{-7}\\ 
25 & 2.333\cdot 10^{-20}+[-1.87,1.87]10^{-7}\\ 
26 & 2.084\cdot 10^{-20}+[-1.651,1.651]10^{-7}\\ 
27 & 3.847\cdot 10^{-8}+[-7.323,7.323]10^{-8}\\ 
28 & 1.265\cdot 10^{-20}+[-1.417,1.417]10^{-7}\\ 
29 & 1.257\cdot 10^{-20}+[-1.366,1.366]10^{-7}\\ 
30 & 3.047\cdot 10^{-16}+[-1.1,1.1]10^{-7}\\ 
31 & 7.491\cdot 10^{-19}+[-1.754,1.754]10^{-7}\\ 
32 & 6.881\cdot 10^{-19}+[-1.7,1.7]10^{-7}\\ 
33 & -7.323\cdot 10^{-10}+[-8.918,8.918]10^{-8}\\ 
34 & 7.818\cdot 10^{-19}+[-8.435,8.435]10^{-8}\\ 
35 & 7.079\cdot 10^{-19}+[-1,1]10^{-7}\\ 
36 & 2.828\cdot 10^{-16}+[-5.943,5.943]10^{-8}\\ 
37 & 7.42\cdot 10^{-19}+[-1.151,1.151]10^{-7}\\ 
38 & 6.551\cdot 10^{-19}+[-1.099,1.099]10^{-7}\\ 
39 & 1.394\cdot 10^{-11}+[-8.788,8.788]10^{-8}\\
40 & -1.015\cdot 10^{-12}+[-2.472,2.472]10^{-10}\\ 
41 & -1.129\cdot 10^{-12}+[-2.747,2.747]10^{-10}\\ 
42 & -6.509\cdot 10^{-13}+[-1.584,1.584]10^{-10}\\ 
43 & -1.184\cdot 10^{-12}+[-2.881,2.881]10^{-10}\\ 
44 & -1.08\cdot 10^{-12}+[-2.627,2.627]10^{-10}\\ 
45 & -1.075\cdot 10^{-12}+[-1.967,1.967]10^{-10}\\ \hline
46-250&\text{small intervals of width }\leq 10^{-11}\\ \hline
\geq 251 & <6.947\cdot 10^{-38}/k^{6}\\\hline\end{array}
\end{matrix}
}\]
The infinite-dimensional cone condition \eqref{ccrec} is satisfied with
\[
\varepsilon = 0.03816 \; ,
\]
which equivalently means that the logarithmic norm estimate is strictly
negative. The interval radii of the first~$39$ modes, separated above
by a horizontal line, of the block~$\ws$ are given in the
eigenbasis of the finite-dimensional derivative at~$\fss$. These
coordinates are being written by the program into the file~\emph{Q.in}.
The actual running time of the program computing this set is again
negligible compared to the last step, i.e., the numerical integration.
As before, we therefore do not report it.

The file \emph{stable\_box\_log.txt} that can be found online \cite{codes}
contains more detailed numerical data related with this computation.

\paragraph{Result of the numerical integration.}
This is the third and final step of the proof of the main theorem from
Section~\ref{secproof}. Some technical parameters of the numerical
integration are chosen in the following way:
\begin{itemize}
\item The Taylor method order is~$16$, 
\item the constant time step is~$0.002$,
\item the Galerkin projection dimension~$m$ is~$15$, the dimension
used by the FFT-algorithm~$N$ is~$32$, see also the description in
Section~\ref{secrigalg}.
\item The whole integration process took \textbf{3m55s} on an
Intel\textregistered Core\texttrademark i7-4610M CPU @ 3.00GHz x 4,
compiled with gcc version 4.9.2 (64-bit Ubuntu 4.9.2-10ubuntu13)
and with flags -O2 -frounding-math.
\item In contrast, if instead of the FFT-algorithm we used the
direct algorithm, then the whole integration process took \textbf{11m55s}.
\end{itemize}
We note that the overall speedup of the FFT algorithm over the direct algorithm
does not seem to be overwhelming in this particular example, we guess
it is due to the fact of moderate Galerkin approximation dimension ($15$)
used. We also guess that some parts of our code may not be optimized,
which results in large constants multiplying the total cost
of the FFT algorithm \eqref{fftcost}. We expect much larger speedups
in case of larger Galerkin approximation dimensions. We will
investigate this as a future research.

The file \emph{num\_integration\_log.txt} that can be found
online~\cite{codes} contains more detailed numerical data related
to this computation. Specifically, the file contains
\begin{itemize}
\item the vector which is the finite-dimensional part of the
propagated polynomial bounds, saved at each step of the integration,
\item each~$100$ steps of the integration the full
infinite-dimensional polynomial bounds,
\item at each step of the integration information at which
coordinates the entry of~$\mathcal{W}_0$ into~$\ws$ 
has not yet been attained.
\end{itemize}
The infinite-dimensional bounds at time $T=3.02$ after
performing~$1510$ numerical integration steps are
\[
\begin{matrix}
\footnotesize{
\begin{array}{|c|c|}\hline\mathbf{k} & \mathbf{a_k} \\ \hline\hline
1 & -1.014\cdot 10^{-8}+[-1.49,1.49]10^{-7}\\ 
2 & 1.454\cdot 10^{-7}+[-2.748,2.748]10^{-7}\\ 
\mathbf{3} & \mathbf{0.2956+[-8.522,8.522]10^{-9}}\\ 
4 & -1.293\cdot 10^{-7}+[-1.183,1.183]10^{-7}\\ 
5 & -7.29\cdot 10^{-9}+[-3.619,3.619]10^{-8}\\ 
6 & -1.589\cdot 10^{-11}+[-1.156,1.156]10^{-9}\\ 
7 & -3.299\cdot 10^{-8}+[-1.501,1.501]10^{-8}\\ 
8 & 1.725\cdot 10^{-8}+[-2.031,2.031]10^{-8}\\ 
9 & -5.613\cdot 10^{-3}+[-1.14,1.14]10^{-9}\\ 
10 & 7.704\cdot 10^{-9}+[-8.354,8.354]10^{-9}\\ 
11 & -7.297\cdot 10^{-10}+[-4.831,4.831]10^{-9}\\ 
12 & 1.083\cdot 10^{-8}+[-4.586,4.586]10^{-9}\\ 
13 & 5.124\cdot 10^{-9}+[-7.282,7.282]10^{-10}\\ 
14 & 1.396\cdot 10^{-9}+[-9.208,9.208]10^{-10}\\ 
15 & 1.063\cdot 10^{-4}+[-4.66,4.66]10^{-9}\\ 
16 & -1.778\cdot 10^{-10}+[-2.702,2.702]10^{-10}\\ 
17 & 6.718\cdot 10^{-12}+[-8.894,8.894]10^{-11}\\ 
18 & -8.661\cdot 10^{-11}+[-1.11,1.11]10^{-10}\\ 
19 & -4.31\cdot 10^{-11}+[-5.222,5.222]10^{-11}\\ 
20 & -5.023\cdot 10^{-12}+[-3.588,3.588]10^{-11}\\ 
21 & -2.021\cdot 10^{-6}+[-4.376,4.376]10^{-10}\\ 
22 & 4.655\cdot 10^{-12}+[-7.155,7.155]10^{-12}\\ 
23 & -2.128\cdot 10^{-13}+[-2.362,2.362]10^{-12}\\ \hline
\end{array}}&
\footnotesize{\begin{array}{|c|c|} \hline
24 & 2.481\cdot 10^{-12}+[-3.129,3.129]10^{-12}\\ 
25 & 1.17\cdot 10^{-12}+[-1.405,1.405]10^{-12}\\ 
26 & 1.492\cdot 10^{-13}+[-9.22,9.22]10^{-13}\\ 
27 & 3.847\cdot 10^{-8}+[-1.241,1.241]10^{-11}\\ 
28 & -1.128\cdot 10^{-13}+[-1.779,1.779]10^{-13}\\ 
29 & 5.294\cdot 10^{-15}+[-5.796,5.796]10^{-14}\\ 
30 & -6.041\cdot 10^{-14}+[-7.655,7.655]10^{-14}\\ 
31 & -2.807\cdot 10^{-14}+[-3.384,3.384]10^{-14}\\ 
32 & -3.631\cdot 10^{-15}+[-2.199,2.199]10^{-14}\\ 
33 & -7.325\cdot 10^{-10}+[-3.017,3.017]10^{-13}\\ 
34 & 2.608\cdot 10^{-15}+[-4.228,4.228]10^{-15}\\ 
35 & -1.224\cdot 10^{-16}+[-1.359,1.359]10^{-15}\\ 
36 & 1.389\cdot 10^{-15}+[-1.775,1.775]10^{-15}\\ 
37 & 6.412\cdot 10^{-16}+[-7.772,7.772]10^{-16}\\ 
38 & 8.304\cdot 10^{-17}+[-5.039,5.039]10^{-16}\\ 
39 & 1.395\cdot 10^{-11}+[-6.956,6.956]10^{-15}\\ 
40 & -5.841\cdot 10^{-17}+[-9.691,9.691]10^{-17}\\ 
41 & 2.718\cdot 10^{-18}+[-3.085,3.085]10^{-17}\\ 
42 & -3.09\cdot 10^{-17}+[-3.994,3.994]10^{-17}\\ 
43 & -1.421\cdot 10^{-17}+[-1.738,1.738]10^{-17}\\ 
44 & -1.839\cdot 10^{-18}+[-1.126,1.126]10^{-17}\\ 
45 & -2.655\cdot 10^{-13}+[-1.557,1.557]10^{-16}\\ \hline
46-200 & \text{small intervals of width }\leq 10^{-17}\\\hline
\geq 201 & <3.295\cdot 10^{-44}/k^{6}\\\hline\end{array}
}\end{matrix}\]
The bounds above are within the set~$\ws$, i.e.,
within the basin of attraction of the stable fixed point.
This implies that the third step of the proof of the main
theorem from Section~\ref{secproof} succeeded.

\subsection{Heteroclinic to the global minimizer}
We present the analogous numerical data from the computer-assisted
existence proof of the heteroclinic between the homogeneous state 
and the \emph{global} minimizer.
\paragraph{Unstable isolating block $\mathcal{W}_u\subset \subspace$ around zero}
This is the first step of the proof of the main theorem from Section~\ref{secproof}). We mark with text in bold the 
numerical values for the eigendirection $e_2$ along which the unstable manifold of the homogeneous state 
connects to the stable \emph{global} minimizer ($\mathcal{W}_u$ is stretched along this direction).
\[\footnotesize{
\begin{array}{|c|c|}\hline\mathbf{k} & \mathbf{a_k}\text{ interval in the form center + radius} \\ \hline\hline
1 & 0+[-1.105,1.105]10^{-15}\\ 
\mathbf{2} & \mathbf{[-0.1,0.1]}\\ 
3 & 0+[-1,1]10^{-12}\\ 
4 & 0+[-7.235,7.235]10^{-15}\\ 
5 & 0+[-1.438,1.438]10^{-15}\\ 
6 & 0+[-6.918,6.918]10^{-16}\\ 
7 & 0+[-4.579,4.579]10^{-16}\\ 
8 & 0+[-2.626,2.626]10^{-16}\\ 
9 & 0+[-2.449,2.449]10^{-16}\\ 
10 & 0+[-1.742,1.742]10^{-16}\\ 
11 & 0+[-1.571,1.571]10^{-16}\\ 
12 & 0+[-1.054,1.054]10^{-16}\\ 
13 & 0+[-1.102,1.102]10^{-16}\\ 
14 & 0+[-8.963,8.963]10^{-17}\\ 
15 & 0+[-8.008,8.008]10^{-17}\\ \hline
16-75&\text{small intervals of width }\leq 10^{-19}\\ \hline
\geq 76 & <6.925\cdot 10^{-59}/k^{6}\\\hline\end{array}
}\]
The infinite-dimensional cone condition \eqref{ccrec} is satisfied with
\[
\varepsilon = 0.2873.
\]
The file \emph{unstable\_box\_log.txt} which can be found online \cite{codes}
contains more detailed numerical data related with this computation.

\paragraph{Stable isolating block $\wsg\subset \subspace$ }
(The second step  of the proof of the main theorem from Section~\ref{secproof})
$\wsg$ is centered at a numerically computed approximation to the actual global stable equilibrium
loaded from the provided file \emph{fixedPointGlobal.in}.

Observe that although the dimension of the finite-dimensional subtail here is larger ($300$),
where by finite-dimensional subtail we mean the distinguished part of the whole tail $T$,
which is represented by explicit intervals.
The diameters of $\wsg$ are now considerably larger, as well as 
the $\varepsilon$ cone condition value, than in the case of the local minimizer.
This is due to the fact that the global minimizer is stronger attracting than the local 
minimizer in the sense that the maximal eigenvalue of the Jacobian matrix at the 
global minimizer is further away from zero than the maximal eigenvalue at the
local minimizer.
\[\begin{matrix}\footnotesize{
\begin{array}{|c|c|}\hline\mathbf{k} & \mathbf{a_k} \\ \hline\hline
1 & -6.339\cdot 10^{-18}+[-2.122,2.122]10^{-5}\\ 
\mathbf{2} & \mathbf{0.3484+[-2.368,2.368]10^{-5}}\\ 
3 & 3.969\cdot 10^{-17}+[-4.191,4.191]10^{-5}\\ 
4 & -1.125\cdot 10^{-15}+[-1.129,1.129]10^{-5}\\ 
5 & 3.558\cdot 10^{-17}+[-1.213,1.213]10^{-5}\\ 
6 & -0.02108+[-7.626,7.626]10^{-6}\\ 
7 & -9.56\cdot 10^{-18}+[-9.823,9.823]10^{-6}\\ 
8 & 4.792\cdot 10^{-16}+[-4.311,4.311]10^{-6}\\ 
9 & -2.52\cdot 10^{-18}+[-4.51,4.51]10^{-6}\\ 
10 & 1.207\cdot 10^{-3}+[-2.687,2.687]10^{-6}\\ 
11 & 4.678\cdot 10^{-18}+[-2.501,2.501]10^{-6}\\ 
12 & -1.251\cdot 10^{-15}+[-1.19,1.19]10^{-6}\\ 
13 & 7.395\cdot 10^{-18}+[-1.571,1.571]10^{-6}\\ 
14 & -7.109\cdot 10^{-5}+[-1.238,1.238]10^{-6}\\ 
15 & 5.117\cdot 10^{-18}+[-1.442,1.442]10^{-6}\\ 
16 & -4.127\cdot 10^{-16}+[-9.132,9.132]10^{-7}\\ 
17 & 2.433\cdot 10^{-18}+[-8.563,8.563]10^{-7}\\ 
18 & 4.197\cdot 10^{-6}+[-4.648,4.648]10^{-7}\\ 
19 & -1.265\cdot 10^{-18}+[-6.274,6.274]10^{-7}\\ 
20 & 4.897\cdot 10^{-17}+[-5.667,5.667]10^{-7}\\ 
21 & -2.285\cdot 10^{-19}+[-7.082,7.082]10^{-7}\\ 
22 & -2.479\cdot 10^{-7}+[-4.735,4.735]10^{-7}\\ 
23 & 1.889\cdot 10^{-18}+[-4.624,4.624]10^{-7}\\ 
24 & -2.62\cdot 10^{-16}+[-2.606,2.606]10^{-7}\\ \hline
\end{array} }&
\footnotesize{
\begin{array}{|c|c|}\hline
25 & 2.058\cdot 10^{-18}+[-2.564,2.564]10^{-7}\\ 
26 & 1.465\cdot 10^{-8}+[-2.352,2.352]10^{-7}\\ 
27 & 3.304\cdot 10^{-19}+[-2.936,2.936]10^{-7}\\ 
28 & -3.412\cdot 10^{-17}+[-1.982,1.982]10^{-7}\\ 
29 & -2.625\cdot 10^{-19}+[-1.998,1.998]10^{-7}\\ 
30 & -8.656\cdot 10^{-10}+[-1.384,1.384]10^{-7}\\ 
31 & 1.09\cdot 10^{-19}+[-2.202,2.202]10^{-7}\\ 
32 & -1.348\cdot 10^{-17}+[-2.093,2.093]10^{-7}\\ 
33 & 1.483\cdot 10^{-19}+[-2.814,2.814]10^{-7}\\ 
34 & 5.115\cdot 10^{-11}+[-2.029,2.029]10^{-7}\\ 
35 & 2.292\cdot 10^{-19}+[-1.962,1.962]10^{-7}\\ 
36 & -4.436\cdot 10^{-17}+[-1.683,1.683]10^{-7}\\ 
37 & 2.326\cdot 10^{-19}+[-2.629,2.629]10^{-7}\\ 
38 & -3.022\cdot 10^{-12}+[-2.325,2.325]10^{-7}\\ 
39 & 3.417\cdot 10^{-20}+[-2.932,2.932]10^{-7}\\ 
40 & -3.135\cdot 10^{-12}+[-7.641,7.641]10^{-10}\\ 
41 & -4.553\cdot 10^{-12}+[-1.109,1.109]10^{-9}\\ 
42 & -3.567\cdot 10^{-12}+[-9.137,9.137]10^{-10}\\ 
43 & -4.418\cdot 10^{-12}+[-1.076,1.076]10^{-9}\\ 
44 & -1.845\cdot 10^{-13}+[-7.486,7.486]10^{-11}\\ 
45 & -3.045\cdot 10^{-13}+[-1.092,1.092]10^{-10}\\ 
46 & -2.93\cdot 10^{-13}+[-9.065,9.065]10^{-11}\\ 
47 & -3.638\cdot 10^{-13}+[-1.075,1.075]10^{-10}\\ \hline
48-300&\text{small intervals of width }10^{-12}\\ \hline
\geq 301 & <5.334\cdot 10^{-44}/k^{6}\\\hline\end{array} }
\end{matrix}
\]
The infinite-dimensional cone condition \eqref{ccrec} is satisfied with
\[
\varepsilon = 0.7869,
\]
which equivalently means that the logarithmic norm estimate is negative.

The file \emph{stable\_box\_log.txt} that can be found online \cite{codes}
contains more detailed numerical data related with this computation.
As in the case of the local minimizer, the interval radii of the first~$39$ modes (separated above with a horizontal line) 
of $\wsg$ are
given in the eigenbasis of the finite-dimensional derivative at $\fssg$.

\paragraph{Result of the numerical  integration}
This is the third and final step of the proof of the main theorem from
Section~\ref{secproof}. Some technical parameters of the numerical
integration are chosen in the following way:
\begin{itemize}
\item The whole integration process took \textbf{6m34s} on an
Intel\textregistered Core\texttrademark i7-4610M CPU @ 3.00GHz x 4,
compiled with gcc version 4.9.2 (64-bit Ubuntu 4.9.2-10ubuntu13)
and with flags -O2 -frounding-math.
\item In contrast, if instead of the FFT-algorithm we used the
direct algorithm, then the whole integration process took \textbf{17m29s}.
\end{itemize}
%
%

The infinite-dimensional bounds at time $T=4.62$ after performing $2310$ numerical integration steps
are
\[
\begin{matrix}
\footnotesize{
\begin{array}{|c|c|}\hline\mathbf{k} & \mathbf{a_k} \\ \hline\hline
1 & -1.13\cdot 10^{-7}+[-5.199,5.199]10^{-8}\\ 
\mathbf{2} & \mathbf{0.3484+[-7.806,7.806]10^{-9}}\\ 
3 & -3.735\cdot 10^{-7}+[-5.197,5.197]10^{-8}\\ 
4 & 2.276\cdot 10^{-7}+[-2.347,2.347]10^{-9}\\ 
5 & 4.499\cdot 10^{-8}+[-1.555,1.555]10^{-8}\\ 
6 & -0.02108+[-1.812,1.812]10^{-9}\\ 
7 & -1\cdot 10^{-7}+[-7.195,7.195]10^{-9}\\ 
8 & -2.038\cdot 10^{-8}+[-4.573,4.573]10^{-10}\\ 
9 & 4.559\cdot 10^{-8}+[-1.703,1.703]10^{-9}\\ 
10 & 1.207\cdot 10^{-3}+[-2.437,2.437]10^{-10}\\ 
11 & -1.722\cdot 10^{-8}+[-7.148,7.148]10^{-10}\\ 
12 & -2.243\cdot 10^{-8}+[-2.337,2.337]10^{-10}\\ 
13 & 1.96\cdot 10^{-8}+[-1.743,1.743]10^{-10}\\ 
14 & -7.098\cdot 10^{-5}+[-6.009,6.009]10^{-11}\\ 
15 & -1.855\cdot 10^{-10}+[-2.586,2.586]10^{-10}\\ 
16 & 8.757\cdot 10^{-11}+[-3.642,3.642]10^{-10}\\ 
17 & -9.373\cdot 10^{-11}+[-3.223,3.223]10^{-10}\\ 
18 & 4.195\cdot 10^{-6}+[-1.181,1.181]10^{-9}\\ 
19 & -4.39\cdot 10^{-12}+[-3.988,3.988]10^{-11}\\ 
20 & -1.142\cdot 10^{-11}+[-3.233,3.233]10^{-11}\\ 
21 & 1.051\cdot 10^{-11}+[-2.777,2.777]10^{-11}\\ 
22 & -2.478\cdot 10^{-7}+[-1.064,1.064]10^{-10}\\ 
23 & 3.143\cdot 10^{-13}+[-2.981,2.981]10^{-12}\\ 
24 & 8.888\cdot 10^{-13}+[-2.403,2.403]10^{-12}\\ \hline
\end{array} } &
\footnotesize{ 
\begin{array}{|c|c|}
\hline
25 & -8.013\cdot 10^{-13}+[-2.045,2.045]10^{-12}\\ 
26 & 1.464\cdot 10^{-8}+[-7.9,7.9]10^{-12}\\ 
27 & -2.201\cdot 10^{-14}+[-2.113,2.113]10^{-13}\\ 
28 & -6.322\cdot 10^{-14}+[-1.699,1.699]10^{-13}\\ 
29 & 5.65\cdot 10^{-14}+[-1.436,1.436]10^{-13}\\ 
30 & -8.65\cdot 10^{-10}+[-5.562,5.562]10^{-13}\\ 
31 & 1.5\cdot 10^{-15}+[-1.458,1.458]10^{-14}\\ 
32 & 4.326\cdot 10^{-15}+[-1.168,1.168]10^{-14}\\ 
33 & -3.846\cdot 10^{-15}+[-9.814,9.814]10^{-15}\\ 
34 & 5.111\cdot 10^{-11}+[-3.809,3.809]10^{-14}\\ 
35 & -1.003\cdot 10^{-16}+[-9.878,9.878]10^{-16}\\ 
36 & -2.894\cdot 10^{-16}+[-7.888,7.888]10^{-16}\\ 
37 & 2.563\cdot 10^{-16}+[-6.592,6.592]10^{-16}\\ 
38 & -3.02\cdot 10^{-12}+[-2.562,2.562]10^{-15}\\ 
39 & 6.611\cdot 10^{-18}+[-6.603,6.603]10^{-17}\\ 
40 & 1.906\cdot 10^{-17}+[-5.255,5.255]10^{-17}\\ 
41 & -1.684\cdot 10^{-17}+[-4.373,4.373]10^{-17}\\ 
42 & 1.785\cdot 10^{-13}+[-1.702,1.702]10^{-16}\\ 
43 & -4.31\cdot 10^{-19}+[-4.366,4.366]10^{-18}\\ 
44 & -1.241\cdot 10^{-18}+[-3.464,3.464]10^{-18}\\ 
45 & 1.094\cdot 10^{-18}+[-2.872,2.872]10^{-18}\\ 
46 & -1.055\cdot 10^{-14}+[-1.118,1.118]10^{-17}\\ 
47 & 2.785\cdot 10^{-20}+[-2.861,2.861]10^{-19}\\ 
48 & 8.009\cdot 10^{-20}+[-2.264,2.264]10^{-19}\\ \hline
49-250 &\text{small intervals of width }10^{-19}\\ \hline
\geq 251 & <3.96\cdot 10^{-64}/k^{6}\\\hline\end{array} }
\end{matrix}
\]

The bounds above are within $\wsg$, i.e., the basin of attraction of the stable fixed point,
which means that the third step of the proof of the main theorem from Section~\ref{secproof} succeeded.

\footnotesize
\bibliography{CW_dbcp_hetero,wanner1a,wanner1b,wanner2a,wanner2b}

\begin{thebibliography}{10}

\bibitem{bahiana:oono:90a}
M.~Bahiana and Y.~Oono.
\newblock Cell dynamical system approach to block copolymers.
\newblock {\em Physical Review A}, 41:6763--6771, 1990.

\bibitem{BW2}
F.~A. Bartha and W.~Tucker.
\newblock Fixed points of a destabilized {Kuramoto–Sivashinsky} equation.
\newblock {\em Applied Mathematics and Computation}, 266:339 -- 349, 2015.

\bibitem{bates:fife:93a}
P.~W. Bates and P.~C. Fife.
\newblock The dynamics of nucleation for the {C}ahn-{H}illiard equation.
\newblock {\em SIAM Journal on Applied Mathematics}, 53(4):990--1008, 1993.

\bibitem{bloemker:etal:10a}
D.~Bl\"omker, B.~Gawron, and T.~Wanner.
\newblock Nucleation in the one-dimensional stochastic {C}ahn-{H}illiard model.
\newblock {\em Discrete and Continuous Dynamical Systems, Series A},
  27(1):25--52, 2010.

\bibitem{CFL}
X.~Cabr\'e, E.~Fontich, and R.~de~la Llave.~T.
\newblock parameterization method for invariant manifolds i: manifolds
  associated to non-resonant subspaces.
\newblock {\em Indiana Univ. Math. J.}, 52:283--328, 2003.

\bibitem{cahn:hilliard:58a}
J.~W. Cahn and J.~E. Hilliard.
\newblock Free energy of a nonuniform system {I}. {I}nterfacial free energy.
\newblock {\em Journal of Chemical Physics}, 28:258--267, 1958.

\bibitem{CGL}
R.~Castelli, M.~Gameiro, and J.-P. Lessard.
\newblock Rigorous numerics for ill-posed {PDEs}: periodic orbits in the
  boussinesq equation.
\newblock Preprint.

\bibitem{choksi:etal:11a}
R.~Choksi, M.~Maras, and J.~F. Williams.
\newblock 2{D} phase diagram for minimizers of a {C}ahn-{H}illiard functional
  with long-range interactions.
\newblock {\em SIAM Journal on Applied Dynamical Systems}, 10(4):1344--1362,
  2011.

\bibitem{choksi:peletier:09a}
R.~Choksi, M.~A. Peletier, and J.~F. Williams.
\newblock On the phase diagram for microphase separation of diblock copolymers:
  {A}n approach via a nonlocal {C}ahn-{H}illiard functional.
\newblock {\em SIAM Journal on Applied Mathematics}, 69(6):1712--1738, 2009.

\bibitem{choksi:ren:03a}
R.~Choksi and X.~Ren.
\newblock On the derivation of a density functional theory for microphase
  separation of diblock copolymers.
\newblock {\em Journal of Statistical Physics}, 113:151--76, 2003.

\bibitem{Cy}
J.~Cyranka.
\newblock Efficient and generic algorithm for rigorous integration forward in
  time of {dPDEs}: {Part I}.
\newblock {\em Journal of Scientific Computing}, 59(1):28--52, 2014.

\bibitem{C1}
J.~Cyranka.
\newblock Existence of globally attracting fixed points of viscous {Burgers}
  equation with constant forcing. a computer assisted proof.
\newblock {\em Topological Methods in Nonlinear Analysis}, 45(2):655–697,
  2015.

\bibitem{codes}
J.~Cyranka.
\newblock Numerical data and source codes from the proof.
\newblock Technical report, 2017.
\newblock
  \href{http://bitbucket.org/dzako/dbcp\_proof}{http://bitbucket.org/dzako/dbcp\_proof}.

\bibitem{CZ}
J.~Cyranka and P.~Zgliczyński.
\newblock Existence of globally attracting solutions for one-dimensional
  viscous {Burgers} equation with nonautonomous forcing---a computer assisted
  proof.
\newblock {\em SIAM Journal on Applied Dynamical Systems}, 14(2):787--821,
  2015.

\bibitem{CzZ}
A.~Czechowski and P.~Zgliczyński.
\newblock Rigorous numerics for {PDEs} with indefinite tail: Existence of a
  periodic solution of the {Boussinesq} equation with time-dependent forcing.
\newblock {\em Schedae Informaticae}, 2015(Volume 24), 2016.

\bibitem{D}
G.~Dahlquist.
\newblock {\em Stability and Error Bounds in the Numerical Intgration of
  Ordinary Differential Equations}.
\newblock Transactions of the Royal Institute of Technology. Almqvist \&
  Wiksells, Uppsala, 1958.

\bibitem{DHMO}
S.~Day, Y.~Hiraoka, K.~Mischaikow, and T.~Ogawa.
\newblock Rigorous numerics for global dynamics: A study of the
  {Swift--Hohenberg} equation.
\newblock {\em SIAM Journal on Applied Dynamical Systems}, 4(1):1--31, 2005.

\bibitem{desi:etal:11a}
J.~P. Desi, H.~Edrees, J.~Price, E.~Sander, and T.~Wanner.
\newblock The dynamics of nucleation in stochastic {C}ahn-{M}orral systems.
\newblock {\em SIAM Journal on Applied Dynamical Systems}, 10(2):707--743,
  2011.

\bibitem{fife:kielhofer:etal:97a}
P.~C. Fife, H.~Kielh\"ofer, S.~Maier-Paape, and T.~Wanner.
\newblock Perturbation of doubly periodic solution branches with applications
  to the {C}ahn-{H}illiard equation.
\newblock {\em Physica D}, 100(3--4):257--278, 1997.

\bibitem{GLP}
M.~Gameiro, J.-P. Lessard, and A.~Pugliese.
\newblock Computation of smooth manifolds via rigorous multi-parameter
  continuation in infinite dimensions.
\newblock {\em Found. Comput. Math.}, 16(2):531--575, 2016.

\bibitem{grinfeld:novickcohen:95a}
M.~Grinfeld and A.~Novick-Cohen.
\newblock Counting stationary solutions of the {C}ahn-{H}illiard equation by
  transversality arguments.
\newblock {\em Proceedings of the Royal Society of Edinburgh}, 125A:351--370,
  1995.

\bibitem{grinfeld:novickcohen:99a}
M.~Grinfeld and A.~Novick-Cohen.
\newblock The viscous {C}ahn-{H}illiard equation: {M}orse decomposition and
  structure of the global attractor.
\newblock {\em Transactions of the American Mathematical Society},
  351(6):2375--2406, 1999.

\bibitem{MJM}
J.~D.~M. James and K.~Mischaikow.
\newblock Rigorous a posteriori computation of (un)stable manifolds and
  connecting orbits for analytic maps.
\newblock {\em SIAM Journal on Applied Dynamical Systems}, 12(2):957--1006,
  2013.

\bibitem{johnson:etal:13a}
I.~Johnson, E.~Sander, and T.~Wanner.
\newblock Branch interactions and long-term dynamics for the diblock copolymer
  model in one dimension.
\newblock {\em Discrete and Continuous Dynamical Systems. Series A},
  33(8):3671--3705, 2013.

\bibitem{KZ}
T.~Kapela and P.~Zgliczyński.
\newblock A {Lohner}-type algorithm for control systems and ordinary
  differential inclusions.
\newblock {\em Discrete and Continuous Dynamical Systems - Series B},
  11(2):365--385, 2009.

\bibitem{Lo}
R.~Lohner.
\newblock {\em Computation of guaranteed enclosures for the solutions of
  ordinary initial and boundary value problems}, chapter Computational Ordinary
  Differential Equations.
\newblock Clarendon Press, Oxford, 1992.

\bibitem{L}
S.~M. Lozinskii.
\newblock Error estimates for the numerical integration of ordinary
  differential equations, part i.
\newblock {\em Izv. Vyss. Uceb. Zaved. Matematica}, 6:52--90, 1958.

\bibitem{maier:etal:08a}
S.~Maier-Paape, U.~Miller, K.~Mischaikow, and T.~Wanner.
\newblock Rigorous numerics for the {C}ahn-{H}illiard equation on the unit
  square.
\newblock {\em Revista Matematica Complutense}, 21(2):351--426, 2008.

\bibitem{maier:etal:07a}
S.~Maier-Paape, K.~Mischaikow, and T.~Wanner.
\newblock Structure of the attractor of the {C}ahn-{H}illiard equation on a
  square.
\newblock {\em International Journal of Bifurcation and Chaos},
  17(4):1221--1263, 2007.

\bibitem{maier:wanner:97a}
S.~Maier-Paape and T.~Wanner.
\newblock Solutions of nonlinear planar elliptic problems with triangle
  symmetry.
\newblock {\em Journal of Differential Equations}, 136(1):1--34, 1997.

\bibitem{maier:wanner:98a}
S.~Maier-Paape and T.~Wanner.
\newblock Spinodal decomposition for the {C}ahn-{H}illiard equation in higher
  dimensions. {P}art {I}: {P}robability and wavelength estimate.
\newblock {\em Communications in Mathematical Physics}, 195(2):435--464, 1998.

\bibitem{maier:wanner:00a}
S.~Maier-Paape and T.~Wanner.
\newblock Spinodal decomposition for the {C}ahn-{H}illiard equation in higher
  dimensions: {N}onlinear dynamics.
\newblock {\em Archive for Rational Mechanics and Analysis}, 151(3):187--219,
  2000.

\bibitem{mccord:88b}
C.~McCord.
\newblock Mappings and homological properties in the {C}onley index theory.
\newblock {\em Ergodic Theory and Dynamical Systems}, 8$^*$(Charles Conley
  Memorial Issue):175--198, 1988.

\bibitem{MM}
K.~Mischaikow and M.~Mrozek.
\newblock {\em Conley Index}, volume~2 of {\em Handbook of Dynamical Systems},
  chapter~9, pages 393--460.
\newblock Elsevier, 2002.

\bibitem{Mo}
R.~Moore.
\newblock {\em Interval Analysis}.
\newblock Prentice Hall, Englewood Cliffs, N.J., 1966.

\bibitem{Ma}
Z.~Mączyńska.
\newblock {\em Rigorous numerics for dissipative partial differential
  equations}.
\newblock PhD thesis, Jagiellonian University, 2011, 2011.

\bibitem{nishiura:ohnishi:95a}
Y.~Nishiura and I.~Ohnishi.
\newblock Some mathematical aspects of the micro-phase separation in diblock
  copolymers.
\newblock {\em Physica D}, 84(1-2):31--39, 1995.

\bibitem{ohta:kawasaki:86a}
T.~Ohta and K.~Kawasaki.
\newblock Equilibrium morphology of block copolymer melts.
\newblock {\em Macromolecules}, 19:2621--2632, 1986.

\bibitem{RMJ}
C.~Reinhardt and J.~Mireles~James.
\newblock {Fourier-Taylor} parameterization of unstable manifolds for parabolic
  partial differential equations: Formalism, implementation and rigorous
  validation.
\newblock preprint arXiv:1601.00307v2, 2016.

\bibitem{ren:wei:06b}
X.~Ren and J.~Wei.
\newblock Droplet solutions in the diblock copolymer problem with skewed
  monomer composition.
\newblock {\em Calculus of Variations and Partial Differential Equations},
  25(3):333--359, 2006.

\bibitem{ren:wei:06a}
X.~Ren and J.~Wei.
\newblock Existence and stability of spherically layered solutions of the
  diblock copolymer equation.
\newblock {\em SIAM Journal on Applied Mathematics}, 66(3):1080--1099, 2006.

\bibitem{ren:wei:07b}
X.~Ren and J.~Wei.
\newblock Many droplet pattern in the cylindrical phase of diblock copolymer
  morphology.
\newblock {\em Reviews in Mathematical Physics}, 19(8):879--921, 2007.

\bibitem{ren:wei:07a}
X.~Ren and J.~Wei.
\newblock Single droplet pattern in the cylindrical phase of diblock copolymer
  morphology.
\newblock {\em Journal of Nonlinear Science}, 17(5):471--503, 2007.

\bibitem{sander:wanner:99a}
E.~Sander and T.~Wanner.
\newblock Monte {C}arlo simulations for spinodal decomposition.
\newblock {\em Journal of Statistical Physics}, 95(5--6):925--948, 1999.

\bibitem{sander:wanner:00a}
E.~Sander and T.~Wanner.
\newblock Unexpectedly linear behavior for the {C}ahn-{H}illiard equation.
\newblock {\em SIAM Journal on Applied Mathematics}, 60(6):2182--2202, 2000.

\bibitem{sander:wanner:16a}
E.~Sander and T.~Wanner.
\newblock Validated saddle-node bifurcations and applications to lattice
  dynamical systems.
\newblock {\em SIAM Journal on Applied Dynamical Systems}, 15(3):1690--1733,
  2016.

\bibitem{temam:88a}
R.~Temam.
\newblock {\em Infinite-Dimensional Dynamical Systems in Mechanics and
  Physics}.
\newblock Springer-Verlag, New York -- Berlin -- Heidelberg, 1988.

\bibitem{T}
W.~Tucker.
\newblock {\em Validated numerics}.
\newblock Princeton University Press, Princeton, NJ, 2011.

\bibitem{BBMJLM}
J.~B. van~den Berg, J.~D. Mireles-James, J.-P. Lessard, and K.~Mischaikow.
\newblock Rigorous numerics for symmetric connecting orbits: Even homoclinics
  of the {Gray–Scott} equation.
\newblock {\em SIAM Journal on Mathematical Analysis}, 43(4):1557--1594, 2011.

\bibitem{vandenBerg2016}
J.~B. van~den Berg, J.~D. Mireles~James, and C.~Reinhardt.
\newblock Computing (un)stable manifolds with validated error bounds:
  Non-resonant and resonant spectra.
\newblock {\em Journal of Nonlinear Science}, 26(4):1055--1095, 2016.

\bibitem{BW}
J.~B. van~den Berg and J.~Williams.
\newblock Validation of the bifurcation diagram in the 2d ohta-kawasaki
  problem.
\newblock preprint, 2016.

\bibitem{wanner:04a}
T.~Wanner.
\newblock Maximum norms of random sums and transient pattern formation.
\newblock {\em Transactions of the American Mathematical Society},
  356(6):2251--2279, 2004.

\bibitem{wanner:16a}
T.~Wanner.
\newblock Topological analysis of the diblock copolymer equation.
\newblock In Y.~Nishiura and M.~Kotani, editors, {\em Mathematical Challenges
  in a New Phase of Materials Science}, volume 166 of {\em Springer Proceedings
  in Mathematics \& Statistics}, pages 27--51. Springer-Verlag, 2016.

\bibitem{wanner:17b}
T.~Wanner.
\newblock Computer-assisted bifurcation diagram validation and applications in
  materials science.
\newblock {\em Proceedings of Symposia in Applied Mathematics}, 74:to appear,
  2017.

\bibitem{wanner:17a}
T.~Wanner.
\newblock Computer-assisted equilibrium validation for the diblock copolymer
  model.
\newblock {\em Discrete and Continuous Dynamical Systems, Series A},
  37(2):1075--1107, 2017.

\bibitem{Wilczak2006}
D.~Wilczak.
\newblock The existence of {Shilnikov} homoclinic orbits in the {Michelson}
  system: A computer assisted proof.
\newblock {\em Foundations of Computational Mathematics}, 6(4):495--535, 2006.

\bibitem{Wilczak2005}
D.~Wilczak and P.~Zgliczy{\'{n}}ski.
\newblock Heteroclinic connections between periodic orbits in planar restricted
  circular three body problem. {Part} ii.
\newblock {\em Communications in Mathematical Physics}, 259(3):561--576, 2005.

\bibitem{WZ}
D.~Wilczak and P.~Zgliczyński.
\newblock Topological method for symmetric periodic orbits for maps with a
  reversing symmetry.
\newblock {\em Discrete and Continuous Dynamical Systems}, 17(3):629--652,
  2007.

\bibitem{Zattr}
P.~Zgliczynski.
\newblock Attracting fixed points for the {Kuramoto--Sivashinsky} equation: A
  computer assisted proof.
\newblock {\em SIAM Journal on Applied Dynamical Systems}, 1(2):215--235, 2002.

\bibitem{ZPer}
P.~Zgliczynski.
\newblock Rigorous numerics for dissipative partial differential equationsii.
  periodic orbit for the {Kuramoto--Sivashinsky} pde---a computer-assisted
  proof.
\newblock {\em Foundations of Computational Mathematics}, 4(2):157--185, 2004.

\bibitem{ZM}
P.~Zgliczynski and K.~Mischaikow.
\newblock Rigorous numerics for partial differential equations: The
  {Kuramoto---Sivashinsky} equation.
\newblock {\em Foundations of Computational Mathematics}, 1(3):255--288, 2001.

\bibitem{Zcc}
P.~Zgliczyński.
\newblock Covering relations, cone conditions and the stable manifold theorem.
\newblock {\em Journal of Differential Equations}, 246(5):1774 -- 1819, 2009.

\bibitem{Z3}
P.~Zgliczyński.
\newblock Rigorous numerics for dissipative {PDEs} {III}. an effective
  algorithm for rigorous integration of dissipative {PDEs}.
\newblock {\em Topological Methods in Nonlinear Analysis}, 36:197--262, 2010.

\bibitem{Zhet}
P.~Zgliczyński.
\newblock Heteroclinic connection in {Kuramoto-Sivashinsky PDE} - a computer
  assisted proof.
\newblock work in progress, 2017.

\bibitem{ZG}
P.~Zgliczyński and M.~Gidea.
\newblock Covering relations for multidimensional dynamical systems.
\newblock {\em Journal of Differential Equations}, 202(1):32 -- 58, 2004.

\bibitem{CM}
A.~Ćwiszewski and M.~Maciejowski.
\newblock Stationary solutions and connecting orbits for $p$-laplace equation.
\newblock preprint, 2016.

\end{thebibliography}
\bibliographystyle{abbrv}

\end{document}